\documentclass[12pt,letterpaper]{article}

\usepackage[utf8]{inputenc}
\usepackage{amsmath}
\usepackage{amsfonts}
\usepackage{amssymb}
\usepackage{hyperref}
\usepackage[left=2cm,right=2cm,top=2cm,bottom=2cm]{geometry}
\usepackage{graphicx}
\usepackage{tikz}
\usepackage{tikz-cd}

\usepackage[ backend=biber, sorting= nyt, style=alphabetic, url = false, isbn = false, doi = false, minalphanames=4, maxalphanames=10, maxbibnames=10]{biblatex}
\DeclareMathOperator{\Aut}{Aut} 
\DeclareMathOperator{\res}{Res} 
\DeclareMathOperator{\End}{End} 
\DeclareMathOperator{\ad}{ad} 
\DeclareMathOperator{\id}{id} 
\DeclareMathOperator{\wt}{wt} 
\DeclareMathOperator{\Span}{Span} 

\DeclareMathOperator{\GL}{GL} 
\newcommand{\one}{\mathbf{1}} 
\newcommand{\iu}{\boldsymbol{i}} 
\newcommand{\lvl}{\kappa} 
\newcommand{\el}{\ell} 
\newcommand{\type}[3]{\left( \begin{smallmatrix} #3\\ #1 \,  #2 \end{smallmatrix}\right)} 
\newcommand{\gtype}{$\left(\!\begin{smallmatrix} g\\ 1 \  g \end{smallmatrix}\!\right)$-type} 
\newcommand{\D}{\mathbb{D}} 
\newcommand{\N}{N} 
\newcommand{\g}{\widetilde{g}} 
\newcommand{\W}{\widetilde{W}} 
\newcommand{\n}{n} 
\newcommand{\cln}{\colon} 
\DeclareMathOperator{\Frac}{Frac} 

\usepackage{amsthm}
\newtheorem{theorem}{Theorem}[section]
\newtheorem{proposition}[theorem]{Proposition}
\newtheorem{corollary}[theorem]{Corollary}
\newtheorem{lemma}[theorem]{Lemma}

\theoremstyle{definition}
\newtheorem{definition}[theorem]{Definition}
\newtheorem{remark}[theorem]{Remark}
\newtheorem{conjecture}[theorem]{Conjecture}
\newtheorem{example}[theorem]{Example}

\usetikzlibrary{decorations.markings,arrows}

\pgfdeclarelayer{edgelayer}
\pgfdeclarelayer{nodelayer}
\pgfsetlayers{edgelayer,nodelayer,main}
\tikzstyle{none}=[]
\tikzstyle{plainedge}=[-, draw=black=]
\tikzstyle{dotted line}=[dashed, draw=black=]

\author{Daniel Tan}
\date{}
\title{\bf 
Differential equations for twisted intertwining operators in abelian orbifold conformal field theory}

\numberwithin{equation}{section}

\begin{document}

\maketitle

\begin{abstract}

Let $ V $ be a vertex operator algebra with two commuting automorphisms $ g_1 $ and $ g_2 $ (possibly of infinite order and  non-semisimple). We derive a Jacobi identity for a twisted intertwining operator $ \mathcal{Y} $ of type $ \type{W_1}{W_2}{W_3} $ when $ W_1 $ is a $ g_1 $-twisted $ V $-module, $ W_2 $ is a $ g_2 $-twisted $ V $-module, and $ W_3 $ is a $ g_1 g_2 $-twisted $ V $-module. 

For any abelian group $ G $ of automorphisms of a CFT-type vertex operator algebra, under certain Noetherian assumptions, we derive differential equations satisfied by the multi-series $ \langle w_0, \mathcal{Y}_1(w_1,z_1) \cdots \mathcal{Y}_\N(w_\N,z_\N) w_{\N+1} \rangle $ when $ \mathcal{Y}_j $ are intertwining operators among $ C_2 $-cofinite and discretely graded $ V $-modules twisted by elements in $ G $. For modules twisted by automorphisms with only real positive eigenvalues, it suffices to assume $ C_1 $-cofiniteness and discrete grading. When all but at most two external modules are twisted by automorphisms having only positive real eigenvalues, it suffices to assume all in-state and internal modules are $ C_1 $-cofinite and discretely graded. 

When $ G $ is a finite abelian group, we use the theory of regular singular points to prove that the multi-series $ \langle w_0, \mathcal{Y}_1(w_1,z_1) \cdots \mathcal{Y}_\N(w_\N,z_\N) w_{\N+1} \rangle $ converges absolutely to a multivalued analytic function when $ |z_1| > \cdots > |z_\N| > 0 $ and analytically extends to the region $ z_i \neq 0 $, $ z_j - z_k \neq 0 $, for $ 1 \leq i \leq \N $ and $ 1 \leq j < k \leq \N $. The vertex operator algebra $ V $ is not assumed to be $ C_2 $-cofinite or simple.
\end{abstract}

\setcounter{section}{-1}

\section{Introduction}
The vertex-operator-algebraic approach to orbifold conformal field theory studies twisted modules for a given vertex operator algebra, and twisted intertwining operators among them. In mathematics, certain elements and structures of twisted modules were first discovered in works such as \cite{LepowskyWilson,FLM84,FLM_E8,LepowskyCalculus,FLM_VOCalc,Lepowsky_VOandMonster}. In physics, orbifold conformal field theory was initiated in \cite{DHVW1,DHVW2} inspired by some of these mathematical insights, namely the construction of the Moonshine Module \cite{FLM84}. Twisted modules are the mathematical counterparts of certain additional sectors added to a conformal field theory before performing a kind of quotient by a group of its symmetries. The Moonshine Module \cite{FLM} came to be interpreted as the first rigorous example of a conformal field theory constructed with this orbifold-conformal-field-theoretic quotient.

Even though the Moonshine Module was one of the first vertex operator algebras constructed, and the first rigorous orbifold conformal field theory, the general theory for the orbifold construction remains underdeveloped. Notably, the convergence of matrix coefficients of products of twisted intertwining operators (these are the mathematical counterparts to the most general chiral genus-zero correlation functions appearing in orbifold conformal field theory) has not been proven in many cases. (See the paper \cite{Reptheoryandorbifoldcft} and the introduction to \cite{TanMostlyUntwisted} for discussions of greater depth about orbifold conformal field theory, twisted representation theory of vertex operator algebras, and the necessity for proving the convergence of products of twisted intertwining operators.)

In this paper, we offer three contributions to the study of the convergence of matrix coefficients of products of twisted intertwining operators among modules twisted by commuting automorphisms of a vertex operator algebra $ V $. First, in Proposition~\ref{prop:Jacobi c5}, we derive a Jacobi identity for a twisted intertwining operator $ \mathcal{Y} $ of type $ \type{W_1}{W_2}{W_3} $ when $ W_1 $ is a $ g_1 $-twisted $ V $-module, $ W_2 $ is a $ g_2 $-twisted $ V $-module, and $ W_3 $ is a $ g_1 g_2 $-twisted $ V $-module. We also show in Proposition~\ref{prop:duality equivalence} that the Jacobi identity can replace the analytic duality axiom in the definition of twisted intertwining operators among modules twisted by commuting automorphisms. Thus, we provide a completely algebraic definition for these twisted intertwining operators. 

Second, in Theorem~\ref{thm: DEs c5} we derive certain differential equations with controlled coefficients satisfied by the multi-series $ \langle w_0, \mathcal{Y}_1(w_1,z_1) \cdots \mathcal{Y}_\N(w_\N,z_\N) w_{\N+1} \rangle $, where $ \mathcal{Y}_j $ are intertwining operators among modules twisted by elements in a general (possibly infinite) abelian subgroup $ G $ of $ \Aut V $. This result requires the vertex operator algebra $ V $ to be of CFT type (i.e., positive energy), the modules to be $ C_2 $-cofinite and discretely graded (i.e., quasi-finite dimensional), and the coefficient ring $ R $ defined in Section~\ref{sec:setting c5} to be Noetherian. For modules twisted by automorphisms with only real positive eigenvalues, it suffices to assume $ C_1 $-cofiniteness and discrete grading. Furthermore, there is a case where all but at most two external modules are twisted by automorphisms having only positive real eigenvalues, for which all in-state and intermediate modules can be assumed to be $ C_1 $-cofinite and discretely graded. Using Huang’s method from \cite{HuangC1}, we do not need to impose cofiniteness conditions for the dual $ W_0 $ of the out-state module. In Theorem~\ref{thm: DEs c c5}, we use further Noetherian assumptions to improve the control on the coefficients of the differential equations.

Third, when $ G $ is a finite abelian group, we use these differential equations in conjunction with the theory of regular singular points to prove Theorem~\ref{thm: convergence c5} which states that the multi-series $ \langle w_0, \mathcal{Y}_1(w_1,z_1) \cdots \mathcal{Y}_\N(w_\N,z_\N) w_{\N+1} \rangle $ converges absolutely to a multivalued analytic function when $ |z_1| > \cdots > |z_\N| > 0 $ and analytically extends to the region $ z_i \neq 0 $, $ z_j - z_k \neq 0 $, for $ 1 \leq i \leq \N $ and $ 1 \leq j < k \leq \N $. When $ \N = 2 $, Proposition~\ref{prop:sing special c5} establishes the regularity of the component-isolated singularities of this function that are obtained by invertible linear changes of coordinates. 

Our convergence results require neither $ C_2 $-cofiniteness nor simplicity of $ V $. We also do not need to establish $ C_2 $-cofiniteness of the fixed-point subalgebra $ V^G $. This should be compared to an argument that uses Miyamoto's fixed-point subalgebra theory \cite{Miyamoto,CM} with Huang's untwisted differential equations theory \cite{HuangDEs, HuangC1}. Section~\ref{sec:examples} provides examples of $ C_1 $-cofinite and $ C_2 $-cofinite twisted modules for which our results can be applied. 

Even though we prove the convergence of products of twisted intertwining operators only when $ G $ is finite, the differential equations that we derive here provide a possible pathway for proving convergence in cases involving automorphisms of infinite order.

The results of this paper are a continuation and extension of some of the results from \cite{TanMostlyUntwisted}; we now provide insight into how the results from \cite{TanMostlyUntwisted} have hinted at and led to the discovery of the present results. In \cite{TanMostlyUntwisted}, the ``twisted KZ equations" of \cite{deBoer:2001nw} were derived in a mathematically rigorous setting and extended to general automorphisms. It is important to note that the twisted KZ equations are for \gtype\ intertwining operators for an automorphism $ g $ of an affine vertex operator algebra, i.e., twisted intertwining operators among an untwisted module and a $ g $-twisted module, with the third module necessarily also being $ g $-twisted. Physically speaking, one may view a matrix coefficient of a product of \gtype\ intertwining operators as a chiral genus-zero correlation function with a $ g $-twisted sector at $ 0 $, a $ g^{-1} $-twisted sector at $ \infty $, and untwisted sectors at all other field insertions. 

In the untwisted story, the KZ equations hint that regular singular differential equations can be derived for products of intertwining operators among $ C_1 $-cofinite discretely graded $ V $-modules for a general vertex operator algebra $ V $ (see \cite{HuangDEs}). Parallel to the untwisted story, the twisted KZ equations hint that differential equations can be derived for products of \gtype\ intertwining operators among $ C_1 $-cofinite discretely graded untwisted and $ g $-twisted $ V $-modules for a general vertex operator algebra $ V $ with automorphism $ g $, as was shown in \cite{TanMostlyUntwisted}. 

As is immediately clear from an explicit inspection of the twisted KZ equations, the coefficients in the differential equations may be multivalued with logarithmic singularities at the points that we would like to be regular singular. In the case that $ g $ has finite order, the logarithmic singularities disappear and the multivaluedness can be removed via a change of variables. With this insight, one can check in the \gtype\ intertwining operator case for a general vertex operator algebra that regular singularities appear when $ g $ has finite order, and hence the convergence and analytic continuation of products of \gtype\ intertwining operators can be proven for a general vertex operator algebra. 

Physical orbifold conformal field theory also suggests that the field insertions of this chiral genus-zero correlation function can move around and change their radial order. For us, this suggests the convergence of the multi-series $ \langle w_0, \mathcal{Y}_1(w_1,x_1) \cdots \mathcal{Y}_\N(w_\N,x_\N) w_{\N+1} \rangle $ defined from products of twisted intertwining operators with $ w_0 \in W_0, \dots, w_{\N+1}\in W_{\N+1} $ being elements of modules where one is $ g $-twisted, another is $ g^{-1} $-twisted, and the others are untwisted---in every possible permutation. We refer to this as the \textit{mostly untwisted case}. Indeed, we are able to prove such convergence when $ g $ has finite order in Section~\ref{sec:convergence}.

To derive the necessary differential equations for the mostly untwisted case, a \textit{new} Jacobi identity is needed. Unlike in \cite{TanMostlyUntwisted}, the first input of the intertwining operator can no longer be assumed to be untwisted, so there are additional complications caused by branches. When dealing with this complication, it is natural to simultaneously consider the case where the modules are twisted by general elements in a possibly infinite abelian subgroup $ G $ of $ \Aut V $. In this case, the commutator formula for $ [Y(u,x_1), \mathcal{Y}(w,x_2)] $ used in \cite{TanMostlyUntwisted} no longer applies. So, a new method for deriving the differential equations is used in Section~\ref{sec:computations}.

However, to use this new method with a general abelian group $ G $, the modules are now assumed to be \textit{$ C_2 $-cofinite} and discretely graded in general (compare this with the $ C_1 $-cofiniteness needed in \cite{TanMostlyUntwisted} and \cite{HuangDEs, HuangC1}). It suffices, however, to assume $ C_1 $-cofiniteness and discrete grading in certain cases involving automorphisms that have only real positive eigenvalues (see Section~\ref{sec:differential equations} for details), such as the mostly untwisted case. We then impose that $ G $ is finite, so that logarithmic terms disappear from the coefficients in the differential equations and suitable changes of variables yield regular singular points. With these regular singular differential equations, we prove convergence in the mostly untwisted case, as well as the case where the modules are twisted by general elements in $ G $.

While the results of this paper extend some of the mostly untwisted results in \cite{TanMostlyUntwisted} and go much further into the abelian group case, we were unable to extend the twisted KZ equations (using the new techniques developed here) to twisted intertwining operators among modules twisted by arbitrary elements of an abelian group---not even when the group is finite. We suspect this is because the twisted modules for an affine vertex operator algebra are not $ C_2 $-cofinite in general. 

We hope that the results of this paper can be applied with the results of \cite{DuHuang} to construct $ G $-crossed tensor categories when $ G $ is finite abelian, similar to the untwisted case in the series \cite{HLZ1,HLZ2,HLZ3,HLZ4,HLZ5,HLZ6,HLZ7,HLZ8}. Compare this with McRae’s fixed-point-subalgebra approach \cite{McRaeTwisted} that starts with a suitable tensor category of $ V^G $-modules. The examples in Section~\ref{sec:examples} provide settings in which to investigate direct constructions of $ G $-crossed tensor categories from twisted $ V $-modules. Our results Theorem~\ref{thm: convergence c5} and Proposition~ \ref{prop:sing special c5}, which provide certain analytic hypotheses for constructions of $ P(z) $-tensor products among twisted modules in \cite{DuHuang}, do not require $ V $ to be $ C_2 $-cofinite or simple.

Even in the mostly untwisted case, for each fixed $ g \in \Aut V $ of finite order, we expect that the category of $ C_1 $-cofinite discretely graded $ g $-twisted modules will naturally have the structure of a bimodule category (see Chapter 7 of \cite{EGNO} for a definition) over the tensor category of $ C_1 $-cofinite discretely graded $ V $-modules. See the end of Section~\ref{sec:convergence} for a discussion, and see Section~\ref{sec:examples} for examples of non-$ C_2 $-cofinite vertex operator algebras with $ C_1 $-cofinite discretely graded twisted and untwisted modules. 

This paper is organized as follows. In Section~\ref{sec:notation}, we state notation and explain the formal logarithm expansion convention to be used in this paper. In Section~\ref{sec:Jacobi identity}, we derive a Jacobi identity for twisted intertwining operators among modules twisted by commuting automorphisms and prove its equivalence with the analytic duality axiom appearing in the definition for twisted intertwining operators. In Section~\ref{sec:setting c5}, we state the data needed and conditions satisfied by the products of twisted intertwining operators among modules twisted by elements of a possibly infinite abelian group that will be considered in the following sections. In Section~\ref{sec:computations}, we first present key identities by heuristic means, and then we rigorously prove these identities using our new Jacobi identity. In Section~\ref{sec:differential equations}, we use the identities from the previous section to derive differential equations for such  products of twisted intertwining operators under certain cofiniteness and Noetherian hypotheses. In Section~\ref{sec:improve coefficients}, we improve upon the form of the coefficients in the differential equations, under further Noetherian hypotheses. In Section~\ref{sec:convergence}, we assume the abelian group is finite and we use the differential equations to prove that such products of twisted intertwining operators converge absolutely in a certain region and analytically extend to a certain larger domain. Furthermore, we prove that certain singularities of such products of two twisted intertwining operators have the form of a regular singularity. In Section~\ref{sec:examples}, we provide examples of $ C_1 $-cofinite and $ C_2 $-cofinite twisted modules for $ C_2 $-cofinite and non-$ C_2 $-cofinite vertex operator algebras of CFT type, demonstrating the capability of our results. 

\paragraph{Acknowledgments} I would like to thank Professor Yi-Zhi Huang and Professor James Lepowsky for their helpful discussions and guidance while I developed these ideas during my PhD.

\section{Notation and conventions}
\label{sec:notation}
We use the notion of a vertex operator algebra from Definition~3.1.22 of \cite{LepowskyLi}. The base field for anything linear will always be $ \mathbb{C} $, which will be important for the theory. We refer the reader to \cite{HLZ2} for the framework of logarithmic formal calculus. To easily distinguish between formal variables (all formal variables commute) and complex numbers, we will use $ x, \log x,  x_i, \log x_i , \xi, \log \xi , \xi_i, \log \xi_i, y_i, \log y_i $ for formal variables and use $ z, z_i, \zeta, \zeta_i, \eta, \eta_i $ for complex numbers. For a vector space $ W $, we denote the space of formal generalized power series $ \sum_{n_1, \dots, n_\N} w_{n_1, \dots, n_\N} x_1^{n_1} \cdots x_\N^{n_\N} $ with coefficients in $ W $ and powers in $ \mathbb{C} $ by $ W\{x_1, \dots, x_\N \} $. We suppress logarithm variables and write $ f(x_1, \dots, x_\N) $ for elements in $ \mathbb{C}\{x_1, \dots ,x_\N\} [[ \log x_1, \dots , \log x_\N ]] $. We use $ \frac{\partial}{\partial x_i} $ to denote the derivation defined by
\begin{align*}
\frac{\partial}{\partial x_i} x_j^n = n x_j^{n-1} \delta_{i,j} \quad \text{and} \quad \frac{\partial}{\partial x_i} \log x_j = x_j^{-1} \delta_{i,j}.
\end{align*}
We sometimes write $ \frac{d}{dx} $ when there is only one pair of variables $ x $ and $ \log x $ involved. 

Given an operator $ \mathcal{A} $ and a locally nilpotent operator $ \mathcal{N} $ on $ W $, and formal variables $ x, x_1, \log x $, we have the following meanings for formal logarithms and exponentials:
\begin{equation}\label{eq:formal meanings}
\begin{gathered}
\log(1 + x) = \sum_{k = 1}^\infty \frac{(-1)^{k+1}}{k} x^k , \qquad \log(x + x_1) = \log x + \log(1 + x_1/x) , \\
e^{\mathcal{A}x} = \sum_{k = 0}^\infty \frac{1}{k!} \mathcal{A}^k x^k, \qquad x^\mathcal{N} = e^{\mathcal{N} \log x}, \qquad (1 + x)^\mathcal{A} = \sum_{k = 0}^\infty \binom{\mathcal{A}}{k} x^k.
\end{gathered}
\end{equation}
These formulas remain valid when we can make a substitution of a formal series (in possibly multiple variables) that is still well defined.

Given a semisimple operator $ \mathcal{S} $ on $ W $, we define $ x^\mathcal{S} u = x^\alpha u $ when $ \mathcal{S} u = \alpha u $, and then extend this definition to all of $ W $. If an operator $ \mathcal{A} $ has commuting semisimple and locally nilpotent parts $ \mathcal{S} $ and $ \mathcal{N} $, we use $ x^\mathcal{A} $ to mean $ x^\mathcal{S} x^\mathcal{N} $. For example,
\begin{align*}
(x + x_1)^\mathcal{A} = (x + x_1)^\mathcal{S}(x + x_1)^\mathcal{N} = x^\mathcal{A} (1 + x_1/x)^\mathcal{A}.
\end{align*}

Given a formal series 
\begin{align*}
f(x_1, \dots, x_\N) = \sum_{\substack{r_1, \dots, r_\N \in \mathbb{C} \\ k_1 , \dots, k_\N \in \mathbb{Z}_{\geq 0}}} a_{r_1,\dots, r_\N, k_1, \dots, k_\N} x_1^{r_1} \cdots x_\N^{r_\N} (\log x_1)^{k_1} \cdots (\log x_\N)^{k_\N},
\end{align*} 
we can specialize a formal variable to a sum of two formal variables by understanding logarithms of formal variables as explained above in \eqref{eq:formal meanings}. For example,
\begin{align*}
f(x + x_0, \dots, x_\N) = \sum_{\substack{r_1, \dots, r_\N \in \mathbb{C} \\ k_1 , \dots, k_\N \in \mathbb{Z}_{\geq 0}}} a_{r_1,\dots, r_\N, k_1, \dots, k_\N} (x + x_0)^{r_1} \cdots x_\N^{r_\N} (\log (x + x_0))^{k_1} \cdots (\log x_\N)^{k_\N}.
\end{align*}

The argument of $ z \in \mathbb{C}^\times $ is defined as the argument of $ z $ satisfying $ 0 \leq \arg z < 2\pi $. For each $ p \in \mathbb{Z} $, we use $ l_p(z) $ to denote the logarithm of $ z $ given by $ \log |z| +  ( \arg z + 2 p \pi) \iu $.  When a formal series $ f(x_1, \dots, x_\N) $ has countable support, it can be specialized to the series 
\begin{align*}
f^{p_1, \dots, p_\N} (z_1, \dots, z_\N) = \sum_{\substack{r_1, \dots, r_\N \in \mathbb{C} \\ k_1 , \dots, k_\N \in \mathbb{Z}_{\geq 0}}} a_{r_1,\dots, r_\N, k_1, \dots, k_\N} e^{r_1 l_{p_1}(z_1) } \cdots e^{r_\N l_{p_\N}(z_\N) } (l_{p_1}(z_1))^{k_1} \cdots (l_{p_\N}(z_\N))^{k_\N},
\end{align*}
for any $ z_1, \dots, z_\N \in \mathbb{C}^\times $. If a series $ f $ does not have the variables explicitly stated, we will write 
\begin{align*}
f|_{x_1 = z_1, \dots , x_\N = z_\N}^{p_1, \dots , p_\N}
\end{align*}
to denote $ f^{p_1, \dots, p_\N} (z_1, \dots, z_\N) $, again suppressing the $ \log $ variables for brevity. If $ f^{p_1, \dots, p_\N} (z_1, \dots, z_\N) $ converges absolutely in some appropriate region, it converges to a multivalued analytic function $ f(z_1, \dots, z_\N) $ with branches given by $ f^{p_1, \dots, p_\N} (z_1, \dots, z_\N) $. We also write expressions such as
\begin{equation*}
\langle w_0, \mathcal{Y}(w_1,z_1) \cdots \mathcal{Y}(w_\N,z_\N) w_{\N+1} \rangle : = 
\langle w_0, \mathcal{Y}(w_1,x_1) \cdots \mathcal{Y}(w_\N,x_\N) w_{\N+1} \rangle|
_{x_1 = z_1, \dots,  x_\N = z_\N} 
\end{equation*}
to denote the (multivalued) series of complex numbers. Each branch will have the value of $ z_i^r = e^{r \log z_i} $ dependent on the choice of value $ \log z_i $.

Suppose we have operators $ g $ and $ \mathcal{L}_g $ on a space $ W $ such that $ g = e^{2 \pi \iu \mathcal{L}_g} $, and $ \mathcal{L}_g $ has commuting semisimple and locally nilpotent parts $ \mathcal{S}_g $ and $ \mathcal{N}_g $. (This is always possible when $ W $ is the direct sum of finite-dimensional subspaces on which $ g $ acts invertibly.) The eigenvalues of $ \mathcal{S}_g $ are determined modulo $\mathbb{Z} $. We will always assume $ \mathcal{S}_g $ is chosen so that its eigenvalues lie in $ \{ z \in \mathbb{C} : 0 \leq \Re(z) < 1 \} $. 

We will use $ V_{(n)} $ to denote the $ L(0) $-eigenspaces with eigenvalue $ n \in \mathbb{Z} $ of a vertex operator algebra $ V $, and use $ V_+ $ to denote $ \coprod_{n \in \mathbb{Z}_{> 0}} V_{(n)} $. If $ g $ is an automorphism of $ V $, we use $ V^{[\alpha]} $, $ V^{[\alpha]}_{(n)} $ and $ V^{[\alpha]}_+ $ to denote the spaces of generalized eigenvectors of $ g $ (acting on $ V $, $ V_{(n)} $ and $ V_+ $, respectively) with eigenvalue $ e^{2 \pi i \alpha} $. We will use $ W_{[n]} $ to denote the $ L(0) $-generalized eigenspaces of a twisted $ V $-module (see Definition~3.1 of \cite{HUANG2018346}) with eigenvalue $ n \in \mathbb{C} $. We say $ \wt w = n $ is the (conformal) weight of a vector $ w \in W_{[n]} $.

Suppose we have a tensor product $ W_0 \otimes \dots \otimes W_{\N+1} $ of $ \mathbb{C} $-graded spaces graded by weight. Given weight-homogeneous elements $ w_i \in W_i $, $ i = 0, \dots, \N + 1 $, we will frequently use $ \sigma $ to denote $ \Re( \wt w_0 + \dots + \wt w_{\N+1} ) $, i.e., the real part of the weight of the tensor product $ w_0 \otimes \cdots \otimes w_{\N+1} $.

\subsection{Expansions of logarithms} \label{subsec:log expand}
We now explain how the formal meaning for $ \log(x+x_1) $ in \eqref{eq:formal meanings} connects with the analytic meaning.

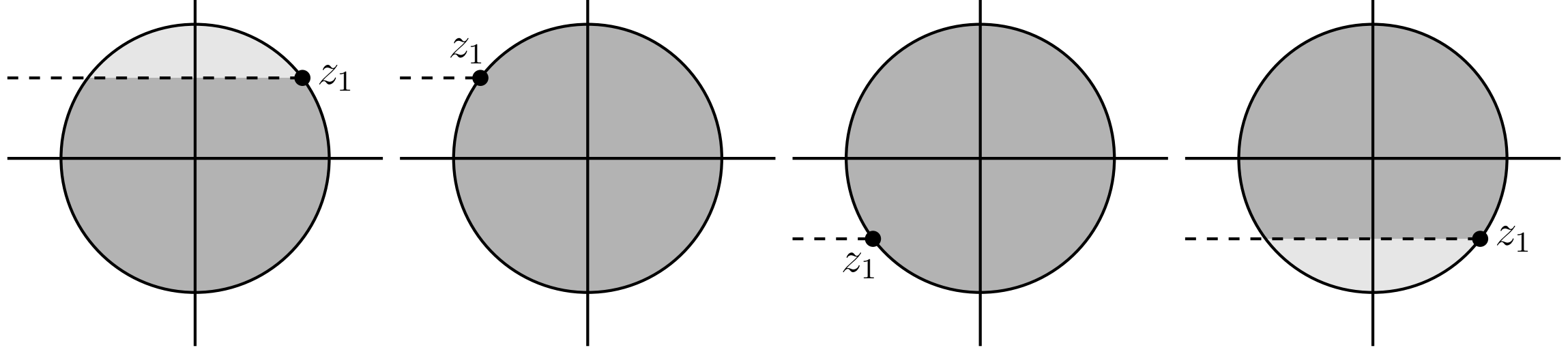
\begin{figure}[b]
\centering
\begin{tikzpicture}[scale=1.3]
\def\R{1}          
\def\ycap{0.6}     
\fill[gray!60] (0,0) circle (\R);
\begin{scope}
  \clip (0,0) circle (\R);
  \fill[gray!20] (-\R,\ycap) rectangle (\R,\R);
\end{scope}
\draw[thick] (0,0) circle (\R);
\draw[thick] (-1.4,0) -- (1.4,0);
\draw[thick] (0,-1.4) -- (0,1.4);
\draw[thick, dashed] (-1.4,\ycap) -- (0.8,\ycap);
\coordinate (z1) at (0.8,0.6);
\fill (z1) circle (0.06);
\node[right] at (z1) {$z_1$};
\end{tikzpicture}
\begin{tikzpicture}[scale=1.3]
\def\R{1}          
\def\ycap{0.6}     
\fill[gray!60] (0,0) circle (\R);
\draw[thick] (0,0) circle (\R);
\draw[thick] (-1.4,0) -- (1.4,0);
\draw[thick] (0,-1.4) -- (0,1.4);
\draw[thick, dashed] (-1.4,\ycap) -- (-0.8,\ycap);
\coordinate (z1) at (-0.8,0.6);
\fill (z1) circle (0.06);
\node at (-0.9,0.8) {$z_1$};
\end{tikzpicture}
\begin{tikzpicture}[scale=1.3]
\def\R{1}          
\def\ycap{0.6}     
\fill[gray!60] (0,0) circle (\R);
\draw[thick] (0,0) circle (\R);
\draw[thick] (-1.4,0) -- (1.4,0);
\draw[thick] (0,-1.4) -- (0,1.4);
\draw[thick, dashed] (-1.4,-\ycap) -- (-0.8,-\ycap);
\coordinate (z1) at (-0.8,-0.6);
\fill (z1) circle (0.06);
\node at (-0.9,-0.8) {$z_1$};
\end{tikzpicture}
\begin{tikzpicture}[scale=1.3]
\def\R{1}          
\def\ycap{0.6}     
\fill[gray!20] (0,0) circle (\R);
\begin{scope}
  \clip (0,0) circle (\R);
  \fill[gray!60] (-\R,-\ycap) rectangle (\R,\R);
\end{scope}
\draw[thick] (0,0) circle (\R);
\draw[thick] (-1.4,0) -- (1.4,0);
\draw[thick] (0,-1.4) -- (0,1.4);
\draw[thick, dashed] (-1.4,-\ycap) -- (0.8,-\ycap);
\coordinate (z1) at (0.8,-0.6);
\fill (z1) circle (0.06);
\node[right] at (z1) {$z_1$};
\end{tikzpicture}

\caption{For fixed $ z_1 \in \mathbb{C}^\times $ (drawn in each quadrant), shaded in dark is the region for $ z_2 $ satisfying $ |z_1| > |z_2| $ and $ -\frac{\pi}{2} < \arg(z_1 - z_2) - \arg z_1 < \frac{\pi}{2} $.} \label{fig:z1>z2 regions}
\end{figure}

Say we want to find the expansion for $ l_p(z_1 - z_2) $ in the region $ |z_1| > |z_2| $. Since $ |z_1| > |z_2| $, we seek to Taylor expand $ l_p(z_1 - z_2) $ in terms of $ z_2 $ at $ 0 $. Fix $ z_1 \in \mathbb{C}^\times $. The function $ f(z_2) =  l_p(z_1 - z_2) $ is discontinuous on the ray $ z_1 - z_2 \in \mathbb{R}_{\geq 0} $, which possibly divides the disk $ \{ z \in \mathbb{C} : |z_1| > |z| \} $ into two connected regions. Suppose that $ z_2 $ lies in the region that contains $ 0 $, that is, $ -\frac{\pi}{2} < \arg(z_1 - z_2) - \arg z_1 < \frac{\pi}{2} $ (see Figure~\ref{fig:z1>z2 regions}). Then, for $ z_1 \notin \mathbb{R}_{>0} $, we have $ f^{(k)}(z_2) = - (k-1)!(z_1 - z_2)^{-k} $ for all $ k \in \mathbb{Z}_{>0} $, and we find that
\begin{align*}
l_p(z_1 - z_2) = f(z_2) = \sum_{k \geq 0} \frac{1}{k!} f^{(k)}(0) z_2^k = l_p(z_1) + \sum_{k > 0} \frac{-1}{k} \left(\frac{z_2}{z_1}\right)^k 
\end{align*}
is the expansion for $ l_p(z_1 - z_2) $ in the region $ |z_1| > |z_2| $, $ -\frac{\pi}{2} < \arg(z_1 - z_2) - \arg z_1 < \frac{\pi}{2} $. Using appropriate one-sided limits, we obtain the same expansions for $ z_1 \in \mathbb{R}_{>0} $.

\begin{figure}[b]
\centering
\begin{tikzpicture}[scale=1.3]
\def\R{1}          
\def\ycap{0.6}     
\fill[gray!20] (0,0) circle (\R);
\draw[thick] (0,0) circle (\R);
\draw[thick] (-1.4,0) -- (1.4,0);
\draw[thick] (0,-1.4) -- (0,1.4);
\draw[thick, dashed] (0.8,\ycap) -- (1.4,\ycap);
\coordinate (z2) at (0.8,0.6);
\fill (z2) circle (0.06);
\node[right] at (0.75,0.8) {$z_2$};
\end{tikzpicture}
\begin{tikzpicture}[scale=1.3]
\def\R{1}          
\def\ycap{0.6}     
\fill[gray!20] (0,0) circle (\R);
\begin{scope}
  \clip (0,0) circle (\R);
  \fill[gray!60] (-\R,\ycap) rectangle (\R,\R);
\end{scope}
\draw[thick] (0,0) circle (\R);
\draw[thick] (-1.4,0) -- (1.4,0);
\draw[thick] (0,-1.4) -- (0,1.4);
\draw[thick, dashed] (-0.8,\ycap) -- (1.4,\ycap);
\coordinate (z2) at (-0.8,0.6);
\fill (z2) circle (0.06);
\node[left] at (z2) {$z_2$};
\end{tikzpicture}
\begin{tikzpicture}[scale=1.3]
\def\R{1}          
\def\ycap{0.6}     
\fill[gray!20] (0,0) circle (\R);
\begin{scope}
  \clip (0,0) circle (\R);
  \fill[gray!60] (-\R,-\ycap) rectangle (\R,\R);
\end{scope}
\draw[thick] (0,0) circle (\R);
\draw[thick] (-1.4,0) -- (1.4,0);
\draw[thick] (0,-1.4) -- (0,1.4);
\draw[thick, dashed] (-0.8,-\ycap) -- (1.4,-\ycap);
\coordinate (z2) at (-0.8,-0.6);
\fill (z2) circle (0.06);
\node at (-0.9,-0.8) {$z_2$};
\end{tikzpicture}
\begin{tikzpicture}[scale=1.3]
\def\R{1}          
\def\ycap{0.6}     
\fill[gray!60] (0,0) circle (\R);
\draw[thick] (0,0) circle (\R);
\draw[thick] (-1.4,0) -- (1.4,0);
\draw[thick] (0,-1.4) -- (0,1.4);
\draw[thick, dashed] (0.8,-\ycap) -- (1.4,-\ycap);
\coordinate (z2) at (0.8,-0.6);
\fill (z2) circle (0.06);
\node[right] at (0.75,-0.8) {$z_2$};
\end{tikzpicture}

\caption{For fixed $ z_2 \in \mathbb{C}^\times $ (drawn in each quadrant), shaded in dark is the region for $ z_1 $ satisfying $ |z_2| > |z_1| $ and $ -\frac{3\pi}{2} < \arg(z_1 - z_2) - \arg z_2 < -\frac{\pi}{2} $.} \label{fig:z2>z1 regions}
\end{figure}
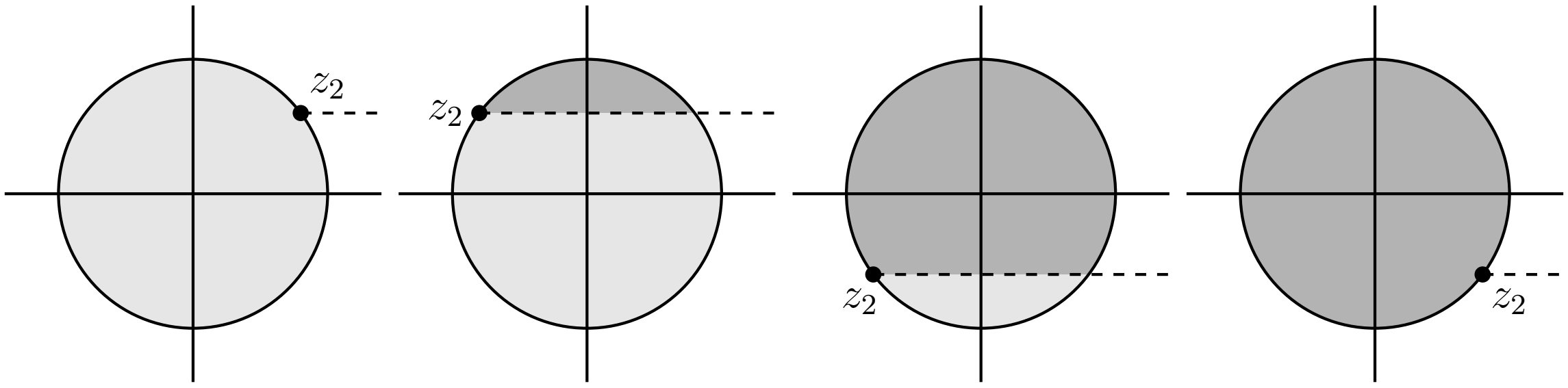

Say we instead want to find the expansion for $ l_p(z_1 - z_2) $ in the region $ |z_2| > |z_1| $. Since $ |z_2| > |z_1| $, we seek to Taylor expand $ l_p(z_1 - z_2) $ in terms of $ z_1 $ at $ 0 $. Fix $ z_2 \in \mathbb{C}^\times $. The function $ f(z_1) =  l_p(z_1 - z_2) $ is discontinuous on the ray $ z_1 - z_2 \in \mathbb{R}_{\geq 0} $, which possibly divides the disk $ \{ z \in \mathbb{C} : |z_2| > |z| \} $ into two connected regions. First consider the case where $ \pi \leq \arg z_2 < 2\pi $, and suppose that $ z_1 $ lies in the region that contains $ 0 $, that is, $ -\frac{3\pi}{2} < \arg(z_1 - z_2) - \arg z_2 < -\frac{\pi}{2} $ (see Figure~\ref{fig:z2>z1 regions}). Then, for $ z_2 \notin \mathbb{R}_{<0} $, $ f^{(k)}(z_1) = (-1)^{k-1} (k-1)!(z_1 - z_2)^{-k} $ for all $ k \in \mathbb{Z}_{>0} $, and 
\begin{align*}
f(z_1) = \sum_{k \geq 0} \frac{1}{k!} f^{(k)}(0) z_1^k = l_p(-z_2) + \sum_{k > 0} \frac{-1}{k} \left(\frac{z_1}{z_2}\right)^k .
\end{align*}
Since $ \pi \leq \arg z_2 < 2\pi $, we have $ l_p(-z_2) = l_p(z_2) - \iu \pi $. Hence,
\begin{align*}
l_p(z_1 - z_2) = l_p(z_2) - \iu \pi + \sum_{k > 0} \frac{-1}{k} \left(\frac{z_1}{z_2}\right)^k 
\end{align*}
is the expansion for $ l_p(z_1 - z_2) $ in the region $ |z_2| > |z_1| $, $ -\frac{3\pi}{2} < \arg(z_1 - z_2) - \arg z_2 < -\frac{\pi}{2} $ and $ \pi \leq \arg z_2 < 2\pi $. Using appropriate one-sided limits, we obtain the same expansions for $ z_2 \in \mathbb{R}_{< 0} $.

We now consider the case  $ 0 \leq \arg z_2 < \pi $. Suppose we want to expand $ l_p(z_1 - z_2) $ in the region $ |z_2| > |z_1| $ with the same restriction $ -\frac{3\pi}{2} < \arg(z_1 - z_2) - \arg z_2 < -\frac{\pi}{2} $ as above. When $ z_2 $ is in the first quadrant, this region is empty. When $ z_2 $ is in the second quadrant, this region does not contain $ 0 $ (see Figure~\ref{fig:z2>z1 regions}). 

In order to Taylor expand $ l_p(z_1 - z_2) $ in terms of $ z_1 $ at $ 0 $, we have to first analytically continue $ f(z_1) =  l_p(z_1 - z_2) $ into the region that contains $ 0 $ along a path that remains in the disk $ \{ z \in \mathbb{C} :  |z_2| > |z|\} $. We find that $ f(z_1) = l_{p-1}(z_1-z_2) $ in the region that contains $ 0 $. So, 
\begin{align*}
f(z_1) = \sum_{k \geq 0} \frac{1}{k!} f^{(k)}(0) z_1^k = l_{p-1}(-z_2) + \sum_{k > 0} \frac{-1}{k} \left(\frac{z_1}{z_2}\right)^k .
\end{align*}
Since $ 0 \leq \arg z_2 < \pi $, we have $ l_{p-1}(-z_2) = l_{p-1}(z_2) + \iu \pi = l_{p}(z_2) - \iu \pi $. Hence,
\begin{align*}
l_p(z_1 - z_2) = l_p(z_2) - \iu \pi + \sum_{k > 0} \frac{-1}{k} \left(\frac{z_1}{z_2}\right)^k 
\end{align*}
is the expansion for $ l_p(z_1 - z_2) $ in the region $ |z_2| > |z_1| $, $ -\frac{3\pi}{2} < \arg(z_1 - z_2) - \arg z_2 < -\frac{\pi}{2} $ and $ 0 \leq \arg z_2 < \pi $.

Similar arguments can be used to find the expansion of $ l_p(z_2 + (z_1-z_2)) $ in the region $ |z_2| > |z_1 - z_2| $, $ - \frac{\pi}{2} < \arg z_1 - \arg z_2 < \frac{\pi}{2} $.

In summary, we have found the following expansions that are consistent with our choice of $ 0 \leq \arg z < 2 \pi $, the definition of $ l_p $, and the formal meaning of $ \log(1+x) $: 
\begin{equation} \label{eq:log x_0 expansion}
\log x_0|_{x_0 = z_1 - z_2}^p = \log(x_1 - x_2)|_{x_1 = z_1, x_2 = z_2}^{p,0} 
\end{equation}
in the region $ |z_1| > |z_2| $, $ - \frac{\pi}{2} < \arg(z_1 - z_2) - \arg z_1 < \frac{\pi}{2} $; 
\begin{equation} \label{eq:log x_0 x2>x1}
\log x_0|_{x_0 = z_1 - z_2}^p = \log(x_2 - x_1)|_{x_1 = z_1, x_2 = z_2}^{0,p} - \iu \pi   
\end{equation}
in the region $ |z_2| > |z_1| $, $ - \frac{3\pi}{2} < \arg(z_1 - z_2) - \arg z_2 < -\frac{\pi}{2} $; and 
\begin{equation}
\log x_1|_{x_1 = z_1}^p = \log(x_2 + x_0)|_{x_0 = z_1 - z_2, x_2 = z_2}^{0,p} 
\end{equation}
in the region $ |z_2| > |z_1 - z_2| $, $ - \frac{\pi}{2} < \arg z_1 - \arg z_2 < \frac{\pi}{2} $.

For all $ n \in \mathbb{C} $, similar arguments show that 
\begin{equation}
x_0^n |_{x_0 = z_1 - z_2}^p = (x_1 - x_2)^n|_{x_1 = z_1, x_2 = z_2}^{p,0} 
\end{equation}
in the region $ |z_1| > |z_2| $, $ - \frac{\pi}{2} < \arg(z_1 - z_2) - \arg z_1 < \frac{\pi}{2} $; 
\begin{equation} \label{eq:power x_0 x2>x1}
x_0^n|_{x_0 = z_1 - z_2}^p = e^{- \iu \pi n} (x_2 - x_1)^n|_{x_1 = z_1, x_2 = z_2}^{0,p}
\end{equation}
in the region $ |z_2| > |z_1| $, $ - \frac{3\pi}{2} < \arg(z_1 - z_2) - \arg z_2 < -\frac{\pi}{2} $; and 
\begin{equation} \label{eq:x1 expansion}
x_1^n|_{x_1 = z_1}^p = (x_2 + x_0)^n|_{x_0 = z_1 - z_2, x_2 = z_2}^{0,p} 
\end{equation}
in the region $ |z_2| > |z_1 - z_2| $, $ - \frac{\pi}{2} < \arg z_1 - \arg z_2 < \frac{\pi}{2} $.

Relations \eqref{eq:log x_0 expansion} to \eqref{eq:x1 expansion} will be implicitly used when defining twisted intertwining operators in Definition~\ref{def:twisted intertwining operator}.

Consistent with \eqref{eq:log x_0 x2>x1} and \eqref{eq:power x_0 x2>x1}, we set the conventions and notations
\begin{gather}
\log( - x_2 + x_1) := \log( x_2 - x_1) - \iu \pi , \\
( - x_2 + x_1)^n := e^{-\iu \pi n} (x_2 - x_1)^n , \quad \text{for all } n \in \mathbb{C}, \\
( - x_2 + x_1)^{\mathcal{A}} := e^{-\iu \pi \mathcal{A}} (x_2 - x_1)^{\mathcal{A}} ,  
\end{gather}
for operators $ \mathcal{A} $ that are sums of semisimple and locally nilpotent parts.

\begin{remark}
When expanding $ l_p(z_1 - z_2) $ and $ (z_1 - z_2)^n $ in the region $ |z_2| > |z_1| $, we find that we also have  
\begin{gather*}
\log x_0|_{x_0 = z_1 - z_2}^p = \log(x_2 - x_1)|_{x_1 = z_1, x_2 = z_2}^{0,p} + \iu \pi , \\
x_0^n|_{x_0 = z_1 - z_2}^p = e^{+ \iu \pi n} (x_2 - x_1)^n|_{x_1 = z_1, x_2 = z_2}^{0,p}.
\end{gather*}
However, this expansion is true when $ |z_2| > |z_1| $ and $ \frac{\pi}{2} < \arg(z_1-z_2) - \arg(z_2) < \frac{3\pi}{2} $, i.e., the lightly shaded region in Figure~\ref{fig:z2>z1 regions}. In \cite{HUANG2018346}, the difference between the two restrictions $ -\frac{3\pi}{2} < \arg(z_1-z_2) - \arg(z_2) < -\frac{\pi}{2} $ and $ \frac{\pi}{2} < \arg(z_1-z_2) - \arg(z_2) < \frac{3\pi}{2} $ is explained. The former restriction is more natural for defining twisted intertwining operators, and, hence, is used when expanding $ l_p(z_1 - z_2) $ and $ (z_1 - z_2)^n $ in the region $ |z_2| > |z_1 | $ throughout this paper.
\end{remark}

\section{Twisted intertwining operators for commuting automorphisms} \label{sec:Jacobi identity}
In this section, we derive a Jacobi identity for intertwining operators among twisted modules twisted by commuting automorphisms (possibly of infinite order and non-semisimple), and show that the Jacobi identity is equivalent to the duality property. We combine techniques from \cite{TanMostlyUntwisted} (which considers the \gtype\ case for a general automorphism $ g $) with techniques from Chapter 5 of \cite{TanThesis} (which considers commuting automorphisms of finite order). 
 
Let $ V = (V, Y, \one, \omega) $ be a vertex operator algebra over $ \mathbb{C} $ (see Definition~3.1.22 of \cite{LepowskyLi}). Recall that an automorphism $ g $ of a vertex operator algebra $ V $ is a linear automorphism of $ V $ such that $ gY(u,x)v = Y(gu,x)gv $ for all $ u,v \in V $, $ g \one = \one $ and $ g\omega = \omega $. Refer to Definition~3.1 of \cite{HUANG2018346} for the notion of a $ g $-twisted $ V $-module. However, we will not require any $ g $-actions (or actions by other automorphisms) on the $ g $-twisted modules here. We recall the definition of an intertwining operator among twisted modules from Definition~4.1 of \cite{HUANG2018346}. 

\begin{definition} \label{def:twisted intertwining operator}
Let $ V $ be a vertex operator algebra with automorphisms $ g_1, g_2, g_3 $. Let $ W_i $ be a $ g_i $-twisted $ V $-module, for $ i = 1,2,3 $. A \emph{twisted intertwining operator of type $\type{W_1}{W_2}{W_3}$} is a linear map
\begin{equation} \label{eq:intertwining operator definition}
\begin{aligned}
\mathcal{Y}( \cdot, x) \cdot : W_1 \otimes W_2 &\to W_3 \{ x \}[\log x], \\
w_1 \otimes w_2 &\mapsto \sum_{k=0}^K \sum_{n \in \mathbb{C}} \mathcal{Y}_{n,k} (w_1) w_2 x^{-n-1} (\log x)^k
\end{aligned}
\end{equation}
satisfying the following conditions:
\begin{enumerate}
\item[(i)] (\emph{Lower truncation}) For all $ w_1 \in W_1 $, $ w_2 \in W_2 $, $ n \in \mathbb{C} $ and $ k = 0, \dots, K $,
\begin{equation*}
\mathcal{Y}_{n+l,k}(w_1)w_2 = 0 \text{ for all sufficiently large } l \in \mathbb{Z}.
\end{equation*} 
\item[(ii)] (\emph{$ L(-1) $-derivative property}) For all $ w_1 \in W_1 $,
\begin{equation*}
\frac{d}{dx} \mathcal{Y}( w_1, x) = \mathcal{Y}( L(-1) w_1, x).
\end{equation*}
\item[(iii)] (\emph{Duality}) For all $ u \in V $, $ w_1 \in W_1 $, $ w_2 \in W_2 $ and $ w_3' \in W_3' $, there exists a formal series
\begin{equation}\label{eq:converge function}
\begin{aligned} 
f(x_0,x_1,x_2) &= \sum_{i,j,k,l,m,n = 0}^N a_{ijklmn} x_0^{r_i} x_1^{s_j} x_2^{t_k}  (\log x_0)^l (\log x_1)^m (\log x_2)^n \\
&\in \mathbb{C}\{x_0,x_1,x_2\}[\log x_0, \log x_1, \log x_2 ] 
\end{aligned}
\end{equation}
such that for $ p_1, p_2, p_{12} \in \mathbb{Z} $, the series
\begin{equation} \label{eq:multiseries}
\begin{aligned}
 \langle w_3', (Y_{W_3})^{p_1}(u,z_1) \mathcal{Y}^{p_2}(w_1,z_2) w_2 \rangle &= \langle w_3', Y_{W_3}(u,x_1) \mathcal{Y}(w_1,x_2) w_2 \rangle |_{x_1 = z_1, x_2 = z_2}^{p_1,p_2} , \\
 \langle w_3',  \mathcal{Y}^{p_2}(w_1,z_2) (Y_{W_2})^{p_1}(u,z_1)w_2 \rangle &= \langle w_3',  \mathcal{Y}(w_1,x_2) Y_{W_2}(u,x_1)w_2 \rangle|_{x_1 = z_1, x_2 = z_2}^{p_1,p_2} , \\
 \langle w_3',  \mathcal{Y}^{p_2}((Y_{W_1})^{p_{12}}(u,z_1 - z_2)w_1,z_2) w_2 \rangle &= \langle w_3',  \mathcal{Y}(Y_{W_1}(u,x_0)w_1,x_2) w_2 \rangle |_{x_0 = z_1-z_2, x_2 = z_2}^{p_{12},p_2} ,
 \end{aligned}
\end{equation}
are absolutely convergent in the regions $ |z_1| > |z_2| > 0 $, $|z_2| > |z_1| > 0 $ and $ |z_2| > |z_1 - z_2| > 0 $, respectively. Moreover, they converge to 
\begin{equation}\label{eq:branches}
\begin{gathered}
f^{p_1,p_1,p_2}(z_1 - z_2,z_1,z_2), \quad f^{p_2,p_1,p_2}(z_1-z_2,z_1,z_2), \quad f^{p_{12}, p_2, p_2}(z_1-z_2,z_1,z_2) 
\end{gathered}
\end{equation} 
in the regions 
\begin{equation}\label{eq:regions}
\begin{gathered}
|z_1| > |z_2| > 0 \text{ and } -\frac{\pi}{2} < \arg(z_1 - z_2) - \arg(z_1)  < \frac{\pi}{2} ,\\
|z_2| > |z_1| > 0 \text{ and } - \frac{3 \pi }{2} < \arg(z_1 - z_2) - \arg(z_2) < -\frac{\pi}{2} ,\\
|z_2| > |z_1 - z_2| > 0 \text{ and } -\frac{\pi}{2} < \arg(z_1) - \arg(z_2)  < \frac{\pi}{2} , \\
\end{gathered}
\end{equation} 
respectively. 
\end{enumerate}
\end{definition}

\begin{remark}
We remark that there is another definition for twisted intertwining operators in \cite{DuHuang} that does not assume the specific form of $ f(x_0,x_1,x_2) $ as in the definition above, and instead assumes a general analytic function on an appropriate domain. However, when the automorphisms $ g_1 $ and $ g_2 $ commute (as will always be the case in this paper), the analytic function in the definition of \cite{DuHuang} has the same form as in Definition~\ref{def:twisted intertwining operator} above. See Chapter 2 of \cite{TanThesis} for more discussion.
\end{remark}

Suppose we have a set $ \{ g_i \}_i $ of commuting automorphisms of $ V $. Since each automorphism $ g_i $ commutes with the action of $ L(0) $ on $ V $, it restricts to a linear automorphism $ g_i|_{V_{(n)}} $ on each finite-dimensional subspace $ V_{(n)} $, $ n \in \mathbb{Z} $. Since $ \{ g_i|_{V_{(n)}} \}_i $ are invertible operators, there exist semisimple operators $ \mathcal{S}_{g_i}|_{V_{(n)}} $ and nilpotent operators $ \mathcal{N}_{g_i}|_{V_{(n)}} $ such that 
\begin{equation*}
g_i|_{V_{(n)}} = e^{2 \pi \iu \left( \mathcal{S}_{g_i}|_{V_{(n)}} + \mathcal{N}_{g_i}|_{V_{(n)}} \right) }.
\end{equation*}  
Furthermore, $ \mathcal{S}_{g_i}|_{V_{(n)}} $ and $ \mathcal{N}_{g_i}|_{V_{(n)}} $ commute, and $ \mathcal{S}_{g_i}|_{V_{(n)}} $ can be chosen to have eigenvalues with real parts in $ [0,1) $. Since $ \{ g_i|_{V_{(n)}}\}_i $ is mutually commuting, the set of operators $ \{\mathcal{S}_{g_i}|_{V_{(n)}}, \mathcal{N}_{g_i}|_{V_{(n)}} \}_i $ is also mutually commuting. We can then define the operators $ \mathcal{S}_{g_i} $ and $ \mathcal{N}_{g_i} $ on $ V $ as the direct sums of these operators over all subspaces $ V_{(n)} $, $ n \in \mathbb{Z} $. Thus, all $ g_i  $ have normalized logarithms $ \mathcal{L}_{g_i} = \mathcal{S}_{g_i} + \mathcal{N}_{g_i} $ with semisimple and locally nilpotent parts that are all mutually commuting.

There are two main insights that build intuition for the following proof. First, $ \{\mathcal{S}_{g_i}, \mathcal{N}_{g_i} \}_i $ is a mutually commuting set. Second, twisted vertex operators $ Y_W(u,x) $ for a $ g $-twisted $ V $-module $ (W, Y_W) $ have single-valued forms $ Y_W(x^{\mathcal{L}_g}u,x) $ (see Proposition~2.3.2 and Remark~2.3.3 of \cite{TanThesis}, Lemma~2.3 of \cite{Huang2016AssociativeAF}, or Proposition~1.2 of \cite{TanCofinite}). 

\begin{proposition}\label{prop:Jacobi c5}
Let $ V $ be a vertex operator algebra. Let $ g_1, g_2, g_3 \in \Aut V $ be commuting automorphisms with $ g_3 = g_1 g_2 $. Let $ W_i $ be a $ g_i $-twisted $ V $-module for each $ i = 1,2,3 $. Let $ \mathcal{Y}: W_1 \otimes W_2 \to W_3 \{x\}[\log x] $ be a twisted intertwining operator. Let $ u \in V $ and $ w_1 \in W_1 $. Then $ \mathcal{Y} $ satisfies the Jacobi identity
\begin{equation}\label{eq:twisted jacobi c5}
\begin{aligned}
&x_0^{-1} \delta \left( \frac{x_1 - x_2}{x_0} \right)  Y_{W_3}\left( \left( \frac{x_1 - x_2}{x_0} \right)^{\mathcal{L}_{g_1}} u,x_1\right) \mathcal{Y}(w_1, x_2) \\
& \qquad - x_0^{-1} \delta \left( \frac{- x_2 + x_1}{x_0} \right)  \mathcal{Y}(w_1, x_2) Y_{W_2}\left( \left( \frac{- x_2 + x_1}{x_0} \right)^{\mathcal{L}_{g_1}} u,x_1\right) \\
&\quad = x_1^{-1} \delta \left( \frac{x_2 + x_0}{x_1} \right)  \mathcal{Y}\left(Y_{W_1}\left( \left( \frac{x_2 + x_0}{x_1} \right)^{\mathcal{L}_{g_2}} u,x_0\right)w_1,x_2\right),
\end{aligned}
\end{equation}
where $ (-x_2 + x_1)^{\mathcal{A}} := (-1)^{\mathcal{A}} (x_2 - x_1)^{\mathcal{A}} $ and $ (-1)^{\mathcal{A}} := e^{- \pi \iu \mathcal{A}} $ for all operators $ \mathcal{A} $ that are sums of semisimple and locally nilpotent parts.
\end{proposition}
\begin{proof}
For each $ u \in V $, $ w_1 \in W_1 $, $w_2 \in W_2 $ and $ w_3' \in W_3' $, we denote by $ f(u,w_1,w_2,w_3'; x_0,x_1,x_2) $ a generalized formal series as in \eqref{eq:converge function} such that 
\begin{gather*}
\langle w_3', Y_{W_3}(u,x_1) \mathcal{Y}(w_1,x_2) w_2 \rangle |_{x_1 = z_1, x_2 = z_2}^{p_1,p_2} , \qquad \langle w_3',  \mathcal{Y}(w_1,x_2) Y_{W_2}(u,x_1)w_2 \rangle|_{x_1 = z_1, x_2 = z_2}^{p_1,p_2} , \\
\langle w_3',  \mathcal{Y}(Y_{W_1}(u,x_0)w_1,x_2) w_2 \rangle |_{x_0 = z_1-z_2, x_2 = z_2}^{p_{12},p_2} ,
\end{gather*}
converge to 
\begin{gather*}
f^{p_1,p_1,p_2}(u,w_1,w_2,w_3';z_1 - z_2,z_1,z_2), \qquad f^{p_2,p_1,p_2}(u,w_1,w_2,w_3';z_1-z_2,z_1,z_2), \\
f^{p_{12}, p_2, p_2}(u,w_1,w_2,w_3';z_1-z_2,z_1,z_2),
\end{gather*}
respectively, in the respective regions as in the definition for twisted intertwining operators. We now fix $ u $, $ w_1 $, $ w_2 $ and $ w_3' $. Let $ u = \sum_i u^{(i)} $ be a decomposition of $ u $ into finitely many simultaneous eigenvectors $ u^{(i)} $ for $ \mathcal{S}_{g_1} $ and $ \mathcal{S}_{g_2} $ with eigenvalues $ \alpha^{(i)}_1 $ and $ \alpha^{(i)}_2 $, respectively. We define the generalized formal series
\begin{equation}
\begin{aligned}
g(x_0,x_1,x_2) &= f\left(x_0^{\mathcal{L}_{g_1}} x_1^{\mathcal{L}_{g_2}} u,w_1,w_2,w_3'; x_0,x_1,x_2 \right) \\
&= \sum_i \sum_{j,k=0}^K \frac{1}{j!k!} x_0^{\alpha^{(i)}_1} x_1^{\alpha^{(i)}_2} (\log x_0)^j (\log x_1)^k  f\left(\mathcal{N}_{g_1}^j \mathcal{N}_{g_2}^k u^{(i)} ,w_1,w_2,w_3'; x_0,x_1,x_2 \right) ,
\end{aligned}
\end{equation}
where $ K $ is sufficiently large so that $ \mathcal{N}_{g_1}^{K+1} u^{(i)} = \mathcal{N}_{g_2}^{K+1} u^{(i)} = 0 $ for all $ i $.
By choosing $ p_{12}, p_1 ,p_2 $ to be equal to some common $ p \in \mathbb{Z} $, we see that
\begin{gather*}
\begin{aligned}
g(x_1 - x_2,x_1,x_2)|_{x_1 = z_1, x_2 = z_2}^{p,p} &= g(x_0,x_1,x_2)|_{x_0 = z_1 - z_2, x_1 = z_1, x_2 = z_2}^{p,p,p} \\
&= \left.\left\langle w_3', Y_{W_3}\left( (x_1 - x_2)^{\mathcal{L}_{g_1}} x_1^{\mathcal{L}_{g_2}} u,x_1 \right) \mathcal{Y}(w_1,x_2) w_2 \right\rangle \right|_{x_1 = z_1, x_2 = z_2}^{p,p}  
\end{aligned} \\
\text{in }  |z_1| > |z_2| > 0 \text{ and } -\frac{\pi}{2} < \arg(z_1 - z_2) - \arg(z_1)  < \frac{\pi}{2} , \\
\begin{aligned}
&\left. \left( \left. e^{- \pi \iu x_0 \frac{\partial}{\partial x_0}} g(x_0,x_1,x_2) \right|_{x_0 = x_2 - x_1}\right) \right|_{x_1 = z_1, x_2 = z_2}^{p,p} \\ 
&\quad = \left. e^{- \pi \iu x_0 \frac{\partial}{\partial x_0}} g(x_0,x_1,x_2) \right|_{x_0 = z_2 - z_1, x_1 = z_1, x_2 = z_2}^{p,p,p} \\
&\quad = \left. \left\langle w_3',  \mathcal{Y}(w_1,x_2) Y_{W_2}\left( e^{- \pi \iu \mathcal{L}_{g_1}} x_0^{\mathcal{L}_{g_1}} x_1^{\mathcal{L}_{g_2}} u, x_1 \right)w_2 \right\rangle \right|_{x_0 = z_2 - z_1, x_1 = z_1, x_2 = z_2}^{p,p,p} \\
&\quad = \left. \left\langle w_3',  \mathcal{Y}(w_1,x_2) Y_{W_2}\left( e^{- \pi \iu \mathcal{L}_{g_1}} (x_2 - x_1)^{\mathcal{L}_{g_1}} x_1^{\mathcal{L}_{g_2}} u, x_1 \right)w_2 \right\rangle \right|_{x_1 = z_1, x_2 = z_2}^{p,p} 
\end{aligned} \\
\text{in } |z_2| > |z_1| > 0, \  - \frac{3 \pi }{2} < \arg(z_1 - z_2) - \arg(z_2) < -\frac{\pi}{2}, \\
\begin{aligned}
g(x_0,x_2 + x_0,x_2)|_{x_0 = z_1 - z_2, x_2 = z_2}^{p,p} &= g(x_0,x_1,x_2)|_{x_0 = z_1 - z_2, x_1 = z_1, x_2 = z_2}^{p,p,p}
\\
&= \left. \left\langle w_3',  \mathcal{Y}\left( Y_{W_1} \left( x_0^{\mathcal{L}_{g_1}} (x_2 + x_0)^{\mathcal{L}_{g_2}}u, x_0 \right)w_1,x_2 \right) w_2 \right\rangle \right|_{x_0 = z_1-z_2, x_2 = z_2}^{p,p} 
\end{aligned} \\
\text{in } |z_2| > |z_1 - z_2| > 0 \text{ and } -\frac{\pi}{2} < \arg(z_1) - \arg(z_2)  < \frac{\pi}{2} .
\end{gather*}
By arguments of unique expansion sets, we see that 
\begin{align*}
g(x_1-x_2,x_1, x_2) &= \left\langle w_3', Y_{W_3}\left( (x_1 - x_2)^{\mathcal{L}_{g_1}} x_1^{\mathcal{L}_{g_2}} u,x_1 \right) \mathcal{Y}(w_1,x_2) w_2 \right\rangle , \\
\left. e^{- \pi \iu x_0 \frac{\partial}{\partial x_0}} g(x_0,x_1, x_2) \right|_{x_0 = x_2 - x_1} &= \left\langle w_3',  \mathcal{Y}(w_1,x_2) Y_{W_2}\left( e^{- \pi \iu \mathcal{L}_{g_1}} (x_2 - x_1)^{\mathcal{L}_{g_1}} x_1^{\mathcal{L}_{g_2}} u, x_1 \right) w_2 \right\rangle  , \\
g(x_0,x_2 + x_0, x_2) &= \left\langle w_3',  \mathcal{Y}\left( Y_{W_1} \left( x_0^{\mathcal{L}_{g_1}} (x_2 + x_0)^{\mathcal{L}_{g_2}}u, x_0 \right)w_1,x_2 \right) w_2 \right\rangle .
\end{align*}
Define a formal series $ h(x_0,x_1,x_2) $ in $ \mathbb{C}\{ x_2 \}[x_0, x_0^{-1}, x_1, x_1^{-1} , \log x_2 ] $, obtained by retaining the terms of $ g(x_0,x_1,x_2) $ with only integral powers of $ x_0 $, $ x_1 $ and no logarithms of $ x_0 $ or $ x_1 $. Since 
\begin{align*}
Y_{W_2}\left( e^{- \pi \iu \mathcal{L}_{g_1}} (x_2 - x_1)^{\mathcal{L}_{g_1}} x_1^{\mathcal{L}_{g_2}} u, x_1 \right)
\end{align*}
has no logarithms of $ x_1 $ and only integral powers of $ x_1 $, and 
\begin{align*}
Y_{W_1} \left( x_0^{\mathcal{L}_{g_1}} (x_2 + x_0)^{\mathcal{L}_{g_2}}u, x_0 \right)
\end{align*}
has no logarithms of $ x_0 $ and only integral powers of $ x_0 $, we see that $ g(z_1-z_2,z_1, z_2) $ is single-valued as $ z_1 $ circles $ z_2 $ and as $ z_1 $ circles $ 0 $. Since $ h(z_1 - z_2, z_1, z_2) - g(z_1 - z_2, z_1,z_2) $ is single-valued as $ z_1 $ circles $ z_2 $ and as $ z_1 $ circles $ 0 $, and is a finite linear combination of terms of the form $ (z_1 - z_2)^\alpha z_1^\beta z_2^\gamma (\log (z_1 - z_2))^i (\log z_1)^j (\log z_2)^k $ with $ \alpha \notin \mathbb{Z} $ or $ \beta \notin \mathbb{Z} $ or $ i > 0 $ or $ j > 0 $, we see that it is equal to $ 0 $. So, $ h(x_1-x_2,x_1, x_2) = g(x_1-x_2,x_1, x_2)$. Similarly $ \left. e^{- \pi \iu x_0 \frac{\partial}{\partial x_0}} h(x_0,x_1, x_2) \right|_{x_0 = x_2 - x_1} = \left. e^{- \pi \iu x_0 \frac{\partial}{\partial x_0}} g(x_0,x_1, x_2) \right|_{x_0 = x_2 - x_1} $, and $ h(x_0,x_2 + x_0, x_2) = g(x_0,x_2 + x_0, x_2) $. Furthermore, we see that $ h(x_0,x_1,x_2) $ is both upper- and lower-truncated in discrete powers of $ x_2 $. Hence, there are only finitely many powers of $ x_2 $. So, we can safely multiply the triple $ \delta $-function identity
\begin{align*}
&x_0^{-1} \delta \left( \frac{x_1 - x_2}{x_0} \right) - x_0^{-1} \delta \left( \frac{- x_2 + x_1}{x_0} \right) = x_1^{-1} \delta \left( \frac{x_2 + x_0}{x_1} \right)  
\end{align*}
by $ h(x_0,x_1,x_2) $ to find
\begin{align*}
&x_0^{-1} \delta \left( \frac{x_1 - x_2}{x_0} \right) h(x_0,x_1,x_2) - x_0^{-1} \delta \left( \frac{- x_2 + x_1}{x_0} \right)  h(x_0,x_1,x_2) \\
&\quad = x_1^{-1} \delta \left( \frac{x_2 + x_0}{x_1} \right)  h(x_0,x_1,x_2), \\
&x_0^{-1} \delta \left( \frac{x_1 - x_2}{x_0} \right) h(x_1-x_2,x_1,x_2) - x_0^{-1} \delta \left( \frac{- x_2 + x_1}{x_0} \right)  h(-x_2 + x_1 ,x_1,x_2), \\
&\quad = x_1^{-1} \delta \left( \frac{x_2 + x_0}{x_1} \right)  h(x_0,x_2 + x_0,x_2), \\
&x_0^{-1} \delta \left( \frac{x_1 - x_2}{x_0} \right) h(x_1-x_2,x_1,x_2) - x_0^{-1} \delta \left( \frac{- x_2 + x_1}{x_0} \right) \left. e^{- \pi \iu x_0 \frac{\partial}{\partial x_0}} h(x_0,x_1, x_2) \right|_{x_0 = x_2 - x_1}  , \\
&\quad = x_1^{-1} \delta \left( \frac{x_2 + x_0}{x_1} \right)  h(x_0,x_2 + x_0,x_2).
\end{align*}
Hence,
\begin{align*}
&x_0^{-1} \delta \left( \frac{x_1 - x_2}{x_0} \right)  \left\langle w_3', Y_{W_3}\left( (x_1 - x_2)^{\mathcal{L}_{g_1}} x_1^{\mathcal{L}_{g_2}} u,x_1 \right) \mathcal{Y}(w_1,x_2) w_2 \right\rangle \\
& \qquad - x_0^{-1} \delta \left( \frac{- x_2 + x_1}{x_0} \right)  \left\langle w_3',  \mathcal{Y}(w_1,x_2) Y_{W_2}\left( e^{- \pi \iu \mathcal{L}_{g_1}} (x_2 - x_1)^{\mathcal{L}_{g_1}} x_1^{\mathcal{L}_{g_2}} u, x_1 \right) w_2 \right\rangle \\
&\quad = x_1^{-1} \delta \left( \frac{x_2 + x_0}{x_1} \right) \left\langle w_3',  \mathcal{Y}\left( Y_{W_1} \left( x_0^{\mathcal{L}_{g_1}} (x_2 + x_0)^{\mathcal{L}_{g_2}}u, x_0 \right)w_1,x_2 \right) w_2 \right\rangle .
\end{align*}
Since this equation is true for all $ w_2 \in W_2 $ and $ w_3' \in W_3' $, it is equivalent to 
\begin{equation}\label{eq:twisted jacobi other form}
\begin{aligned}
&x_0^{-1} \delta \left( \frac{x_1 - x_2}{x_0} \right) Y_{W_3}\left( (x_1 - x_2)^{\mathcal{L}_{g_1}} x_1^{\mathcal{L}_{g_2}} u,x_1 \right) \mathcal{Y}(w_1,x_2) \\
& \qquad - x_0^{-1} \delta \left( \frac{- x_2 + x_1}{x_0} \right)  \mathcal{Y}(w_1,x_2) Y_{W_2}\left( e^{- \pi \iu \mathcal{L}_{g_1}} (x_2 - x_1)^{\mathcal{L}_{g_1}} x_1^{\mathcal{L}_{g_2}} u, x_1 \right)  \\
&\quad = x_1^{-1} \delta \left( \frac{x_2 + x_0}{x_1} \right) \mathcal{Y}\left( Y_{W_1} \left( x_0^{\mathcal{L}_{g_1}} (x_2 + x_0)^{\mathcal{L}_{g_2}}u, x_0 \right)w_1,x_2 \right) .
\end{aligned}
\end{equation}
This is equivalent to the desired equation \eqref{eq:twisted jacobi c5} after replacing $ u $ with $ x_0^{- \mathcal{L}_{g_1} } x_1^{- \mathcal{L}_{g_2}} u $.
\end{proof}

\begin{remark}
Note that the mutual commutativity of $ \mathcal{S}_{g_1} $, $ \mathcal{S}_{g_2} $, $ \mathcal{N}_{g_1} $ and $ \mathcal{N}_{g_2} $ was used implicitly throughout the proof above. It remains a necessary property for the derivation of a Jacobi identity of this form. In the case that $ g_1 $ and $ g_2 $ do not commute, the insertion of $ x_0^{\mathcal{L}_{g_1}} x_1^{\mathcal{L}_{g_2}} $ does not in general obtain single-valuedness. Single-valuedness is important in the proof above because we have essentially performed a contour integration in the variable $ z_1 $ on the Riemann sphere with punctures at $ 0 $, $ z_2 $ and $ \infty $, which is in fact how we get the ``Jacobi identity shape" in the Jacobi identity. So, when $ g_1 $ and $ g_2 $ do not commute, we cannot perform such contour integration of a single-valued function on the Riemann sphere. To perform such contour integration in this case, one would need to lift to a single-valued meromorphic form on a different Riemann surface to obtain an identity analogous to the Jacobi identity.
\end{remark}

\begin{remark}
When we assume all automorphisms have finite order, the Jacobi identity \eqref{eq:twisted jacobi c5} implies the Jacobi identity that appears in Definition~3.1 from \cite{Xu}. Even in the finite order case, instead of simply using this definition without justification, we must prove that \eqref{eq:twisted jacobi c5} follows from the definitions of a twisted intertwining operator in \cite{HUANG2018346} and \cite{DuHuang}, and determine the correct branch for the meaning of $ (-1)^{m} $ when $ m \notin \mathbb{Z} $. 

When we allow the automorphisms to possibly have infinite order, we derive the Jacobi identity assumed in the definition for twisted intertwining operators among modules twisted by commuting elements found in Definition~2.3 of \cite{HuangNahmInequality} (where a different choice was made for the meaning of $ (-1)^{m} $, which does not affect the results). 
\end{remark}

Since $ \omega $ is fixed by $ g_1 $ and $ g_2 $, we have the usual relation
\begin{align*}
&x_0^{-1} \delta \left( \frac{x_1 - x_2}{x_0} \right) Y_{W_3}(\omega,x_1) \mathcal{Y}(w_1, x_2) - x_0^{-1} \delta \left( \frac{- x_2 + x_1}{x_0} \right) \mathcal{Y}(w_1, x_2) Y_{W_2}(\omega,x_1) \\
&\quad = x_1^{-1} \delta \left( \frac{x_2 + x_0}{x_1} \right) \mathcal{Y}(Y_{W_1}(\omega,x_0)w_1,x_2),
\end{align*}
which, together with the $ L(-1) $-derivative property, implies the usual commutator relations for the action of the Virasoro algebra, such as 
\begin{equation*}
[L(0), \mathcal{Y}(w_1, x_2)] = \mathcal{Y}(L(0)w_1, x_2) + x_2 \frac{d}{d x_2}\mathcal{Y}(w_1, x_2) ,
\end{equation*}
directly from the Jacobi identity. This implies that for all $ u \in V $, $ w_1 \in W_1 $, $ w_2 \in W_2 $, $ w_3' \in W_3' $ of homogeneous weight,  we have
\begin{equation}\label{eq:L(0) power relation c5}
\left(x_1 \frac{\partial}{\partial x_1} + x_2 \frac{\partial}{\partial x_2} - \Delta\right)^k \langle w_3', Y_{W_3}(u,x_1) \mathcal{Y}(w_1, x_2) w_2 \rangle = 0
\end{equation}
for all sufficiently large $ k \in \mathbb{Z}_{>0} $, where $ \Delta = \wt w_3' - \wt u - \wt w_1 - \wt w_2 $. Note that this property follows directly from the Jacobi identity \eqref{eq:twisted jacobi c5} without the need to assume the duality property. We now show that the duality property follows from the Jacobi identity \eqref{eq:twisted jacobi c5} in the presence of the other axioms in the definition of twisted intertwining operators.

\begin{proposition}\label{prop:duality equivalence}
Let $ V $ be a vertex operator algebra. Let $ g_1, g_2, g_3 \in \Aut V $ be commuting automorphisms with $ g_3 = g_1 g_2 $. Let $ W_i $ be a $ g_i $-twisted $ V $-module for each $ i = 1,2,3 $. Let $ \mathcal{Y}: W_1 \otimes W_2 \to W_3 \{x\}[\log x] $ be a linear map satisfying the lower truncation condition and $ L(-1) $-derivative property (as in Definition~\ref{def:twisted intertwining operator}). Assume that for all $ u \in V $ and $ w_1 \in W_1 $, we have
\begin{equation} \label{eq:twisted jacobi as an assumption}
\begin{aligned}
&x_0^{-1} \delta \left( \frac{x_1 - x_2}{x_0} \right)  Y_{W_3}\left( \left( \frac{x_1 - x_2}{x_0} \right)^{\mathcal{L}_{g_1}} u,x_1\right) \mathcal{Y}(w_1, x_2) \\
& \qquad - x_0^{-1} \delta \left( \frac{- x_2 + x_1}{x_0} \right)  \mathcal{Y}(w_1, x_2) Y_{W_2}\left( \left( \frac{- x_2 + x_1}{x_0} \right)^{\mathcal{L}_{g_1}} u,x_1\right) \\
&\quad = x_1^{-1} \delta \left( \frac{x_2 + x_0}{x_1} \right)  \mathcal{Y}\left(Y_{W_1}\left( \left( \frac{x_2 + x_0}{x_1} \right)^{\mathcal{L}_{g_2}} u,x_0\right)w_1,x_2\right),
\end{aligned}
\end{equation}
where $ (-x_2 + x_1)^{\mathcal{A}} := (-1)^{\mathcal{A}} (x_2 - x_1)^{\mathcal{A}} $ and $ (-1)^{\mathcal{A}} := e^{- \pi \iu \mathcal{A}} $ for all operators $ \mathcal{A} $ that are sums of semisimple and locally nilpotent parts. Then $ \mathcal{Y} $ satisfies the duality condition for a twisted intertwining operator (as in Definition~\ref{def:twisted intertwining operator}).
\end{proposition}

\begin{proof}
Let $ U $ be a subspace of $ V $ containing $ u $ that is invariant under the actions of $ \mathcal{L}_{g_1} $ and $ \mathcal{L}_{g_2} $. We can always choose $ U $ to be finite dimensional. 

By the lower truncation property of $ Y_{W_1} $, we can pick some $ m \in \mathbb{Z} $ sufficiently high so that $ \res_{x_0} x_0^{m+l} Y_{W_1}\left( x_0^{\mathcal{L}_{g_1}}v,x_0 \right) w_1 = 0 $ for all $ l \in \mathbb{Z}_{\geq 0} $ and $ v \in U $.  Then \eqref{eq:twisted jacobi as an assumption} (or its equivalent form \eqref{eq:twisted jacobi other form}) implies
\begin{align*}
0 &= \res_{x_0} x_0^{m} x_1^{-1} \delta \left( \frac{x_2 + x_0}{x_1} \right) \mathcal{Y}\left( Y_{W_1} \left( x_0^{\mathcal{L}_{g_1}} (x_2 + x_0)^{\mathcal{L}_{g_2}}u, x_0 \right)w_1,x_2 \right)  \\
&= \res_{x_0} x_0^{m} \left( x_0^{-1} \delta \left( \frac{x_1 - x_2}{x_0} \right) Y_{W_3}\left( (x_1 - x_2)^{\mathcal{L}_{g_1}} x_1^{\mathcal{L}_{g_2}} u,x_1 \right) \mathcal{Y}(w_1,x_2) \right. \\
& \qquad \left. - x_0^{-1} \delta \left( \frac{- x_2 + x_1}{x_0} \right)  \mathcal{Y}(w_1,x_2) Y_{W_2}\left( (-x_2 + x_1)^{\mathcal{L}_{g_1}} x_1^{\mathcal{L}_{g_2}} u, x_1 \right)  \ \right) \\
&=  Y_{W_3}\left( (x_1 - x_2)^{\mathcal{L}_{g_1}+m} x_1^{\mathcal{L}_{g_2}} u,x_1 \right) \mathcal{Y}(w_1,x_2)  \\
& \qquad  - \mathcal{Y}(w_1,x_2) Y_{W_2}\left( (-x_2 + x_1)^{\mathcal{L}_{g_1}+m} x_1^{\mathcal{L}_{g_2}} u, x_1 \right) 
\end{align*}
Hence, for all $ w_2 \in W_2 $ and $ w_3' \in W_3' $, we have
\begin{align*}
 & \left\langle w_3', Y_{W_3}\left( (x_1 - x_2)^{\mathcal{L}_{g_1}+m} x_1^{\mathcal{L}_{g_2}} u,x_1 \right) \mathcal{Y}(w_1,x_2) w_2 \right\rangle \\
 & \qquad =  \left\langle w_3', \mathcal{Y}(w_1,x_2) Y_{W_2}\left( (-x_2 + x_1)^{\mathcal{L}_{g_1}+m} x_1^{\mathcal{L}_{g_2}} u, x_1 \right) w_2 \right\rangle .
\end{align*}
This series has no logarithms of $ x_1 $ and only integral powers in $ x_1 $. Furthermore, it is lower truncated in $ x_1 $ and $ x_2 $, and hence also upper truncated in $ x_1 $ and $ x_2 $ by \eqref{eq:L(0) power relation c5}. So, we see that it is equal to some series
\begin{align*}
h_1(u;x_1,x_2) = \sum_{i,j,k = 0}^N a_{ijk} x_1^{s_i} x_2^{t_j} (\log x_2)^k ,
\end{align*}
for some $ s_i \in \mathbb{Z} $ and $ a_{ijk}, t_j \in \mathbb{C} $, implicitly dependent on $ w_1 $, $ w_2 $ and $ w_3' $.  Define $ f_1(x_0,x_1,x_2) = h_1 \left(x_0^{-\mathcal{L}_{g_1}-m} x_1^{-\mathcal{L}_{g_2}} u ; x_1,x_2 \right) $, which we can do since the choice of $ m $ is uniform over $ U $. Then
\begin{equation}
\begin{aligned}
 \langle w_3', Y_{W_3}(u,x_1) \mathcal{Y}(w_1,x_2) w_2 \rangle |_{x_1 = z_1, x_2 = z_2}^{p_1,p_2} &= f_1(x_1 - x_2,x_1,x_2)|_{x_1 = z_1, x_2 = z_2}^{p_1,p_2}, \\
 \langle w_3',  \mathcal{Y}(w_1,x_2) Y_{W_2}(u,x_1)w_2 \rangle|_{x_1 = z_1, x_2 = z_2}^{p_1,p_2} &= \left( e^{- \pi \iu x_0 \frac{\partial}{\partial x_0}} f_1(x_0,x_1, x_2)|_{x_0 = x_2 - x_1} \right)|_{x_1 = z_1, x_2 = z_2}^{p_1,p_2} , 
 \end{aligned}
\end{equation}
are absolutely convergent in the regions $ |z_1| > |z_2| > 0 $ and $|z_2| > |z_1| > 0 $, respectively. Moreover, they converge to 
\begin{equation}
\begin{gathered}
f_1^{p_1,p_1,p_2}(z_1 - z_2,z_1,z_2) \quad \text{and} \quad f_1^{p_2,p_1,p_2}(z_1-z_2,z_1,z_2) 
\end{gathered}
\end{equation} 
in the regions 
\begin{equation}
\begin{gathered}
|z_1| > |z_2| > 0 \text{ and } -\frac{\pi}{2} < \arg(z_1 - z_2) - \arg(z_1)  < \frac{\pi}{2} ,\\
|z_2| > |z_1| > 0 \text{ and } - \frac{3 \pi }{2} < \arg(z_1 - z_2) - \arg(z_2) < -\frac{\pi}{2} ,
\end{gathered}
\end{equation} 
respectively. For the same fixed $ w_1 $, $w_2 $ and $ w_3' $, we pick some $ q \in \mathbb{Z} $ sufficiently high so that $ \res_{x_1} x_1^{q+l} Y_{W_2}(x_1^{\mathcal{L}_{g_2}} v,x_1) w_2 = 0 $ for all $ l \in \mathbb{Z}_{\geq 0} $ and $ v \in U $, by the lower truncation property of $ Y_{W_2} $. Then
\begin{align*}
0 &= \res_{x_1} x_1^{q} x_0^{-1} \delta \left( \frac{- x_2 + x_1}{x_0} \right)  \mathcal{Y}(w_1,x_2) Y_{W_2}\left(  (-x_2 + x_1)^{\mathcal{L}_{g_1}} x_1^{\mathcal{L}_{g_2}} u, x_1 \right) w_2 \\ 
&= \res_{x_1} x_1^{q} \left( x_0^{-1} \delta \left( \frac{x_1 - x_2}{x_0} \right) Y_{W_3}\left( (x_1 - x_2)^{\mathcal{L}_{g_1}} x_1^{\mathcal{L}_{g_2}} u,x_1 \right) \mathcal{Y}(w_1,x_2) w_2 \right.\\
& \qquad - \left. x_1^{-1} \delta \left( \frac{x_2 + x_0}{x_1} \right) \mathcal{Y}\left( Y_{W_1} \left( x_0^{\mathcal{L}_{g_1}} (x_2 + x_0)^{\mathcal{L}_{g_2}}u, x_0 \right)w_1,x_2 \right) w_2 \right) \\
&= \res_{x_1} x_1^{q} \left( x_1^{-1} \delta \left( \frac{x_0 + x_2}{x_1} \right) Y_{W_3}\left( (x_1 - x_2)^{\mathcal{L}_{g_1}} x_1^{\mathcal{L}_{g_2}} u,x_1 \right) \mathcal{Y}(w_1,x_2) w_2\right. \\
& \qquad - \left.x_1^{-1} \delta \left( \frac{x_2 + x_0}{x_1} \right) \mathcal{Y}\left( Y_{W_1} \left( x_0^{\mathcal{L}_{g_1}} (x_2 + x_0)^{\mathcal{L}_{g_2}}u, x_0 \right)w_1,x_2 \right) w_2 \right) \\
&= Y_{W_3}\left( \left((x_0 + x_2) - x_2 \right)^{\mathcal{L}_{g_1}} (x_0 + x_2)^{\mathcal{L}_{g_2}+q} u, x_0 + x_2 \right) \mathcal{Y}(w_1,x_2) w_2  \\
& \qquad -  \mathcal{Y}\left( Y_{W_1} \left( x_0^{\mathcal{L}_{g_1}} (x_2 + x_0)^{\mathcal{L}_{g_2} + q}u, x_0 \right)w_1,x_2 \right) w_2 .
\end{align*}
Hence, 
\begin{align*}
& \left\langle w_3', Y_{W_3}\left( \left((x_0 + x_2) - x_2 \right)^{\mathcal{L}_{g_1}} (x_0 + x_2)^{\mathcal{L}_{g_2}+q} u, x_0 + x_2 \right) \mathcal{Y}(w_1,x_2) w_2 \right\rangle \\
& \qquad = \left\langle w_3', \mathcal{Y}\left( Y_{W_1} \left( x_0^{\mathcal{L}_{g_1}} (x_2 + x_0)^{\mathcal{L}_{g_2} + q}u, x_0 \right)w_1,x_2 \right) w_2 \right\rangle .
\end{align*}
This series has no logarithms of $ x_0 $ and only integral powers in $ x_0 $. Furthermore, it is lower truncated in $ x_0 $ and $ x_2 $, and hence also upper truncated in $ x_0 $ and $ x_2 $ by \eqref{eq:L(0) power relation c5}. So, we see that it is equal to some series
\begin{align*}
h_2(u;x_0,x_2) = \sum_{i,j,k = 0}^M b_{ijk} x_0^{c_i} x_2^{d_j} (\log x_2)^k ,
\end{align*}
for some $ c_i \in \mathbb{Z} $ and $ b_{ijk}, d_j \in \mathbb{C} $, implicitly dependent on $ w_1 $, $ w_2 $ and $ w_3' $. Define $ f_2(x_0,x_1,x_2) =  h_2\left(x_0^{-\mathcal{L}_{g_1}} x_1^{-\mathcal{L}_{g_2} - q}u;x_0,x_2 \right) $, which we can do since the choice of $ q $ is uniform over $ U $. Then
\begin{equation}
\begin{aligned}
 \langle w_3', Y_{W_3}(u,x_1) \mathcal{Y}(w_1,x_2) w_2 \rangle |_{x_1 = z_1, x_2 = z_2}^{p_1,p_2} &= f_2(x_1 - x_2,x_1,x_2)|_{x_1 = z_1, x_2 = z_2}^{p_1,p_2}, \\
 \langle w_3',  \mathcal{Y}(Y_{W_1}(u,x_0)w_1,x_2) w_2 \rangle |_{x_0 = z_1-z_2, x_2 = z_2}^{p_{12},p_2} &= f_2(x_0,x_2 + x_0,x_2)|_{x_0 = z_1-z_2, x_2 = z_2}^{p_{12},p_2} ,
 \end{aligned}
\end{equation}
are absolutely convergent in the regions $ |z_1| > |z_2| > 0 $ and $ |z_2| > |z_1 - z_2| > 0 $, respectively. Moreover, they converge to 
\begin{equation}
\begin{gathered}
f_2^{p_1,p_1,p_2}(z_1 - z_2,z_1,z_2) \quad \text{and} \quad f_2^{p_{12}, p_2, p_2}(z_1-z_2,z_1,z_2) 
\end{gathered}
\end{equation} 
in the regions 
\begin{equation}
\begin{gathered}
|z_1| > |z_2| > 0 \text{ and } -\frac{\pi}{2} < \arg(z_1 - z_2) - \arg(z_1)  < \frac{\pi}{2} ,\\
|z_2| > |z_1 - z_2| > 0 \text{ and } -\frac{\pi}{2} < \arg(z_1) - \arg(z_2)  < \frac{\pi}{2} , 
\end{gathered}
\end{equation} 
respectively. 

Set $ f(x_0,x_1,x_2) = f_1(x_0,x_1,x_2) $ to be the series in the duality property. By looking at the convergence of $ \langle w_3', Y_{W_3}(u,x_1) \mathcal{Y}(w_1,x_2) w_2 \rangle |_{x_1 = z_1, x_2 = z_2}^{p_1,p_2} $, we see that $ f_1^{p_1,p_1,p_2}(z_1 - z_2,z_1,z_2) = f_2^{p_1,p_1,p_2}(z_1 - z_2,z_1,z_2) $ in an open region.  Since these are complex analytic functions in $ z_1 $ and $ z_2 $, we have $ f_1^{p_{12},p_1,p_2}(z_1 - z_2,z_1,z_2) = f_2^{p_{12},p_1,p_2}(z_1 - z_2,z_1,z_2) $. Thus, the series \eqref{eq:multiseries} converge to \eqref{eq:branches} on the regions \eqref{eq:regions}.
\end{proof}

Thus, when a twisted intertwining operator is among modules twisted by commuting automorphisms, the Jacobi identity \eqref{eq:twisted jacobi c5} is \textit{equivalent} to the duality property appearing in the definition of twisted intertwining operators. By replacing the duality property with the Jacobi identity, we obtain a completely algebraic definition for twisted intertwining operators among modules twisted by commuting automorphisms. While the original analytic definition in \cite{HUANG2018346} has an insightfully transparent motivation (which reveals deep connections to orbifold conformal field theory), a completely algebraic definition is sometimes more suited for certain proofs such as the results in this paper and \cite{HuangNahmInequality}.

\section{The setting for deriving differential equations} \label{sec:setting c5}
Let $ V $ be a vertex operator algebra and let $ G $ be an \emph{abelian} subgroup of $ \Aut V  $ (with possibly infinitely many elements). The Jacobi identity \eqref{eq:twisted jacobi c5} is sufficient to fully describe the intertwining operators among modules twisted by elements in $ G $. In what follows, our aim is to explore how the methods of \cite{HuangDEs}, \cite{HuangC1} and \cite{TanMostlyUntwisted} can generalize to this setting. We will see that the logarithms of the eigenvalues for $ g \in G $ play an important role since they shift the weights of elements when performing residue calculus with the Jacobi identity. To compensate for this shifting, the assumption of $ C_\n $-cofiniteness will be used, as opposed to just $ C_1 $-cofiniteness as in \cite{HuangDEs}, \cite{HuangC1} and \cite{TanMostlyUntwisted}.

We recall the definition of $ C_n $-cofiniteness for twisted modules as introduced in \cite{TanCofinite}. First, recall that twisted vertex operators for twisted modules $ (W,Y_W) $ have the form $ Y_W(u,x) = \sum_{n \in \mathbb{C}} \sum_{k = 0}^K u_{n,k} x^{-n-1} (\log x)^k $. We write $ (Y_W)_0(u,x) $ for the series $ \sum_{n \in \mathbb{C}} u_{n,0} x^{-n-1} $ of $ (\log x) $-free terms. 

\begin{definition} \label{def:Cn-cofinite twisted modules}
Let $ W $ be a $ g $-twisted $ V $-module. For each $ n \in \mathbb{Z}_{>0} $, we define $ C_n(W) $ to be the subspace of $ W $ spanned by the elements
\begin{equation}
u_{\alpha -  n ,0} w , \qquad \text{for all } u \in V^{[\alpha]}, \  \alpha \in \{z \in \mathbb{C} : 0 \leq \Re(z) < 1\}, \ w \in W .
\end{equation}
When $ n = 1 $, we impose the extra condition that $ u \in V_+ = \coprod_{m \in \mathbb{Z}_{> 0}} V_{(m)} $.  We say that $ W $ is \emph{$ C_n $-cofinite} if $ \dim W / C_n(W) < \infty $.
\end{definition}

Also recall the notion of a discretely graded (i.e., quasi-finite) vector space. 

\begin{definition}
Let $ W = \coprod_{n \in \mathbb{C}} W_{[n]} $ be a $ \mathbb{C} $-graded vector space. We say that $ W $ is \emph{discretely graded} (also known as \emph{quasi-finite dimensional}) if 
\begin{equation}\label{eq:discrete graded}
\dim \coprod_{\substack{ n \in \mathbb{C} \\ \Re(n) < r}} W_{[n]} < \infty , \qquad \text{ for all } r \in \mathbb{R}.
\end{equation} 
\end{definition}

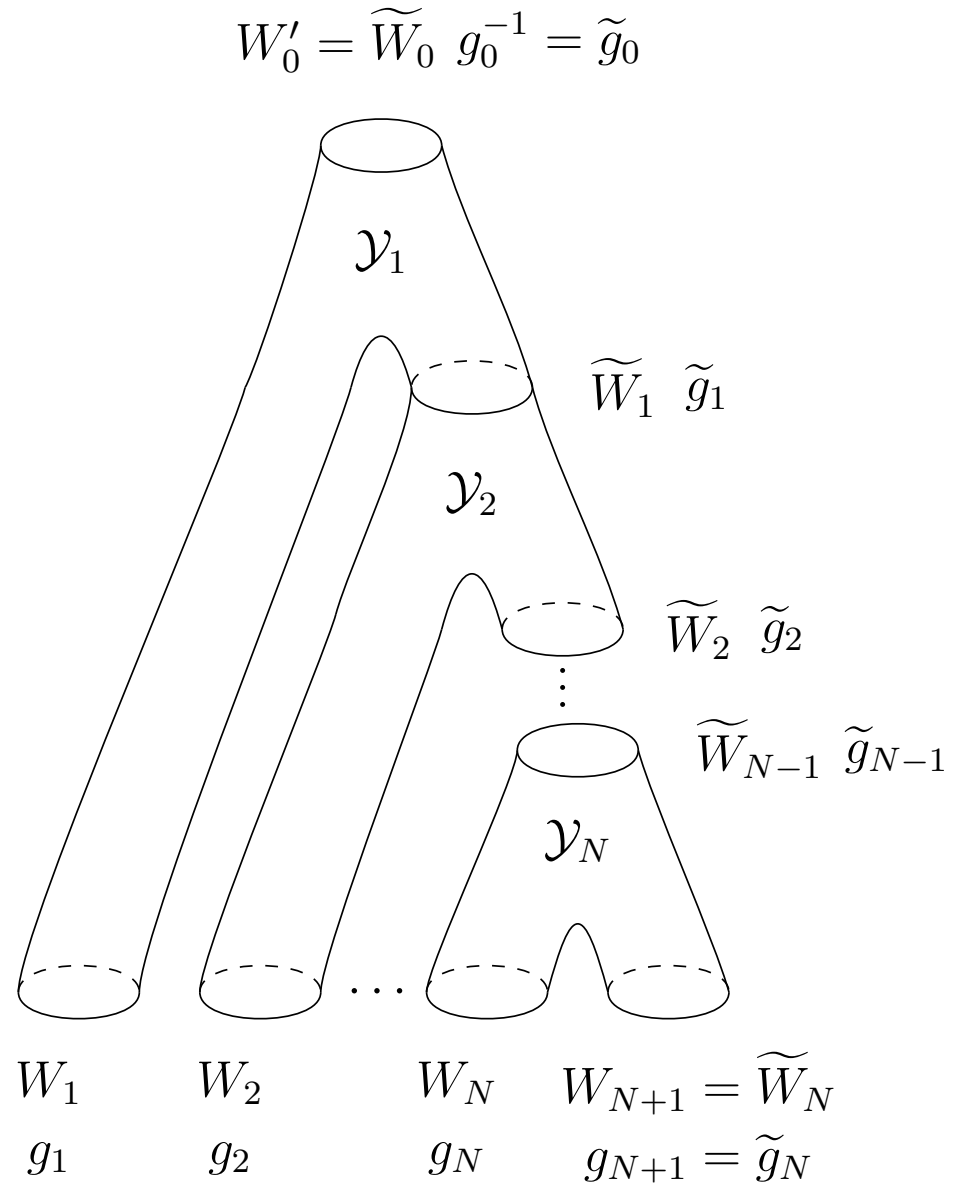
\begin{figure}[b!]
\centering

\scalebox{0.8}{
\begin{tikzpicture}[scale=0.5]
	\begin{pgfonlayer}{nodelayer}
		\node [style=none] (0) at (-5, 0) {};
		\node [style=none] (1) at (-3, 0) {};
		\node [style=none] (2) at (-2, 0) {};
		\node [style=none] (3) at (0, 0) {};
		\node [style=none] (4) at (1.75, 0) {};
		\node [style=none] (5) at (3.75, 0) {};
		\node [style=none] (6) at (4.75, 0) {};
		\node [style=none] (7) at (6.75, 0) {};
		\node [style=none] (8) at (3.25, 4) {};
		\node [style=none] (9) at (5.25, 4) {};
		\node [style=none] (19) at (2, 6) {};
		\node [style=none] (20) at (0.25, 6.25) {};
		\node [style=none] (21) at (2, 6) {};
		\node [style=none] (22) at (3, 6) {};
		\node [style=none] (23) at (5, 6) {};
		\node [style=none] (24) at (1.5, 10) {};
		\node [style=none] (25) at (3.5, 10) {};
		\node [style=none] (26) at (-1.25, 10) {};
		\node [style=none] (27) at (0.5, 10) {};
		\node [style=none] (28) at (1.5, 10) {};
		\node [style=none] (29) at (3.5, 10) {};
		\node [style=none] (30) at (0, 14) {};
		\node [style=none] (31) at (2, 14) {};
		\node [style=none] (32) at (-4.5, -1.5) {$W_1$};
		\node [style=none] (33) at (-1.5, -1.5) {$W_2$};
		\node [style=none] (34) at (2.25, -1.5) {$W_\N$};
		\node [style=none] (35) at (6.25, -1.5) {$W_{\N+1} = \W_{\N}$};
		\node [style=none] (36) at (-4.5, -2.75) {$g_1$};
		\node [style=none] (37) at (-1.5, -2.75) {$g_2$};
		\node [style=none] (38) at (2.25, -2.75) {$g_{\N}$};
		\node [style=none] (39) at (6.25, -2.75) {$g_{\N+1} = \g_{\N}$};
		\node [style=none] (40) at (7.25, 4) {$\W_{\N-1}$};
		\node [style=none] (41) at (6.25, 6) {$\W_2$};
		\node [style=none] (42) at (5, 10) {$\W_1$};
		\node [style=none] (43) at (0.25, 15.75) {$W_0' = \W_0$};
		\node [style=none] (44) at (9.25, 4) {\ \ $\g_{\N-1}$};
		\node [style=none] (45) at (7.5, 6) {\ $\g_2$};
		\node [style=none] (46) at (6.25, 10) {\ $\g_1$};
		\node [style=none] (47) at (3.5, 15.75) {\ \ $g_0^{-1} = \g_0$};
		\node [style=none] (48) at (1, 0) {$\cdots$};
		\node [style=none] (49) at (4, 5.25) {$\vdots$};
		\node [style=none] (50) at (1, 12.25) {$\mathcal{Y}_1$};
		\node [style=none] (51) at (2.5, 8.25) {$\mathcal{Y}_2$};
		\node [style=none] (52) at (4.25, 2.5) {$\mathcal{Y}_{\N}$};
	\end{pgfonlayer}
	\begin{pgfonlayer}{edgelayer}
		\draw [style=plainedge, bend right=90, looseness=0.75] (6.center) to (7.center);
		\draw [style=plainedge, bend right=90, looseness=0.75] (4.center) to (5.center);
		\draw [style=plainedge, bend right=90, looseness=0.75] (2.center) to (3.center);
		\draw [style=plainedge, bend right=90, looseness=0.75] (0.center) to (1.center);
		\draw [style=dotted line, bend left=90, looseness=0.75] (6.center) to (7.center);
		\draw [style=dotted line, bend left=90, looseness=0.75] (4.center) to (5.center);
		\draw [style=dotted line, bend left=90, looseness=0.75] (2.center) to (3.center);
		\draw [style=dotted line, bend left=90, looseness=0.75] (0.center) to (1.center);
		\draw [style=plainedge, bend right=90, looseness=0.75] (8.center) to (9.center);
		\draw [style=plainedge, bend left=270, looseness=0.75] (9.center) to (8.center);
		\draw [style=plainedge, in=105, out=-75, looseness=0.50] (9.center) to (7.center);
		\draw [style=plainedge, in=105, out=75, looseness=4.00] (5.center) to (6.center);
		\draw [style=plainedge, in=-90, out=90, looseness=0.25] (4.center) to (8.center);
		\draw [style=plainedge, in=90, out=-105, looseness=0.25] (19.center) to (3.center);
		\draw [style=plainedge, bend right=90, looseness=0.75] (22.center) to (23.center);
		\draw [style=dotted line, bend left=90, looseness=0.75] (22.center) to (23.center);
		\draw [style=plainedge, bend right=90, looseness=0.75] (24.center) to (25.center);
		\draw [style=dotted line, bend left=270, looseness=0.75] (25.center) to (24.center);
		\draw [style=plainedge, in=105, out=-75, looseness=0.75] (25.center) to (23.center);
		\draw [style=plainedge, in=105, out=75, looseness=3.25] (21.center) to (22.center);
		\draw [style=plainedge, in=-90, out=75, looseness=0.50] (20.center) to (24.center);
		\draw [style=plainedge, bend right=90, looseness=0.75] (30.center) to (31.center);
		\draw [style=plainedge, in=105, out=-75, looseness=0.75] (31.center) to (29.center);
		\draw [style=plainedge, in=105, out=75, looseness=3.00] (27.center) to (28.center);
		\draw [style=plainedge, in=-90, out=60, looseness=0.25] (26.center) to (30.center);
		\draw [style=plainedge, in=90, out=-105, looseness=0.25] (20.center) to (2.center);
		\draw [style=plainedge, in=-105, out=90, looseness=0.25] (1.center) to (27.center);
		\draw [style=plainedge, in=-105, out=90, looseness=0.25] (0.center) to (26.center);
		\draw [style=plainedge, bend left=90, looseness=0.75] (30.center) to (31.center);
	\end{pgfonlayer}
\end{tikzpicture}
}

\caption{A diagram summarizing the twisted $ V $-modules, the automorphism in the abelian group $ G $ twisting each module, and the intertwining operators defining the correlation functions. The modules $ W_1, \dots, W_{\N+1} $ contain the \emph{in-states} of the correlation functions. The module $ \W_0 $ contains the \emph{out-states}. The modules $ \W_1, \dots, \W_{\N-1} $ contain the \emph{intermediate} states.} \label{fig:diagram c5}
\end{figure}

We now provide the setting for the following sections of this paper. We have a fixed vertex operator algebra $ V $ and an abelian subgroup $ G \subseteq \Aut V  $. Let $ \N \in \mathbb{Z}_{>0} $ be the number of twisted intertwining operators being considered. Let $ g_0, \dots, g_{\N+1} , \g_0, \dots, \g_N \in G $. Assume that $ g_0^{-1} = \g_0 $, $ g_{\N+1} = \g_\N $, and $ g_i \g_i = \g_{i-1} $ for $ i = 1, \dots, \N $. Let $ W_i $ be $ g_i $-twisted $ V $-modules, for $ i = 1, \dots, \N $. Let $ \widetilde{W}_i $ be $ \g_i $-twisted $ V $-modules, for $ i = 0, \dots, \N $. Let $ \mathcal{Y}_i : W_i \otimes \widetilde{W}_i \to  \widetilde{W}_{i-1} \{x\} [ \log x ] $ be twisted intertwining operators, for $ i = 1, \dots , \N $. To aid with notation later on, we also denote $ \widetilde{W}_0' $ by $ W_0 $ (a $ g_0 $-twisted $ V $-module), and denote $ \widetilde{W}_{\N} $ by $ W_{\N+1} $. This data is summarized in Figure~\ref{fig:diagram c5}. 

This completes the specification of the data needed to define chiral genus-zero correlation functions in a $ G $-equivariant conformal field theory. The conditions on the data that we need to obtain our results are stated as follows.

Recall that $ G $ is assumed to be abelian, but may be infinite, and there are no restrictions such as semisimplicity of the elements. The modules $ W_0, \dots, W_{\N+1}, \W_1, \dots ,\W_{N-1} $ are assumed to be discretely graded. The \emph{in-state} modules $ W_1, \dots, W_{\N+1} $ and the \emph{intermediate} modules $ \W_1, \dots ,\W_{N-1} $ have the following cofiniteness conditions. Let $ \n \in \mathbb{Z}_{>0} $ be a parameter to control the strength of the $ C_\n $-cofiniteness condition. For each $ i = 1, \dots, \N +1 $, we assume that 
\begin{equation} \label{eq:cofinite conditions 1}
\begin{aligned}
&\circ \ \text{ $ W_i $ is $ C_{\n} $-cofinite if $ g_i $ has an eigenvalue not in $ \mathbb{R}_{>0} $; and  } \\
&\circ \ \text{ $ W_i $ is $ C_1 $-cofinite if $ g_i $ only has eigenvalues in $ \mathbb{R}_{>0} $. }
\end{aligned}
\end{equation}

Similarly, for each $ i = 1, \dots, \N -1 $, we assume that 
\begin{equation} \label{eq:cofinite conditions 2}
\begin{aligned}
&\circ \ \text{ $ \W_i $ is $ C_{\n} $-cofinite if $ \g_i $ has an eigenvalue not in $ \mathbb{R}_{>0} $; and  } \\
&\circ \ \text{ $ \W_i $ is $ C_1 $-cofinite if $ \g_i $ only has eigenvalues in $ \mathbb{R}_{>0} $. }
\end{aligned}
\end{equation}

We will assume that $ V $ is of CFT type, i.e. $ V_{(m)} = 0 $ when $ m < 0 $ and $ V_{(0)} = \mathbb{C} \one $. 

Many calculations below will not require the cofiniteness conditions, discrete grading conditions, nor the CFT type conditions. We will make it clear when these conditions are needed to be assumed.

\begin{remark}
Note that there is no cofiniteness assumption for $ W_0 $ or $ \W_0 $, since we will use the techniques discovered in \cite{HuangC1}. We have taken the idea and method used to remove the cofiniteness conditions for $ W_0 $ or $ \W_0 $ directly from Huang. 
\end{remark}

\begin{remark}
When an automorphism $ g $ has finite order, the statement ``$ g $ only has eigenvalues in $ \mathbb{R}_{>0} $" becomes equivalent to the statement ``$ g = \id_V $". We will use this in Section~\ref{sec:convergence} when we assume that $ G $ is finite. Until Section~\ref{sec:convergence}, we will not assume that automorphisms have finite order, so the former statement still remains more general. 
\end{remark}

\begin{remark}
The conditions for when a module is assumed to be $ C_{\n} $-cofinite or $ C_1 $-cofinite may at first seem strange, so we will present some insight for now before it becomes clear in Section~\ref{sec:differential equations}. If an automorphism $ g $ of $ V $ has an eigenvalue $ e^{2 \pi \iu \alpha }$ corresponding to some generalized eigenvector $ u \in V $, it is in fact the normalized logarithm $ \alpha $ (when $ \alpha $ is chosen to have real part in $ [0,1) $) that plays a role in shifting the weights of the action of $ u $ on a $ g $-twisted module. This is due to the fact that the action of such a generalized eigenvector $ u $ will always have the form
\begin{align*}
Y_W(u,x) = \sum_{m \in \alpha + \mathbb{Z}} \sum_{k = 0}^K u_{m,k} x^{-m-1} (\log x)^k,
\end{align*}
recalling that $ \wt( u_{m,k}) = \wt u - m - 1 $. In Section~\ref{sec:differential equations}, only the real part of such a shift will affect arguments that classically employ $ C_1 $-cofiniteness conditions, and hence introduce the need for $ C_{\n} $-cofiniteness conditions. Since no such shift of the weights occurs when $ \Re(\alpha) = 0 $, this corresponds to $ g $-eigenvalues of the form $ e^{- 2 \pi \Im(\alpha) } \in \mathbb{R}_{>0} $.
\end{remark}

\begin{remark}
We will impose cofiniteness and discrete grading assumptions for the intermediate modules $ \W_1, \dots, \W_{\N-1} $ to be used in Section~\ref{sec:convergence} when proving convergence of the series $ \langle w_0, \mathcal{Y}_1(w_1,x_1) \cdots \mathcal{Y}_\N(w_\N,x_\N) w_{\N+1} \rangle $. If all in-state modules $ W_1, \dots, W_{\N+1} $ are $ C_n $-cofinite, for $ n \in \mathbb{Z}_{\geq 2} $, then the twisted Nahm inequality from \cite{HuangNahmInequality} provides $ C_n $-cofiniteness of the intermediate modules (when $ \W_{i-1} $ is a \textit{surjective product} of $ W_i $ and $ \W_i $). At the present time there is no analogous twisted Nahm inequality involving $ C_1 $-cofinite twisted modules. It will be interesting to see if such an inequality can be proven, even if only in the case that $ C_1 $-cofinite modules are twisted by elements with only positive real eigenvalues.
\end{remark}

Let $ A $ be the additive subgroup of $ \mathbb{C} $ generated by the union of $ \{1\} $ and the set of eigenvalues for $ 
\{ \mathcal{S}_g \}_{g \in G} $. Let $ R $ be the subring of 
\begin{equation*}
\mathbb{C}\{ x_i , (x_j - x_k) : i, j,k = 1, \dots, \N \text{ with } j < k \} [ \log x_i, \log (x_j - x_k)  : i,j,k = 1, \dots, \N \text{ with } j < k]
\end{equation*}
generated by $ \{ x_i^\alpha, (x_j - x_k)^\alpha , \log x_i , \log (x_j - x_k) : i,j,k = 1, \dots, \N \text{ with } j < k ; \alpha \in A \} $. Observe that to compress notation, we can define $ x_{\N+1} = 0 $ and write the generators of $ R $ as
\begin{equation*}
\{ (x_i - x_j)^\alpha , \log (x_i - x_j) : i < j \in  \{1, \dots, \N +1\} ; \alpha \in A \},
\end{equation*}
as we will do later on. We will write the elements of $ R $ in function notation as $ f(x_1, \dots, x_\N) $ with explicit dependence on $ x_1, \dots, x_\N $ only. 

In the case that each $ g \in G $ acts semisimply on $ V $, we can remove all of the $ \log x_i $ and $ \log (x_j - x_k) $ from $ R $. This occurs in the important case that each $ g \in G $ has finite order (see Section~\ref{sec:convergence}). 

We have the map
\begin{equation}
\iota_{|x_1| > \dots > |x_\N|}: R \to \mathbb{C} \{x_1, \dots, x_\N\}[\log x_1, \dots, \log x_\N]  
\end{equation}
that expands $ (x_j-x_k)^{\alpha} $ and $ \log (x_j - x_k) $ in non-negative powers of $ x_k $ when $ j < k $. 

Note that $ R $ is Noetherian if $ A $ is finitely generated. In the case that $ V $ and $ G $ are both finitely generated, we can project the finitely many generators onto the subspaces $ V_{(n)} $ and consider the finitely many subspaces $ V_{(n)} $ where the image of the projection is non-trivial. Each generator $ g $ for $ G $ decomposes these $ V_{(n)} $ into eigenspaces for $ \mathcal{S}_g $. We can then take a basis of eigenvectors for $ \mathcal{S}_g $ in each of these subspaces to obtain a finite generating set for $ V $ consisting of eigenvectors for $ \mathcal{S}_g $. Since we can choose finitely many generators for $ G $, $ A $ is finitely generated when $ V $ and $ G $ are both finitely generated. 

Regardless of whether $ V $ is finitely generated, when $ G $ is finite, then $ A = \frac{1}{t_0} \mathbb{Z} $ for some $ t_0 \in \mathbb{Z}_{>0} $ dividing $ |G| $. This is the important case when proving convergence (see Section~\ref{sec:convergence}).

Define
\begin{equation}
T = R \otimes W_0 \otimes \cdots \otimes W_{\N+1} ,
\end{equation}
which is to be viewed as an $ R $-module by multiplication on the left. The $ \mathbb{C} $-gradings by conformal weight on $ W_i $, together with the trivial grading on $ R $, induce a $ \mathbb{C} $-grading on $ T $ called the \emph{weight}. We denote the subspace of weight-homogeneous elements of weight with real part $ r \in \mathbb{R} $ by $ T_{(r)} $ so that $ T = \coprod_{r \in \mathbb{R}} T_{(r)} $. Observe that the discrete gradings on the modules $ W_i $ ensure that each subspace $ F_r(T) = \coprod_{s \leq r} T_{(s)} $, $ r\in \mathbb{R} $, is a finitely generated $ R $-module. Hence, we have a filtration $ F_r(T) \subseteq F_s(T) $, for $ r \leq s \in \mathbb{R} $, of finitely generated submodules of $ T $. 

Later, we will see that we can work with the submodules 
\begin{equation}
T_{K_0} = R \otimes \left(  \coprod_{\substack{ m \in \mathbb{C} \\ \Re(m) \leq K_0}} (W_0)_{[m]} \right) \otimes W_1 \otimes \cdots \otimes W_{\N+1}
\end{equation}
of $ T $, for each $ K_0 \in \mathbb{R} $. There are similar gradings and filtrations on $ T_{K_0} $.

We are interested in the series
\begin{equation}\label{eq:correlator c5}
\langle w_0, \mathcal{Y}_1(w_1,x_1) \cdots \mathcal{Y}_\N(w_\N,x_\N) w_{\N+1} \rangle, 
\end{equation}
for each $(\N+2) $-tuple of $ w_0 \in W_0, \dots, w_{\N+1} \in W_{\N+1} $. Consider the linear map
\begin{equation}  \label{eq:correlation map}
\begin{gathered}
\phi : T \to \mathbb{C}[\log x_1, \dots, \log x_\N] \{x_1, \dots, x_\N \}, \\
\begin{aligned}
&\phi ( f(x_1,\dots, x_\N) w_0 \otimes \cdots \otimes w_{\N+1} )   \\
&\qquad = \iota_{|x_1| > \cdots > |x_\N|} f(x_1,\dots, x_\N) \langle w_0, \mathcal{Y}_1(w_1,x_1) \cdots \mathcal{Y}_\N(w_\N,x_\N) w_{\N+1} \rangle
\end{aligned}
\end{gathered}
\end{equation}
that maps elements of $ T $ to chiral correlation functions defined by the product of $ \mathcal{Y}_1 , \dots, \mathcal{Y}_\N $ formally expanded in the domain $ |x_1| > \cdots > |x_\N| > 0 $. This map of course restricts to $ T_{K_0} $, for each $ K_0 \in \mathbb{R} $.

We seek to find formal differential equations satisfied by \eqref{eq:correlator c5} in this formal region. Furthermore, when $ G $ is assumed to be finite, these differential equations will have regular singularities at certain points.

\section{Computations} \label{sec:computations}

In this section, we compute the generators $ \mathcal{A}_\el $ for a certain submodule of $ T_{K_0} $ that will be used to derive the desired differential equations. The algebraic techniques with commutators of $ Y^+ $ and $ Y^- $ used in \cite{TanMostlyUntwisted} do not extend beyond \gtype\ intertwining operators. So, a different method will be used in this paper.

In what follows, the subscripts of $ Y_{W_i} $ and $ \mathcal{Y}_i $ will be dropped for brevity. The subscripts can be deduced by determining what module $ Y $ is acting on and by determining what module the argument of $ \mathcal{Y} $ lies in. Let $ u \in V $ be a simultaneous generalized eigenvector for the elements in $ G $. Let $ \alpha_i \in \mathbb{C} $ with $ \Re(\alpha_i) \in [0,1) $ be such that $ \mathcal{S}_{g_i} u = \alpha_i u $. Let $ w_i \in W_i $. 

\subsection{A heuristic argument}
We first give a heuristic argument for the main result of this section before we proceed with the rigorous computations.
Assume that all products of intertwining operators converge and analytically continue to the maximally defined regions. (Note that we cannot logically assume this convergence and analytic continuation since it is what we are aiming to prove.)

Let $ |z_1| > \cdots > |z_\N| > 0 $. Let $ \xi_0, \dots, \xi_{\N+1} $ be formal variables. Let $ \zeta $ be a complex variable. We consider the multivalued function 
\begin{align*}
& f(w_1, \dots, w_\N; u ; z_1, \dots, z_\N; \zeta) \\
	&\quad= F \left( \langle w_0, Y(u,\xi_0) \mathcal{Y}(w_1, x_1) \cdots \mathcal{Y}(w_\N, x_\N) w_{\N+1} \rangle|_{x_i = z_i,\, \xi_0 = \zeta} \right) \\
	&\quad= F\left( \langle w_0, \mathcal{Y}(Y(u,\xi_1 ) w_1, x_1) \cdots \mathcal{Y}(w_\N, x_\N) w_{\N+1} \rangle|_{x_i = z_i,\, \xi_1 = \zeta - z_1 } \right)\\
	&\quad= \cdots \\
	&\quad= F \left( \langle w_0, \mathcal{Y}( w_1, x_1) \cdots \mathcal{Y}(Y(u,\xi_\N)w_\N, x_\N) w_{\N+1} \rangle|_{x_i = z_i, \,  \xi_\N = \zeta - z_N} \right)\\
	&\quad= F \left( \langle w_0, \mathcal{Y}( w_1, x_1) \cdots \mathcal{Y}(w_\N, x_\N) Y(u,\xi_{\N+1}) w_{\N+1} \rangle|_{x_i = z_i, \,  \xi_{\N+1} = \zeta } \right),
\end{align*}
where $ F $ denotes the analytic continuation of a function.

The function $ f(w_1, \dots, w_\N; u ; z_1, \dots, z_\N; \zeta) $ will not necessarily be single-valued as a function in $ \zeta $ since we are working with twisted modules. (Compare this to the case where we are working with only untwisted modules. Then $ f $ is immediately single-valued in $ \zeta $.) To perform contour integration, we seek a single-valued function. 

The key observation is the following: To obtain a single-valued function (as a function in $ \zeta $) we can replace $ u $ with $ (\zeta - z_1)^{\mathcal{L}_{g_1}} \cdots (\zeta - z_\N)^{\mathcal{L}_{g_\N}} \zeta^{\mathcal{L}_{g_{\N+1}}} u $. For example, when we circle around infinity, we have $ |\zeta | > |z_1| > \cdots > |z_\N| $, so $ Y( (\zeta - z_1)^{\mathcal{L}_{g_1}} \cdots (\zeta - z_\N)^{\mathcal{L}_{g_\N}} \zeta^{\mathcal{L}_{g_{\N+1}}} u , \zeta ) $ expands to a sum of $ Y( \zeta^{\mathcal{L}_{g_1} + \dots + \mathcal{L}_{g_{\N+1}}} v , \zeta ) $ multiplied by a power series in $ \zeta^{-1} $ (where $ v $ is an element of $ V $), which is single-valued. Also, when we follow an arbitrarily small circle around each $ z_i $, then $ |z_i|, |z_i - z_j| > |\zeta - z_i| $. So, $ Y( (\zeta - z_1)^{\mathcal{L}_{g_1}} \cdots (\zeta - z_\N)^{\mathcal{L}_{g_\N}} \zeta^{\mathcal{L}_{g_{\N+1}}}  u , \zeta - z_i) $ expands to a sum of $ Y( (\zeta - z_i)^{\mathcal{L}_{g_i}} v, \zeta - z_i) $ multiplied by a power series in $ (\zeta - z_i) $ (where $ v $ is an element of $ V $), which is single-valued. Finally, we can further observe that the same is true when we set $ z_{\N+1} = 0 $, and we follow an arbitrarily small circle around the origin. 

Let $ n_1, \dots, n_{\N+1} \in \mathbb{Z} $, and let $ z_{\N+1} = 0 $ for convenience. For fixed $ w_1, \dots, w_\N , u, z_1, \dots, z_\N $, we define
\begin{align*}
h(\zeta) =  f \left( w_1, \dots, w_\N; \prod_{j=1}^{\N+1} (\zeta - z_j)^{\mathcal{L}_{g_j} - n_j} u ; z_1, \dots, z_\N; \zeta \right).
\end{align*}
Observe that $ h(\zeta) $ is a  single-valued function of $ \zeta $, so we can integrate it along closed contours and apply the residue theorem. Contour integration gives
\begin{align*}
\oint_{|\zeta| = r } h(\zeta) d \zeta = \sum_{i = 1}^{\N+1} \oint_{|\zeta-z_i| = \epsilon_i} h(\zeta) d \zeta,
\end{align*}
 for some $ r > |z_1| $ and sufficiently small $ \epsilon_i > 0 $. Hence,
\begin{align*}
& \oint_{|\zeta| = r } F \left( \langle w_0, Y \left( \prod_{j=1}^{\N+1}(\zeta - z_j)^{\mathcal{L}_{g_j} - n_j} u,\xi_0 \right) \mathcal{Y}(w_1, x_1) \cdots \mathcal{Y}(w_\N, x_\N) w_{\N+1} \rangle|_{x_k = z_k, \xi_0 = \zeta} \right) d \zeta \\
&= \sum_{i = 1}^{\N} \oint_{|\zeta-z_i| = \epsilon_i} F \left( \langle w_0,\mathcal{Y}(w_1, x_1) \cdots \mathcal{Y}\left( Y\left( \prod_{j=1}^{\N+1}((z_i - z_j) + (\zeta - z_i))^{\mathcal{L}_{g_j} - n_j}    u,\xi_i \right)w_i, x_i\right) \cdots \right. \\
&\quad  \left. \vphantom{\prod_{j=1}^{\N+1}} \mathcal{Y}(w_\N, x_\N) w_{\N+1} \rangle|_{x_k = z_k, \, \xi_i = \zeta - z_i} \right)  d (\zeta-z_i) \\
& + \oint_{|\zeta| = \epsilon_{\N+1}} F \left( \langle w_0,\mathcal{Y}(w_1, x_1)  \cdots  \mathcal{Y}(w_\N, x_\N)  Y\left( \prod_{j=1}^{\N+1}(\zeta - z_j )^{\mathcal{L}_{g_j} - n_j}    u,\xi_{\N+1} \right) w_{\N+1} \rangle|_{x_k = z_k, \, \xi_{\N+1}  = \zeta} \right)  d \zeta, 
\end{align*}
which implies
\begin{align*}
& \sum_{k_1, \dots, k_\N = 0}^\infty F \left( (-x_1)^{k_1} \cdots (-x_\N)^{k_\N}  \langle w_0, \left( \prod_{j=1}^{\N} \binom{\mathcal{L}_{g_j} - n_j}{k_j} u \right)_{\alpha_1 + \cdots + \alpha_{\N+1} - n_1 - \cdots - n_{\N+1} - k_1 - \cdots - k_\N, 0} \right. \\
&\quad \left. \vphantom{ \binom{\mathcal{L}_{g_1} - n_1}{k_1} } \mathcal{Y}(w_1, x_1) \cdots \mathcal{Y}(w_\N, x_\N) w_{\N+1} \rangle|_{x_i = z_i} \right) \\
&= \sum_{i = 1}^{\N} \sum_{k_1, \dots, \widehat{k_i}, \dots, k_{\N+1} = 0}^\infty F \left( \langle w_0,\mathcal{Y}(w_1, x_1) \cdots \right. \\
& \quad \left. \mathcal{Y}\left( \left( \prod_{\substack{j = 1 \\ j \neq i} }^{\N + 1} \binom{\mathcal{L}_{g_j} - n_j}{k_j} (z_i-z_j)^{\mathcal{L}_{g_j} - n_j-k_j} u \right)_{\alpha_i - n_i + k_1 + \cdots + \widehat{k_i} + \cdots + k_{\N+1},0} w_i, x_i \right) \cdots \mathcal{Y}(w_\N, x_\N) w_{\N+1} \rangle|_{x_i = z_i} \right) \\
& \quad + \sum_{k_1, \dots, k_{\N} = 0}^\infty F \left( \langle w_0,\mathcal{Y}(w_1, x_1) \cdots \mathcal{Y}(w_\N, x_\N) \right. \\
& \qquad \left. \cdot \left( \prod_{j = 1}^{\N} \binom{\mathcal{L}_{g_j} - n_j}{k_j} (-z_j)^{\mathcal{L}_{g_j} - n_j-k_j} u \right)_{\alpha_{\N+1} - n_{\N+1} + k_1 + \cdots + k_{\N},0} w_{\N+1} \rangle|_{x_i = z_i} \right) .
\end{align*}
The notation $ \widehat{k_i} $ means that we omit $ k_i $ from the sum. Converting this identity to formal variables shows that the following element is in the kernel of $ \phi $:
\begin{equation} \label{eq:correlator identity c5}
\begin{aligned}
& - \sum_{k_1, \dots, k_\N = 0}^\infty (-x_1)^{k_1} \cdots (-x_\N)^{k_\N}  \left(  \prod_{j=1}^{\N} \binom{\mathcal{L}_{g_j} - n_j}{k_j} u \right)_{\alpha_1 + \cdots + \alpha_{\N+1} - n_1 - \cdots - n_{\N+1} - k_1 - \cdots - k_\N, 0}^* w_0 \\
&\qquad \otimes w_1 \otimes \cdots \otimes w_{\N+1} \\
& + \sum_{i = 1}^{\N+1} \sum_{k_1, \dots, \widehat{k_i}, \dots, k_{\N+1} = 0}^\infty  w_0 \otimes w_1 \otimes \cdots \otimes  w_{i-1} \\
& \quad \otimes \left( \prod_{\substack{j = 1 \\ j \neq i} }^{\N + 1} \binom{\mathcal{L}_{g_j} - n_j}{k_j} \cln(x_i-x_j)^{\mathcal{L}_{g_j} - n_j-k_j}\cln u \right)_{\alpha_i - n_i + k_1 + \cdots + \widehat{k_i} + \cdots + k_{\N+1},0} w_i \otimes w_{i+1} \otimes \cdots \otimes w_{\N+1},
\end{aligned}
\end{equation}

where $ x_{\N+1} = 0 $, $ -^* $ denotes the adjoint of an operator, and we have used the notation
\begin{align*}
\cln(x_i-x_j)^{\mathcal{A}}\cln &= \begin{cases}
(x_i-x_j)^{\mathcal{A}} = \sum_{k = 0}^\infty \binom{\mathcal{A}}{k} x_i^{\mathcal{A}-k} (-x_j)^k  &\text{if } i < j, \\
(-x_j+x_i)^{\mathcal{A}} = \sum_{k = 0}^\infty \binom{\mathcal{A}}{k} (-x_j)^{\mathcal{A}-k} x_i^k  &\text{if } j < i,
\end{cases} 
\end{align*}
for operators $ \mathcal{A} $ that are sums of semisimple and locally nilpotent parts. 

We emphasize that equation \eqref{eq:correlator identity c5} was derived above heuristically to provide intuition for the reader. We must prove equation \eqref{eq:correlator identity c5} through other means since we cannot assume that products of intertwining operators converge without proof. The only convergence assumed in the definition of a twisted intertwining operator is the convergence stated in the duality axiom which, in this case, is equivalent to the Jacobi identity \eqref{eq:twisted jacobi c5}. 

\subsection{Algebraic computations}
We now use the Jacobi identity together with purely algebraic techniques to derive the identity \eqref{eq:correlator identity c5}. We can substitute new formal variables $ \xi $ and $ y_1, \dots, y_\N $ into \eqref{eq:twisted jacobi c5} to obtain
\begin{equation}\label{eq:Jacobi with y and xi}
\begin{aligned}
& y_i^{-1} \delta \left( \frac{\xi - x_i}{y_i} \right)  Y((\xi - x_i)^{\mathcal{L}_{g_i}} \xi^{\mathcal{L}_{\g_i}} u,\xi) \mathcal{Y}(w_i, x_i) \\
& \quad - y_i^{-1} \delta \left( \frac{- x_i + \xi}{y_i} \right) \mathcal{Y}(w_i, x_i) Y((- x_i + \xi)^{\mathcal{L}_{g_i}} \xi^{\mathcal{L}_{\g_i}} u,\xi) \\
\qquad &= \xi^{-1} \delta \left( \frac{x_i + y_i}{\xi} \right) \mathcal{Y}(Y( y_i^{\mathcal{L}_{g_i}} (x_i + y_i)^{\mathcal{L}_{\g_i}} u,y_i)w_i,x_i) 
\end{aligned}
\end{equation}
for each $ i = 1, \dots, \N $.

For the $ i^{\text{th}} $ summand below, we apply \eqref{eq:Jacobi with y and xi} with $ u $ replaced with 
\begin{align*}
\prod_{j = 1}^{i-1} (- x_j + \xi)^{\mathcal{L}_{g_j}} \xi^{\mathcal{L}_{\g_j}} \cdot \prod_{j = i+1}^{\N} (\xi - x_j)^{\mathcal{L}_{g_j}} \xi^{\mathcal{L}_{\g_j}} \cdot u.
\end{align*} 
Then we obtain
\begin{equation} \label{eq:intertwining operator identity c5}
\begin{aligned}
& \sum_{i = 1}^{\N}  
	\prod_{j = 1}^{i-1}  y_j^{-1} \delta \left( \frac{- x_j + \xi}{y_j} \right) \cdot \xi^{-1} \delta \left( \frac{x_i + y_i}{\xi} \right)\cdot \prod_{j = i+1}^{\N}  y_j^{-1} \delta \left( \frac{\xi - x_j}{y_j} \right) \cdot \mathcal{Y}(w_1, x_1) \cdots \\
	&\quad \cdot \mathcal{Y}\left( Y \left(
	\prod_{j = 1}^{i-1} (- x_j + \xi)^{\mathcal{L}_{g_j}} \xi^{\mathcal{L}_{\g_j}} \cdot y_i^{\mathcal{L}_{g_i}} (x_i + y_i)^{\mathcal{L}_{\g_i}} \cdot \prod_{j = i+1}^{\N} (\xi - x_j)^{\mathcal{L}_{g_j}} \xi^{\mathcal{L}_{\g_j}} \cdot
	u,y_i \right) w_i,x_i \right) \cdots \mathcal{Y}(w_{\N}, x_{\N})  \\
& = \sum_{i = 1}^{\N} 
	\prod_{j = 1}^{i-1} y_j^{-1} \delta \left( \frac{- x_j + \xi}{y_j} \right) \cdot
	\prod_{j = i}^{\N} y_j^{-1} \delta \left( \frac{\xi - x_j}{y_j} \right) \\
&\qquad \cdot \mathcal{Y}(w_1, x_1) \cdots Y\left( \prod_{j = 1}^{i-1} (- x_j + \xi)^{\mathcal{L}_{g_j}} \xi^{\mathcal{L}_{\g_j}} \cdot \prod_{j = i}^{\N} (\xi - x_j)^{\mathcal{L}_{g_j}} \xi^{\mathcal{L}_{\g_j}} \cdot u,\xi \right) \mathcal{Y}(w_i, x_i) \cdots \mathcal{Y}(w_{\N}, x_{\N})  \\ 	
&\quad - \sum_{i = 1}^{\N} 
	\prod_{j = 1}^{i} y_j^{-1} \delta \left( \frac{- x_j + \xi}{y_j} \right) \cdot
	\prod_{j = i+1}^{\N} y_j^{-1} \delta \left( \frac{\xi - x_j}{y_j} \right) \\
& \quad\qquad  \cdot \mathcal{Y}(w_1, x_1) \cdots \mathcal{Y}(w_i, x_i) Y \left( \prod_{j = 1}^{i} (- x_j + \xi)^{\mathcal{L}_{g_j}} \xi^{\mathcal{L}_{\g_j}} \cdot \prod_{j = i+1}^{\N} (\xi - x_j)^{\mathcal{L}_{g_j}} \xi^{\mathcal{L}_{\g_j}} \cdot u,\xi \right) \cdots \mathcal{Y}(w_{\N}, x_{\N})  \\
& = \prod_{j = 1}^{\N} y_j^{-1} \delta \left( \frac{\xi - x_j}{y_j} \right) \cdot
	 Y \left( \prod_{j = 1}^{\N} (\xi - x_j)^{\mathcal{L}_{g_j}} \xi^{\mathcal{L}_{\g_j}} \cdot u,\xi \right) \mathcal{Y}(w_1,x_1) \cdots \mathcal{Y}(w_{\N},x_{\N}) \\
&\quad - \prod_{j = 1}^{\N} y_j^{-1} \delta \left( \frac{- x_j + \xi}{y_j} \right) \cdot \mathcal{Y}(w_1,x_1) \cdots \mathcal{Y}(w_{\N},x_{\N}) Y \left(  \prod_{j = 1}^{\N} (- x_j + \xi)^{\mathcal{L}_{g_j}} \xi^{\mathcal{L}_{\g_j}} \cdot u,\xi \right).
\end{aligned}
\end{equation} 

To use the property of $ C_n $-cofiniteness to derive differential equations, we want to explore the effects of changing $ u_{\alpha_\el - n,0} w_\el $ to $ w_\el $ in a chiral correlation function, for each $ \el = 1, \dots, \N +1 $ separately. We can replace $ u $ with
\begin{equation*}
\xi^{-\mathcal{L}_{\g_1} - \dots - \mathcal{L}_{\g_{\N - 1}}} u  
\end{equation*}
to make every power of $ \xi $ in \eqref{eq:intertwining operator identity c5} integral, furthermore with no logarithms of $ \xi $. We then have
\begin{equation} \label{eq:intertwining operator identity c6}
\begin{aligned}
& \sum_{i = 1}^{\N}  
	\prod_{j = 1}^{i-1}  y_j^{-1} \delta \left( \frac{- x_j + \xi}{y_j} \right) \cdot \xi^{-1} \delta \left( \frac{x_i + y_i}{\xi} \right)\cdot \prod_{j = i+1}^{\N}  y_j^{-1} \delta \left( \frac{\xi - x_j}{y_j} \right) \cdot \mathcal{Y}(w_1, x_1) \cdots \\
	&\quad \cdot \mathcal{Y}\left( Y \left(
	\prod_{j = 1}^{i-1} (- x_j + \xi)^{\mathcal{L}_{g_j}} \cdot y_i^{\mathcal{L}_{g_i}} (x_i + y_i)^{\mathcal{L}_{\g_i}} \cdot \prod_{j = i+1}^{\N} (\xi - x_j)^{\mathcal{L}_{g_j}} \cdot \xi^{-\mathcal{L}_{\g_i} + \mathcal{L}_{\g_\N}} u,y_i \right) w_i,x_i \right) \cdots \mathcal{Y}(w_{\N}, x_{\N})  \\
& = \prod_{j = 1}^{\N} y_j^{-1} \delta \left( \frac{\xi - x_j}{y_j} \right) \cdot
	 Y \left( \prod_{j = 1}^{\N} (\xi - x_j)^{\mathcal{L}_{g_j}} \cdot \xi^{\mathcal{L}_{\g_\N}}  u,\xi \right) \mathcal{Y}(w_1,x_1) \cdots \mathcal{Y}(w_{\N},x_{\N}) \\
&\quad - \prod_{j = 1}^{\N} y_j^{-1} \delta \left( \frac{- x_j + \xi}{y_j} \right) \cdot \mathcal{Y}(w_1,x_1) \cdots \mathcal{Y}(w_{\N},x_{\N}) Y \left(  \prod_{j = 1}^{\N} (- x_j + \xi)^{\mathcal{L}_{g_j}} \cdot \xi^{\mathcal{L}_{\g_\N}}  u,\xi \right).
\end{aligned}
\end{equation} 

After this change, triple $ \delta $-functions with $ \xi $ in the denominator can safely exchange $ \xi $'s for the numerator. Furthermore, to extract the most general relation for a product of twisted intertwining operators, we will apply $ \res_\xi \res_{y_1} \cdots \res_{y_\N} y_1^{-n_1} \cdots y_\N^{-n_\N} \xi^{ - n_{\N+1}} $ to both sides of \eqref{eq:intertwining operator identity c6}, for $ n_1, \dots, n_{\N+1} \in \mathbb{Z} $. 

First, the first term on the right-hand side of \eqref{eq:intertwining operator identity c6} can be computed as
\begin{align*}
&\res_\xi \res_{y_1} \cdots \res_{y_\N}y_1^{-n_1} \cdots y_\N^{-n_\N} \xi^{-n_{\N+1}}	\\
& \qquad \cdot \prod_{j = 1}^{\N} y_j^{-1} \delta \left( \frac{\xi - x_j}{y_j} \right) \cdot
	 Y \left( \prod_{j = 1}^{\N} (\xi - x_j)^{\mathcal{L}_{g_j}} \cdot \xi^{\mathcal{L}_{\g_\N}}  u,\xi \right) \mathcal{Y}(w_1,x_1) \cdots \mathcal{Y}(w_{\N},x_{\N}) \\
&= \res_\xi \xi^{- n_{\N+1}}  Y \left( \prod_{j = 1}^{\N} (\xi - x_j)^{\mathcal{L}_{g_j} - n_j} \cdot \xi^{\mathcal{L}_{\g_\N}}  u,\xi \right) \mathcal{Y}(w_1,x_1) \cdots \mathcal{Y}(w_{\N},x_{\N})   \\
&= \res_\xi \xi^{- n_{\N+1}}  Y \left( \sum_{k_1, \dots, k_\N = 0}^\infty \prod_{j = 1}^{\N}  \binom{\mathcal{L}_{g_j} - n_j}{k_j} \xi^{\mathcal{L}_{g_j} - n_j - k_j} (-x_j)^{k_j} \cdot \xi^{\mathcal{L}_{\g_\N}} u,\xi \right) \mathcal{Y}(w_1,x_1) \cdots \mathcal{Y}(w_{\N},x_{\N})   \\
&= \res_\xi \sum_{k_1, \dots, k_\N = 0}^\infty \xi^{\alpha_1 + \cdots + \alpha_{\N+1} - n_1 - \cdots - n_{\N+1} - k_1 - \cdots - k_\N} Y_0 \left(  \prod_{j = 1}^{\N}  \binom{\mathcal{L}_{g_j} - n_j}{k_j} (-x_j)^{k_j} \cdot u,\xi \right) \mathcal{Y}(w_1,x_1) \cdots \mathcal{Y}(w_{\N},x_{\N})   \\
&= \sum_{k_1, \dots, k_\N = 0}^\infty \left(  \prod_{j = 1}^{\N} \binom{\mathcal{L}_{g_j} - n_j}{k_j} (-x_j)^{k_j} u \right)_{\alpha_1 + \cdots + \alpha_{\N+1} - n_1 - \cdots - n_{\N+1} - k_1 - \cdots - k_\N,0} \mathcal{Y}(w_1,x_1) \cdots \mathcal{Y}(w_{\N},x_{\N}).
\end{align*}
The second term on the right-hand side can be computed with
\begin{align*}
&\res_\xi \res_{y_1} \cdots \res_{y_\N}y_1^{-n_1} \cdots y_\N^{-n_\N} \xi^{-n_{\N+1}}  \\
& \qquad \cdot \prod_{j = 1}^{\N} y_j^{-1} \delta \left( \frac{- x_j + \xi}{y_j} \right) \cdot \mathcal{Y}(w_1,x_1) \cdots \mathcal{Y}(w_{\N},x_{\N}) Y \left(  \prod_{j = 1}^{\N} (- x_j + \xi)^{\mathcal{L}_{g_j}} \cdot \xi^{\mathcal{L}_{\g_\N}}  u,\xi \right) \\
&= \res_\xi \xi^{-n_{\N+1}} \mathcal{Y}(w_1,x_1) \cdots \mathcal{Y}(w_{\N},x_{\N}) Y \left(  \prod_{j = 1}^{\N} (- x_j + \xi)^{\mathcal{L}_{g_j} - n_j} \cdot \xi^{\mathcal{L}_{\g_\N}}  u,\xi \right) \\
&= \res_\xi \xi^{-n_{\N+1}} \mathcal{Y}(w_1,x_1) \cdots \mathcal{Y}(w_{\N},x_{\N}) Y \left( \sum_{k_1, \dots, k_\N = 0}^\infty \prod_{j = 1}^{\N} \binom{\mathcal{L}_{g_j} -n_j}{k_j}(- x_j)^{\mathcal{L}_{g_j}-n_j-k_j} \xi^{k_j}  \cdot \xi^{\mathcal{L}_{\g_\N}}  u,\xi \right) \\
&= \res_\xi \sum_{k_1, \dots, k_\N = 0}^\infty \xi^{\alpha_{\N+1} - n_{\N+1} + k_1 + \cdots + k_\N}  \mathcal{Y}(w_1,x_1) \cdots \mathcal{Y}(w_{\N},x_{\N}) Y_0 \left(  \prod_{j = 1}^{\N} \binom{\mathcal{L}_{g_j} -n_j}{k_j}(- x_j)^{\mathcal{L}_{g_j}-n_j-k_j} \cdot u,\xi \right) \\
&= \sum_{k_1, \dots, k_\N = 0}^\infty  \mathcal{Y}(w_1,x_1) \cdots \mathcal{Y}(w_{\N},x_{\N}) \left(  \prod_{j = 1}^{\N} \binom{\mathcal{L}_{g_j} -n_j}{k_j}(- x_j)^{\mathcal{L}_{g_j}-n_j-k_j} \cdot  u \right)_{\alpha_{\N+1} - n_{\N+1} + k_1 + \cdots + k_\N,0} .
\end{align*}

To compute the left-hand side, we compute the following, for each $ i = 1, \dots, \N $.

\begingroup
\allowdisplaybreaks
\begin{align*}
&\res_\xi \res_{y_1} \cdots \res_{y_\N}y_1^{-n_1} \cdots y_\N^{-n_\N} \xi^{-n_{\N+1}}  \\
 &\quad \cdot \prod_{j = 1}^{i-1}  y_j^{-1} \delta \left( \frac{- x_j + \xi}{y_j} \right) \cdot \xi^{-1} \delta \left( \frac{x_i + y_i}{\xi} \right)\cdot \prod_{j = i+1}^{\N}  y_j^{-1} \delta \left( \frac{\xi - x_j}{y_j} \right) \cdot \mathcal{Y}(w_1, x_1) \cdots \\
	&\quad \cdot \mathcal{Y}\left( Y \left(
	\prod_{j = 1}^{i-1} (- x_j + \xi)^{\mathcal{L}_{g_j}} \cdot y_i^{\mathcal{L}_{g_i}} (x_i + y_i)^{\mathcal{L}_{\g_i}} \cdot \prod_{j = i+1}^{\N} (\xi - x_j)^{\mathcal{L}_{g_j}} \cdot \xi^{-\mathcal{L}_{\g_i} + \mathcal{L}_{\g_\N}} u,y_i \right) w_i,x_i \right) \cdots \mathcal{Y}(w_{\N}, x_{\N})  \\
&=\res_\xi \res_{y_i} y_i^{-n_i} \xi^{-n_{\N+1}} \xi^{-1} \delta \left( \frac{x_i + y_i}{\xi} \right) \mathcal{Y}(w_1, x_1) \cdots \\
	&\quad \cdot \mathcal{Y}\left( Y \left(
	\prod_{j = 1}^{i-1} (- x_j + \xi)^{\mathcal{L}_{g_j}-n_j} \cdot y_i^{\mathcal{L}_{g_i}} (x_i + y_i)^{\mathcal{L}_{\g_i}} \cdot \prod_{j = i+1}^{\N} (\xi - x_j)^{\mathcal{L}_{g_j}-n_j} \cdot \xi^{-\mathcal{L}_{\g_i} + \mathcal{L}_{\g_\N}} u,y_i \right) w_i,x_i \right) \cdots \mathcal{Y}(w_{\N}, x_{\N})  \\
&=\res_\xi \res_{y_i} y_i^{-n_i} \mathcal{Y}(w_1, x_1) \cdots  \mathcal{Y}\left( \left\{ \xi^{-1} \delta \left( \frac{x_i + y_i}{\xi} \right) Y \left(
	\prod_{j = 1}^{i-1} (- x_j + \xi)^{\mathcal{L}_{g_j}-n_j} \cdot y_i^{\mathcal{L}_{g_i}} (x_i + y_i)^{\mathcal{L}_{\g_i}} \right. \right. \right. \\
&\qquad \left. \left. \left. \cdot \prod_{j = i+1}^{\N} (\xi - x_j)^{\mathcal{L}_{g_j}-n_j} \cdot \xi^{-\mathcal{L}_{\g_i} + \mathcal{L}_{\g_\N} - n_{\N+1}} u,y_i \right) w_i \right\} ,x_i \right) \cdots \mathcal{Y}(w_{\N}, x_{\N})  \\
&=\res_\xi \res_{y_i} y_i^{-n_i} \mathcal{Y}(w_1, x_1) \cdots  \mathcal{Y}\left( e^{y_i \frac{\partial}{\partial x_i}} \left\{ \xi^{-1} \delta \left( \frac{x_i}{\xi} \right) Y \left(
	\prod_{j = 1}^{i-1} (- x_j + \xi)^{\mathcal{L}_{g_j}-n_j} \cdot y_i^{\mathcal{L}_{g_i}} x_i^{\mathcal{L}_{\g_i}} \right. \right. \right. \\
&\qquad \left. \left. \left. \cdot \prod_{j = i+1}^{\N} (\xi - x_j)^{\mathcal{L}_{g_j}-n_j} \cdot \xi^{-\mathcal{L}_{\g_i} + \mathcal{L}_{\g_\N} - n_{\N+1}} u,y_i \right) w_i \right\} ,x_i \right) \cdots \mathcal{Y}(w_{\N}, x_{\N})  \\
&=\res_\xi \res_{y_i} y_i^{-n_i} \mathcal{Y}(w_1, x_1) \cdots  \mathcal{Y}\left( e^{y_i \frac{\partial}{\partial x_i}} \left\{ \xi^{-1} \delta \left( \frac{x_i}{\xi} \right) Y \left(
	\prod_{j = 1}^{i-1} (- x_j + x_i)^{\mathcal{L}_{g_j}-n_j} \cdot y_i^{\mathcal{L}_{g_i}} x_i^{\mathcal{L}_{\g_i}} \right. \right. \right. \\
&\qquad \left. \left. \left. \cdot \prod_{j = i+1}^{\N} (x_i - x_j)^{\mathcal{L}_{g_j}-n_j} \cdot x_i^{-\mathcal{L}_{\g_i} + \mathcal{L}_{\g_\N} - n_{\N+1}} u,y_i \right) w_i \right\} ,x_i \right) \cdots \mathcal{Y}(w_{\N}, x_{\N})  \\
&= \res_{y_i} y_i^{-n_i} \mathcal{Y}(w_1, x_1) \cdots  \mathcal{Y}\left( e^{y_i \frac{\partial}{\partial x_i}} \left\{  Y \left(
	\prod_{j = 1}^{i-1} (- x_j + x_i)^{\mathcal{L}_{g_j}-n_j} \cdot y_i^{\mathcal{L}_{g_i}} x_i^{\mathcal{L}_{\g_i}} \right. \right. \right. \\
&\qquad \left. \left. \left. \cdot \prod_{j = i+1}^{\N} (x_i - x_j)^{\mathcal{L}_{g_j}-n_j} \cdot x_i^{-\mathcal{L}_{\g_i} + \mathcal{L}_{\g_\N} - n_{\N+1}} u,y_i \right) w_i \right\} ,x_i \right) \cdots \mathcal{Y}(w_{\N}, x_{\N})  \\
&= \res_{y_i} y_i^{-n_i} \mathcal{Y}(w_1, x_1) \cdots  \mathcal{Y}\left(  Y \left(
	\prod_{\substack{j = 1 \\ j \neq i} }^{\N + 1}  \sum_{k_j = 0}^\infty y_i^{k_j} \binom{\mathcal{L}_{g_j}-n_j}{k_j} \cln(x_i - x_j)^{\mathcal{L}_{g_j}-n_j-k_j} \cln \cdot y_i^{\mathcal{L}_{g_i}}  u,y_i \right) w_i ,x_i \right) \\
	&\qquad \cdots \mathcal{Y}(w_{\N}, x_{\N})  \\
&= \res_{y_i} y_i^{-n_i} \sum_{k_1, \dots, \widehat{k_i}, \dots, k_{\N+1} = 0}^\infty  \mathcal{Y}(w_1, x_1) \cdots \\
	&\qquad  \cdot \mathcal{Y}\left(  Y \left(
	\prod_{\substack{j = 1 \\ j \neq i} }^{\N + 1} y_i^{k_j} \binom{\mathcal{L}_{g_j}-n_j}{k_j} \cln(x_i - x_j)^{\mathcal{L}_{g_j}-n_j-k_j} \cln \cdot y_i^{\mathcal{L}_{g_i}}  u,y_i \right) w_i ,x_i \right) \cdots \mathcal{Y}(w_{\N}, x_{\N}) \\
&= \res_{y_i} y_i^{-n_i} \sum_{k_1, \dots, \widehat{k_i}, \dots, k_{\N+1} = 0}^\infty  \mathcal{Y}(w_1, x_1) \cdots \\
	&\qquad  \cdot \mathcal{Y}\left(  Y_0 \left(
	\prod_{\substack{j = 1 \\ j \neq i} }^{\N + 1} y_i^{k_j} \binom{\mathcal{L}_{g_j}-n_j}{k_j} \cln(x_i - x_j)^{\mathcal{L}_{g_j}-n_j-k_j} \cln \cdot y_i^{\alpha_i}  u,y_i \right) w_i ,x_i \right) \cdots \mathcal{Y}(w_{\N}, x_{\N}) \\
&= \sum_{k_1, \dots, \widehat{k_i}, \dots, k_{\N+1} = 0}^\infty  \mathcal{Y}(w_1, x_1) \cdots \\
	&\qquad  \cdot \mathcal{Y}\left(  \left(
	\prod_{\substack{j = 1 \\ j \neq i} }^{\N + 1} \binom{\mathcal{L}_{g_j}-n_j}{k_j} \cln(x_i - x_j)^{\mathcal{L}_{g_j}-n_j-k_j} \cln \cdot u \right)_{\alpha_i - n_i + k_1 + \cdots + \widehat{k_i} + \cdots + k_{\N+1} ,0} w_i ,x_i \right) \cdots \mathcal{Y}(w_{\N}, x_{\N}).
\end{align*}
\endgroup

We have continued to use the notation $ x_{\N+1} = 0 $ and
\begin{align*}
\cln(x_i-x_j)^{A}\cln &= \begin{cases}
(x_i-x_j)^{A} = \sum_{k = 0}^\infty \binom{A}{k} x_i^{A-k} (-x_j)^k  &\text{if } i < j, \\
(-x_j+x_i)^{A} = \sum_{k = 0}^\infty \binom{A}{k} (-x_j)^{A-k} x_i^k  &\text{if } j < i,
\end{cases} 
\end{align*}
for operators $ A $ where the right-hand side is defined.

Thus, we have
\begin{equation}\label{eq:correlator identity 2 c5}
\begin{aligned}
 & \sum_{ i = 1}^\N \sum_{k_1, \dots, \widehat{k_i}, \dots, k_{\N+1} = 0}^\infty  \mathcal{Y}(w_1, x_1) \cdots \\
	&\qquad  \cdot \mathcal{Y}\left(  \left(
	\prod_{\substack{j = 1 \\ j \neq i} }^{\N + 1} \binom{\mathcal{L}_{g_j}-n_j}{k_j} \cln(x_i - x_j)^{\mathcal{L}_{g_j}-n_j-k_j} \cln \cdot u \right)_{\alpha_i - n_i + k_1 + \cdots + \widehat{k_i} + \cdots + k_{\N+1} ,0} w_i ,x_i \right) \cdots \mathcal{Y}(w_{\N}, x_{\N}) \\
	&= \sum_{k_1, \dots, k_\N = 0}^\infty \left(  \prod_{j = 1}^{\N} \binom{\mathcal{L}_{g_j} - n_j}{k_j} (-x_j)^{k_j} u \right)_{\alpha_1 + \cdots + \alpha_{\N+1} - n_1 - \cdots - n_{\N+1} - k_1 - \cdots - k_\N,0} \mathcal{Y}(w_1,x_1) \cdots \mathcal{Y}(w_{\N},x_{\N}) \\
&  \quad  - \sum_{k_1, \dots, k_\N = 0}^\infty  \mathcal{Y}(w_1,x_1) \cdots \mathcal{Y}(w_{\N},x_{\N}) \left(  \prod_{j = 1}^{\N} \binom{\mathcal{L}_{g_j} -n_j}{k_j}(- x_j)^{\mathcal{L}_{g_j}-n_j-k_j} \cdot  u \right)_{\alpha_{\N+1} - n_{\N+1} + k_1 + \cdots + k_\N,0} ,
\end{aligned}
\end{equation}
which shows that \eqref{eq:correlator identity c5} is indeed in the kernel of $ \phi $.

The identity \eqref{eq:correlator identity 2 c5} tells us how to define the generators as elements of $ T $. For each $ \el = 1, \dots, \N+1 $,
we replace $ u $ in \eqref{eq:correlator identity 2 c5} with
\begin{equation*}
u^{(\el)} = \left( \prod_{\substack{j = 1 \\ j \neq \el} }^{\N + 1} \cln(x_\el-x_j)^{-\mathcal{L}_{g_j} + n_j}\cln u \right).
\end{equation*}
We then define
\begin{equation} \label{eq:A_k defn c5}
\begin{aligned}
& \mathcal{A}_\el(u, w_0, \dots, w_{\N+1}; n_1, \dots, n_{\N+1}) \\
&= -\sum_{k_1, \dots, k_\N = 0}^\infty (-x_1)^{k_1} \cdots (-x_\N)^{k_\N}  \left(  \prod_{j=1}^{\N} \binom{\mathcal{L}_{g_j} - n_j}{k_j} u^{(\el)} \right)_{\alpha_1 + \cdots + \alpha_{\N+1} - n_1 - \cdots - n_{\N+1} - k_1 - \cdots - k_\N, 0}^* w_0 \\
&\qquad \otimes w_1 \otimes \cdots \otimes w_{\N+1} \\
& + \sum_{i = 1}^{\N+1} \sum_{k_1, \dots, \widehat{k_i}, \dots, k_{\N+1} = 0}^\infty  w_0 \otimes w_1 \otimes \cdots \otimes  w_{i-1} \\
& \quad \otimes \left( \prod_{\substack{j = 1 \\ j \neq i} }^{\N + 1} \binom{\mathcal{L}_{g_j} - n_j}{k_j} \cln(x_i-x_j)^{\mathcal{L}_{g_j} - n_j-k_j}\cln u^{(\el)} \right)_{\alpha_i - n_i + k_1 + \cdots + \widehat{k_i} + \cdots + k_{\N+1},0} w_i \otimes w_{i+1} \otimes \cdots \otimes w_{\N+1},
\end{aligned}
\end{equation}
 For all $ \el = 1, \dots, \N + 1$, for all simultaneous generalized eigenvectors $ u \in V $ for $ G $, for all $ w_i \in W_i $, and for all $ n_i \in \mathbb{Z} $, equation \eqref{eq:correlator identity 2 c5} shows that the element $ \mathcal{A}_\el(u, w_0, \dots, w_{\N+1}; n_1, \dots, n_{\N+1}) $ is in the kernel of $ \phi $. Furthermore, we have the following lemma.

\begin{lemma}\label{lem:A_l in T_K_0 c5}
Let $ K_0 \in \mathbb{R} $. Let $ u \in V_+ $ be a simultaneous generalized eigenvector for $ G $. Let $ \alpha_i \in \mathbb{C} $ with $ \Re(\alpha_i) \in [0,1) $ be such that $ \mathcal{S}_{g_i} u = \alpha_i u $. Let $ w_0 \in \coprod_{\substack{ m \in \mathbb{C} \\ \Re(m) \leq K_0}} (W_0)_{[m]} $, $ w_1 \in W_1 , \dots, w_{\N+1} \in W_{\N+1} $. Then 
\begin{align*}
\mathcal{A}_\el(u, w_0, \dots, w_{\N+1}; n_1, \dots, n_{\N+1}) \in T_{K_0}
\end{align*}
whenever $ n_1 + \cdots + n_{\N+1} \geq \Re(\alpha_1 + \cdots + \alpha_{\N+1}) $.
\end{lemma}

\begin{proof}
Assume that $ u, w_0, \dots, w_{\N+1} $ are all of homogeneous weight. By explicitly looking at \eqref{eq:A_k defn c5}, we see that all of the terms in $ \mathcal{A}_\el(u, w_0, \dots, w_{\N+1}; n_1, \dots, n_{\N+1}) $ are clearly in $ T_{K_0} $, except for possibly the terms
\begin{align*}
&-(-x_1)^{k_1} \cdots (-x_\N)^{k_\N}  \left(  \prod_{j=1}^{\N} \binom{\mathcal{L}_{g_j} - n_j}{k_j} u^{(\el)} \right)_{\alpha_1 + \cdots + \alpha_{\N+1} - n_1 - \cdots - n_{\N+1} - k_1 - \cdots - k_\N, 0}^* w_0 \\
&\qquad \otimes w_1 \otimes \cdots \otimes w_{\N+1} .
\end{align*}
Since the operators $ \mathcal{L}_{g_j} $ do not shift the weight of $ u $, each of the terms above has the same weight as 
\begin{align*}
 u_{\alpha_1 + \cdots + \alpha_{\N+1} - n_1 - \cdots - n_{\N+1} - k_1 - \cdots - k_\N, 0}^* w_0 \otimes w_1 \otimes \cdots \otimes w_{\N+1} .
\end{align*}
So we need to check that 
\begin{align*}
\Re( \wt  u_{\alpha_1 + \cdots + \alpha_{\N+1} - n_1 - \cdots - n_{\N+1} - k_1 - \cdots - k_\N,0}^* w_0  ) \leq K_0 ,
\end{align*}
for some fixed $ n_1, \dots,n_{\N+1} $ and for all $k_1, \dots, k_\N \geq 0 $. Indeed, for all $ k_1, \dots, k_\N \geq 0 $, we have
\begin{align*}
& \Re( \wt  u_{\alpha_1 + \cdots + \alpha_{\N+1} - n_1 - \cdots - n_{\N+1} - k_1 - \cdots - k_\N,0}^* w_0  ) \\
&\quad = - \wt u + \Re(\alpha_1 + \cdots + \alpha_{\N+1}) - n_1 - \cdots - n_{\N+1} - k_1 - \cdots - k_\N + 1 + \Re(\wt w_0) \\
&\quad \leq  \Re(\alpha_1 + \cdots + \alpha_{\N+1}) - n_1 - \cdots - n_{\N+1} + K_0 \\
&\quad \leq  K_0
\end{align*} 
whenever $ \Re(\alpha_1 + \cdots + \alpha_{\N+1}) - n_1 - \cdots - n_{\N+1} \leq 0 $.
\end{proof}

\begin{remark}
The method used here to derive the generators $ \mathcal{A}_\el $ is different from that used in \cite{TanMostlyUntwisted}, which uses a commutator formula for $ [Y(u,x_1), \mathcal{Y}(w,x_2)] $. This commutator was also used to derive the twisted KZ equations. Such a commutator does not exist in the current setting since the closest identity available from \eqref{eq:twisted jacobi c5} would be for $ Y \left( (x_1 - x_2)^{\mathcal{L}_{g_1}} u,x_1 \right) \mathcal{Y}(w,x_2) -  \mathcal{Y}(w,x_2) Y \left( (-x_2 + x_1)^{\mathcal{L}_{g_1}} u,x_1 \right) $, which does not provide a clean commutator formula when $ g_1 \neq \id_V $. This is why the method used here appears more complicated than that in \cite{TanMostlyUntwisted}. One should observe that the Jacobi identity always encodes some analogue of a commutation formula, so the two methods are in fact still closely related.
\end{remark}

\begin{remark}
Note that in papers such as \cite{HuangDEs} and \cite{TanMostlyUntwisted} that use Huang's method of deriving differential equations for products of intertwining operators, the submodule $ J $ had additional generators $ \mathcal{A}_0 $. The method used to derive the differential equations introduced in \cite{HuangC1} does not need the $ \mathcal{A}_0 $ elements as generators. This is an improvement to the method used in \cite{HuangDEs} and \cite{TanMostlyUntwisted}. We will use this new method in this paper, which is why we have only computed $ \mathcal{A}_\el $ for $ \el = 1, \dots, \N +1 $, corresponding to the in-state modules. As a consequence, we do not need to assume that $ W_0 $ is $ C_n $-cofinite.
\end{remark}

\begin{remark}
Formula \eqref{eq:A_k defn c5} looks quite complicated compared to the analogous generators $ \mathcal{A} $, $ \mathcal{B} $ and $ \mathcal{C} $ in the original untwisted case in \cite{HuangDEs}. Our generators have ``infinitely" more additional terms in the sums over $ k_j $ that do not appear in the untwisted case. Of course, by the lower truncation property or the fact that discretely graded modules are lower bounded, all sums are actually finite.

By setting $ \mathcal{L}_{g_i} = 0 $ and $ n_j = 0 $ for all $  i,j = 1, \dots, \N+1 $ with $ j \neq \el $, and setting $ n_\el = 1 $, we see that most of the binomial coefficients in  \eqref{eq:A_k defn c5} vanish. We then recover the operators  $ \mathcal{A} $, $ \mathcal{B} $ and $ \mathcal{C} $ (up to an overall negative sign, which is just a definition convention). 

A closer inspection of our derivation reveals that the extra terms come from the factors $ (x_1 - x_2)^{\alpha_1} $, $ (- x_2 + x_1)^{\alpha_1} $ and $ (x_2 + x_0)^{\alpha_2} $ in the Jacobi identity \eqref{eq:twisted jacobi c5}. Alternatively, in the heuristic derivation, we can see that the analogous factors are introduced to obtain a single-valued correlation function. After single-valuedness is achieved, we must expand these factors during the contour integration, producing the extra terms. 

It is interesting to observe that these additional terms do not make that much of a difference beyond increasing the complexity of the calculations. In the following section, we will see that the additional terms all have a sufficiently low weight to not complicate the derivation of differential equations with regular singularities.

The additional factors $ (x_1 - x_2)^{\alpha_1} $, $ (- x_2 + x_1)^{\alpha_1} $ and $ (x_2 + x_0)^{\alpha_2} $ do, however, shift the weights in a manner that motivates the use of $ C_2 $-cofinite modules as an assumption.
\end{remark}

\section{Differential equations}\label{sec:differential equations}

In this section, we follow a method similar to those used in \cite{HuangDEs}, \cite{HuangC1} and \cite{TanMostlyUntwisted} to derive differential equations with certain coefficients satisfied by products of $ \N $ twisted intertwining operators among discretely graded $ C_n $-cofinite twisted modules twisted by elements in $ G $. For modules twisted by elements with only positive real eigenvalues, we can relax their assumptions to discretely grading and $ C_1 $-cofiniteness. In Section~\ref{sec:convergence}, when we further assume that $ G $ is finite, we will use these differential equations to prove the convergence (and extension) property for products of said twisted intertwining operators.

For each $ K_0 \in \mathbb{R} $, define $ J_{K_0} $ to be the $ R $-submodule of $ T_{K_0} $ generated by the elements $ \mathcal{A}_{\el}(u,w_0, \dots, w_{\N+1}; n_1, \dots, n_{\N+1}) $, for all $ \el = 1, \dots, \N + 1 $, for all simultaneous generalized eigenvectors $ u \in V_+ $ of $ G $, for all $ w_0 \in \coprod_{\substack{ m \in \mathbb{C} \\ \Re(m) \leq K_0}} (W_0)_{[m]} $ and $ w_1 \in W_1 , \dots, w_{\N+1} \in W_{\N+1} $, and for all $ (n_1, \dots, n_{\N+1}) \in \mathbb{Z}^{\N+1} $ such that $ n_1 + \cdots + n_{\N+1} \geq \Re(\alpha_1 + \cdots + \alpha_{\N+1} )$. By Lemma~\ref{lem:A_l in T_K_0 c5}, we see that these generators are indeed elements of $ T_{K_0} $.

Define $ F_r(J_{K_0}) = F_r(T_{K_0}) \cap J_{K_0} $. We then have the filtration $ F_r(J_{K_0}) \subseteq F_s(J_{K_0}) $, $ r \leq s \in \mathbb{R}$. 

Recall that in Section~\ref{sec:setting c5}, we introduced a parameter $ n $ to control the $ C_n $-cofiniteness conditions for certain twisted modules. We claim that $ C_2 $-cofiniteness of these modules, i.e. $ n = 2 $, is a strong enough condition to derive differential equations when $ R $ is assumed to be Noetherian. Furthermore, these differential equations have coefficients defined on the domain $ M^\N = \{ (z_1, \dots, z_\N) \in \mathbb{C}^\N : z_i \neq 0,\ z_j \neq z_k \text{ when } j \neq k \} $. In Remark~\ref{rem:C_1-cofinite general c5}, we will explain why $ C_1 $-cofiniteness is not enough in general, but is enough in special cases. This suggests that $ C_2 $-cofiniteness is really needed for our method (however, we do not have a counterexample of a product of twisted intertwining operators among $ C_1 $-cofinite twisted modules that does not converge). 

\begin{lemma} \label{lem: generator filter c5}
Let $ K_0 \in \mathbb{R} $. Let $ u \in V_+ $ be a weight-homogeneous simultaneous generalized eigenvector for $ G $ and let $ w_i \in W_i $, $ i = 0, \dots, \N+1 $, be weight-homogeneous with $ \Re(\wt w_0) \leq K_0 $. Let 
\begin{align*}
w^{(\el)} &=  w_0 \otimes w_1 \otimes \cdots \otimes u_{\alpha_\el - n_\el,0}  w_\el \otimes \cdots \otimes w_{\N+1}, \qquad \text{ for } \el = 1, \dots, \N + 1,
\end{align*} 
where $ n_\el = n $ if $ g_\el $ has an eigenvalue not in $ \mathbb{R}_{>0} $, and $ n_\el = 1 $ if $ g_\el $ only has  eigenvalues in $ \mathbb{R}_{>0} $. Let $ \Re (\wt w^{(\el)}) = s $. If $ n \geq  2 $, then $ w^{(\el)} \in F_s(J_{ K_0}) + F_r(T_{ K_0}) $ for some $ r < s $.
\end{lemma}

\begin{proof}
Let $ \el = 1, \dots, \N + 1 $. We show that there exists $ (n_1, \dots, n_{\N+1}) \in \mathbb{Z}^{\N+1} $ with $ n_1 + \cdots + n_{\N+1} \geq \Re( \alpha_1 + \cdots + \alpha_{\N+1} ) $ such that $ \mathcal{A}_\el(u, w_0, \dots, w_{\N+1}; n_1, \dots, n_{\N+1}) - w^{(\el)} \in F_r(T_{K_0}) $ for some $ r < s $. Choose $ n_\el = n $ or $ n_\el = 1 $ depending on $ g_\el $, and leave the remaining $ n_i \in \mathbb{Z} $, $ i = 1, \dots, \widehat{\el}, \dots, \N+1 $, unchosen for now. We check the requirements for the real part of the weight of each term in $ \mathcal{A}_{\el}(u,w_0, \dots, w_{\N+1}; n_1, \dots, n_{\N+1}) - w^{(\el)} $ to be less than $ s $. Since the action of $ \mathcal{L}_{g_i} $ does not shift weights, it suffices to show that the real part of the weight of
\begin{align*}
w_0 \otimes w_1 \otimes \cdots \otimes u_{\alpha_\el - n_\el,0}  w_\el \otimes \cdots \otimes w_{\N+1}
\end{align*}
is greater than the real parts of the weights of
\begin{gather*}
u_{\alpha_1 + \cdots + \alpha_{\N+1} - n_1 - \cdots - n_{\N+1} - k_1 - \cdots - k_\N,0}^* w_0 \otimes w_1 \otimes \cdots  \otimes w_{\N+1} , \text{ and }  \\
w_0 \otimes w_1 \otimes \cdots \otimes u_{\alpha_i - n_i + k_1 + \cdots + \widehat{k_i} + \cdots + k_{\N+1},0}  w_i \otimes \cdots  \otimes w_{\N+1} , \ \  \text{for } i = 1, \dots, \N+1, 
\end{gather*}
for all $ (k_1, \dots, k_{\N+1} ) \in \mathbb{Z}_{\geq 0}^{\N+1} $ with some $ k_j > 0 $, $ j \neq \el $, when $ i = \el $, or for all $ (k_1, \dots, k_{\N+1} ) \in \mathbb{Z}_{\geq 0}^{\N+1} $, otherwise. Using $ \sigma $ to denote $ \Re( \wt w_0 + \cdots + \wt w_{\N+1}) $, we see that the condition
\begin{equation}
\begin{aligned}
&\sigma - \wt u + \Re(\alpha_1 + \cdots + \alpha_{\N+1}) - n_1 - \cdots - n_{\N+1} - k_1 - \cdots - k_\N + 1 \\
&\quad < \sigma + \wt u - \Re(\alpha_\el) + n_\el - 1,  \\
&\sigma + \wt u - \Re(\alpha_i) + n_i - k_1 - \cdots - \widehat{k_i} - \cdots - k_{\N+1} - 1 \\
&\quad < \sigma + \wt u - \Re(\alpha_\el) + n_\el - 1 , \quad \text{ for } i = 1, \dots, \N+1,
\end{aligned}
\end{equation}
for all $ (k_1, \dots, k_{\N+1} ) \in \mathbb{Z}_{\geq 0}^{\N+1} $ with some $ k_j > 0 $, $ j \neq \el $, when $ i = \el $, or for all $ (k_1, \dots, k_{\N+1} ) \in \mathbb{Z}_{\geq 0}^{\N+1} $, otherwise, is equivalent to
\begin{equation}\label{eq:inequality condition c5}
\begin{gathered}
\Re(\alpha_1 + \cdots + \alpha_{\N+1} +  \alpha_\el)  < 2 \wt u + n_1 + \cdots + n_{\N+1} + n_\el - 2,  \\
\Re(\alpha_\el - \alpha_i) < n_\el - n_i , \quad \text{ for } i \neq \el.
\end{gathered}
\end{equation}
First, consider the case where $ g_\el $ has an eigenvalue not in $ \mathbb{R}_{>0} $. If we let $ \Re(\alpha_j) \in [0,1) $ and $ \wt u \geq 1 $ be arbitrary, we see that a necessary and sufficient condition is
\begin{equation}\label{eq:inequality arbitray eigenvalues c5}
\begin{gathered}
N  + 2 \leq n_1 + \cdots + n_{\N+1} + n_\el, \text{ and } \\
n_i \leq n_\el - 1 , \quad \text{ for } i \neq \el.
\end{gathered}
\end{equation}
This is easily satisfied if $ n_\el = n \geq 2 $. For example, we could pick $ n_i = 1 $ for all $ i \neq \el $. Furthermore, $ n_1 + \cdots + n_{\N+1} = N + n > \Re(\alpha_1 + \cdots + \alpha_{\N+1}) $, so that $ \mathcal{A}_\el(u, w_0, \dots, w_{\N+1}; n_1, \dots, n_{\N+1}) \in T_{K_0} $.

Second, consider the case where $ g_\el $ only has eigenvalues in $ \mathbb{R}_{>0} $, so we choose $ n_\el = 1 $. Then $ \Re(\alpha_\el) = 0 $. There are now two cases when choosing $ n_i $, $ i \neq \el $. If $ \Re(\alpha_\el - \alpha_i) < 0 $ (which implies $ 0 < \Re(\alpha_i) < 1 $), we choose $ n_i = 1 $. If $ \Re(\alpha_\el - \alpha_i) = 0 $ (which implies $ \Re(\alpha_i) = 0 $), we choose $ n_i = 0 $. In both cases, we have
$ \Re(\alpha_\el - \alpha_i) < n_\el - n_i $ and $ \Re(\alpha_i) \leq n_i $ when $ i \neq \el $. So, 
\begin{align*}
\Re(\alpha_1 + \cdots + \widehat{\alpha_\el} + \cdots + \alpha_{\N+1}) \leq n_1 + \cdots + \widehat{n_\el} + \cdots + n_{\N+1}.
\end{align*}
Since $ \Re(\alpha_\el) = 0 < 1 = n_\el $ and $ 0 \leq 2 \wt u - 2 $, this implies
\begin{align*}
\Re(\alpha_1 + \cdots + \alpha_{\N+1})  &<  n_1 + \cdots + n_{\N+1} , \\
\Re(\alpha_1 + \cdots + \alpha_{\N+1} + \alpha_\el ) &< 2 \wt u + n_1 + \cdots + n_{\N+1} + n_\el - 2. \qedhere
\end{align*}
\end{proof}

\begin{remark}\label{rem:choices c5}
When $ g_\el $ has an eigenvalue not in $ \mathbb{R}_{>0} $, there is some flexibility in how we pick $ (n_1, \dots, n_{\N+1} ) $. For example, if $ n_\el = n = 2 $, we could also pick $ n_j = 0 $ for one $ j \neq \el $ and the rest to be $ n_i = 1 $. With this choice, we still have
\begin{align*} 
\Re(\alpha_1 + \cdots + \alpha_{\N + 1} ) < \N + 1 &= n_1 + \cdots + n_{\N+1}, \\
N  + 2 \leq N + 3 &= n_1 + \cdots + n_{\N+1} + n_\el. 
\end{align*}
There is also some flexibility in the choice of $ (n_1, \dots, n_{\N+1} ) $ when $ g_\el $ only has eigenvalues in $ \mathbb{R}_{>0} $. We can pick $ n_j = 0 $ for one $ j \neq \el $ with $ \Re(\alpha_j) > 0 $. Observe that $ \Re(\alpha_\el - \alpha_j) < 0 < 1 = n_\el - n_j $ is still satisfied. We only have $ \Re(\alpha_i) \leq n_i $ for $ i \neq \el,j $, giving 
\begin{align*}
\Re(\alpha_1 + \cdots + \widehat{\alpha_\el} + \cdots + \widehat{\alpha_j} + \cdots + \alpha_{\N+1}) \leq n_1 + \cdots + \widehat{n_\el} + \cdots + \widehat{n_j} + \cdots + n_{\N+1}.
\end{align*}
But $ \Re(\alpha_\el + \alpha_j) = \Re(\alpha_j) < 1 = n_\el + n_j $, so we still obtain
\begin{align*}
\Re(\alpha_1 + \cdots + \alpha_{\N+1})  &<  n_1 + \cdots + n_{\N+1} , \\
\Re(\alpha_1 + \cdots + \alpha_{\N+1} + \alpha_\el ) &< 2 \wt u + n_1 + \cdots + n_{\N+1} + n_\el - 2. 
\end{align*}
\end{remark}

\begin{remark} \label{rem:C_1-cofinite general c5}
Suppose that $ n = 1 $, i.e., the modules $ W_1 , \dots, W_{\N+1} $ are $ C_1 $-cofinite. To follow Huang's method, we must use $ \mathcal{A}_{\el}(u,w_0, \dots, w_{\N+1}; n_1, \dots, n_{\N+1}) $ to express 
\begin{align*}
w_0 \otimes w_1 \otimes \cdots \otimes u_{\alpha_\el - 1,0 } w_\el \otimes  \cdots \otimes w_{\N+1} 
\end{align*}
in terms of lower weight elements in $ T_{K_0} $ up to $ J_{K_0} $. When choosing $ (n_1, \dots, n_{\N+1}) $, we must have $ n_\el = n = 1 $ and we are left to choose $ n_i $, $ i \neq \el $, so that \eqref{eq:inequality condition c5} holds. If we allow for $ \Re(\alpha_j) \in [0,1) $ and $ \wt u \geq 1 $ to be arbitrary, then we must satisfy \eqref{eq:inequality arbitray eigenvalues c5}. Supposing these inequalities hold, we must have $ n_i \leq 0 $ for $ i \neq \el $, which implies the contradiction $ N+2 \leq 2 $. Thus, $ C_1 $-cofiniteness for general twisted modules is not a strong enough assumption to use the present method of finding $ n_i $ satisfying \eqref{eq:inequality condition c5} uniformly for arbitrary eigenvalues. We do not have a counterexample  showing that there exists a product of twisted intertwining operators among $ C_1 $-cofinite twisted modules that does not converge. So, at the moment, we cannot rule out another method for proving convergence.

Note that there may also be special cases where $ n_i $ can be chosen to satisfy \eqref{eq:inequality condition c5} dependent on specific eigenvalues (as opposed to a uniform choice for arbitrary eigenvalues). However, we will not go into a detailed exploration of all such special cases here. 
\end{remark}

\begin{lemma} \label{lem: generator filter C1 case c5}
Let $ K_0 \in \mathbb{R} $. Let $ u \in V_+ $ be a weight-homogeneous simultaneous generalized eigenvector for $ G $ and let $ w_i \in W_i $, $ i = 0, \dots, \N+1 $, be weight-homogeneous with $ \Re(\wt w_0) \leq K_0 $. Let 
\begin{align*}
w^{(\el)} &=  w_0 \otimes w_1 \otimes \cdots \otimes u_{\alpha_\el - 1,0}  w_\el \otimes \cdots \otimes w_{\N+1}, \qquad \text{ for } \el = 1, \dots, \N + 1
\end{align*} 
and denote $ \Re (\wt w^{(\el)}) = s $. Assume $ n = 1 $ and at least $ \N $ of the automorphisms $ g_0, \dots, g_{\N+1} $ only have eigenvalues in $ \mathbb{R}_{>0} $. Then $ w^{(\el)} \in F_s(J_{ K_0}) + F_r(T_{ K_0}) $ for some $ r < s $.
\end{lemma}

\begin{proof}
Following the start of the proof for Lemma~\ref{lem: generator filter c5} and choosing $ n_\el = 1 $, we see that we need to  choose the remaining $ n_i $ so that the inequalities
\begin{equation}
\begin{gathered}
\Re(\alpha_1 + \cdots + \alpha_{\N+1} +  \alpha_\el ) < 2 \wt u + n_1 + \cdots + n_{\N+1} + n_\el - 2,  \\
\Re(\alpha_\el - \alpha_i )< n_\el - n_i , \quad \text{ for } i \neq \el, \\
\Re(\alpha_1 + \cdots + \alpha_{\N+1)} \leq n_1 + \cdots + n_{\N+1},
\end{gathered}
\end{equation}
are satisfied. Say there are $ g_j $ and $ g_k $ possibly with eigenvalues not in $ \mathbb{R}_{>0} $. Then $ \Re(\alpha_i) = 0 $ for all $ i \neq j,k $. We also observe that $ g_j g_k =  \prod_{i\neq j,k} g_i^{-1} $, hence $ \Re(\alpha_j + \alpha_k) = 0 $ or $ 1 $. From this, we see that we can choose $ n_i = 0 $, for all $ i \neq \el $.
\end{proof}

From here on, the cofiniteness conditions, discrete grading conditions, and CFT type assumptions are needed. Recall that these assumptions were stated in Section~\ref{sec:setting c5}.

\begin{lemma}\label{lem: F_r(T) c5}
Let $ K_0 \in \mathbb{R} $ and assume that either 
\begin{enumerate}
\item[(i)] $ n \geq 2 $ and $ g_0, \dots, g_{\N+1} \in G $ are general; or  
\item[(ii)] $ n = 1 $ and at least $ \N $ of the automorphisms $ g_0, \dots, g_{\N+1} $ only have eigenvalues in $ \mathbb{R}_{>0} $.
\end{enumerate} 
Then there exists $ M \in \mathbb{Z} $ such that $ F_r(T_{K_0}) \subseteq F_r(J_{K_0}) + F_M(T_{K_0}) $, for all $ r \in \mathbb{R} $. Furthermore, $ T_{K_0} = J_{K_0} + F_M(T_{K_0}) $ and the $ R $-module $ T_{K_0}/J_{K_0} $ is finitely generated.
\end{lemma} 

\begin{proof}
Let $ n_i = n $ if $ g_i $ has an eigenvalue not in $ \mathbb{R}_{>0} $, and $ n_i = 1 $ if $ g_i $ only has eigenvalues in $ \mathbb{R}_{>0} $. Since $ \dim W_i / C_{n_i}(W_i) < \infty $, for $ i = 1,\dots , \N + 1 $, there is a finite-dimensional subspace $ U_i \subseteq W_i $ such that $ W_i = U_i + C_{n_i}(W_i) $. Hence, we can find some $ M_i \in \mathbb{Z} $ sufficiently large so that any element in $ (W_i)_{[m]} $ is necessarily in $  C_{n_i}(W_i) $ whenever $ \Re(m) > M_i $. Taking $ M = \lceil K_0 \rceil + \sum_{i = 1}^{\N+1} M_i $, we see that any weight-homogeneous element in $ T_{K_0} $ that has real part of its weight larger than $ M $ must be in
\begin{align*}
\sum_{i = 1}^{\N + 1} R \otimes \left(  \coprod_{\substack{ m \in \mathbb{C} \\ \Re(m) \leq K_0}} (W_0)_{[m]} \right) \otimes W_1 \otimes \cdots \otimes C_{n_i}(W_i) \otimes \cdots \otimes W_{\N+1} .
\end{align*}
This space is $ R $-spanned by elements in Lemma~\ref{lem: generator filter c5} or Lemma~\ref{lem: generator filter C1 case c5} since $ V $ is assumed to be of CFT type.

For all $ r \leq M $, we immediately have 
\begin{align*}
F_r(T_{K_0}) \subseteq F_r(T_{K_0}) \subseteq F_M(T_{K_0}) \subseteq F_r(J_{K_0}) + F_M(T_{K_0}) .
\end{align*} 
Since $ T_{K_0} $ is $ R $-spanned by a discretely graded space, we can use induction on $ r \geq M $. Assume that we have some $ s > M $ such that $ F_r(T_{K_0}) \subseteq F_r(J_{K_0}) + F_M(T_{K_0}) $ for all $ r < s $. 

Consider an element $ w \in (T_{K_0})_{(s)} $, which we can assume is an $ R $-linear combination of the elements from Lemma~\ref{lem: generator filter c5} or Lemma~\ref{lem: generator filter C1 case c5}. Hence, we have $ w \in F_{s}(J_{K_0}) + F_r(T_{K_0}) $, for some $ r < s $. By the induction hypothesis, we have $ w \in F_{s}(J_{K_0}) + F_r(T_{K_0}) \subseteq  F_{s}(J_{K_0}) + F_{r}(J_{K_0}) + F_M(T_{K_0}) $. So, we have shown that $ F_s(T_{K_0}) \subseteq  F_{s}(J_{K_0}) + F_M(T_{K_0}) $.

Furthermore, we have
\begin{align*}
T_{K_0} = \bigcup_{r \in \mathbb{R}} F_r(T_{K_0}) \subseteq \bigcup_{r \in \mathbb{R}} \left( F_r(J_{K_0}) + F_M(T_{K_0}) \right) \subseteq \bigcup_{r \in \mathbb{R}} F_r(J_{K_0}) + F_M(T_{K_0}) = J_{K_0} + F_M(T_{K_0}).
\end{align*}
Also, $ J_{K_0} + F_M(T_{K_0}) \subseteq T_{K_0} $, so we have $ T_{K_0} = J_{K_0} + F_M(T_{K_0}) $. So, the image of $ F_M(T_{K_0}) $ in $ T_{K_0}/J_{K_0} $ is all of $ T_{K_0}/J_{K_0} $. Hence, $ T_{K_0}/J_{K_0} $ is finitely generated since $ F_M(T_{K_0}) $ is finitely generated.
\end{proof}

Consider a change of variables from $ (z_1, \dots, z_\N) $ to $ (\zeta_1, \dots, \zeta_\N) $ given by
\begin{equation}
z_j = \sum_{i=1}^\N b_{ij} \zeta_i + \gamma_j,
\end{equation}
for $ B = (b_{ij})_{i,j = 1}^\N \in \GL(\N, \mathbb{C}) $ and $ (\gamma_1, \dots, \gamma_\N) \in \mathbb{C}^\N $.
Then 
\begin{equation}
\frac{\partial}{\partial \zeta_i} = \sum_{j = 1}^\N \frac{\partial z_j}{\partial \zeta_i} \frac{\partial}{\partial z_j} = \sum_{j = 1}^\N b_{ij} \frac{\partial}{\partial z_j} .
\end{equation}
Define the operators $ \widehat{L}_{\zeta_\el} $ on $ W_0 \otimes \cdots \otimes W_{\N+1} $ by 
\begin{equation}
\widehat{L}_{\zeta_\el} (w_0 \otimes \cdots \otimes w_{\N+1}) = \sum_{j = 1}^\N b_{\el j} w_0 \otimes \cdots \otimes L(-1) w_j \otimes \cdots \otimes w_{\N+1} ,
\end{equation}
for each $ \el = 1, \dots, \N $.

For $ w \in T_{K_0} $, we use $ [w] $ to denote the equivalence class $ w + J_{K_0} $ in $ T_{K_0}/J_{K_0} $.

\begin{theorem}  \label{thm: DEs c5}
Assume that either 
\begin{enumerate}
\item[(i)] $ n \geq 2 $ and $ g_0, \dots, g_{\N+1} \in G $ are general; or  
\item[(ii)] $ n = 1 $ and at least $ \N $ of the automorphisms $ g_0, \dots, g_{\N+1} $ only have eigenvalues in $ \mathbb{R}_{>0} $.
\end{enumerate} 
Assume that $ R $ is Noetherian. Then, for each $ (N+2)$-tuple of $ w_i \in W_i $, there exists a system of differential equations
\begin{equation}\label{eq:DE}
\begin{aligned}
& \left(\frac{\partial}{\partial \zeta_\el} \right)^{m_\el} \psi + \sum_{k = 1}^{m_\el} a_{k,\el}(z_1, \dots, z_\N)  \left(\frac{\partial}{\partial \zeta_\el} \right)^{m_\el-k} \psi  = 0, \qquad \el = 1, \dots, \N , 
\end{aligned} 
\end{equation}
for some $ m_{\el} \in \mathbb{Z}_{>0} $ and $ a_{k,\el}(x_1, \dots, x_\N) \in R $, formally satisfied by the series
\begin{align*}
\langle w_0, \mathcal{Y}_1(w_1,z_1) \cdots \mathcal{Y}_\N(w_\N,z_\N) w_{\N+1} \rangle
\end{align*}
in the region $ |z_1| > \dots > | z_\N | > 0 $.
\end{theorem}

\begin{proof}
Choose $ K_0 \in \mathbb{R} $ sufficiently large so that $ w_0 $ is an element of $ \coprod_{\substack{ m \in \mathbb{C} \\ \Re(m) \leq K_0}} (W_0)_{[m]} $. Consider an ascending chain $ (M_j)_{j = 0}^\infty $ of submodules $ M_j $ of $ T_{K_0}/J_{K_0} $ generated by 
\begin{equation}
\{ [\widehat{L}_{\zeta_\el}^k (w_0 \otimes \dots \otimes w_{\N+1})] : k \in \mathbb{Z}, \ 0 \leq k \leq j \} .
\end{equation} 
Since $ R $ is a Noetherian ring and $ T_{K_0}/J_{K_0} $ is a finitely generated $ R $-module, this chain stabilizes. So, there is some $ m_\el \in \mathbb{Z}_{> 0} $ such that $ M_{m_\el} = M_{m_\el-1} $. Hence, $ [\widehat{L}_{\zeta_\el}^{m_\el} (w_0 \otimes \dots \otimes w_{\N+1})] \in M_{m_\el-1} $, and we obtain 
\begin{equation}\label{eq:L(-1)-eq}
\begin{aligned}
&[\widehat{L}_{\zeta_\el}^{m_\el} (w_0 \otimes \dots \otimes w_{\N+1})] + a_{1,\el} (x_1, \dots, x_\N) [\widehat{L}_{\zeta_\el}^{m_\el-1} (w_0 \otimes \dots \otimes w_{\N+1})] + \\
&\quad \cdots + a_{m_\el,\el} (x_1, \dots, x_\N)[w_0 \otimes \dots \otimes w_{\N+1}] = 0.
\end{aligned}
\end{equation}

Since $ J_{K_0} $ is contained in the kernel of $ \phi $, we have the map
\begin{align*}
\overline{\phi} : T_{K_0}/J_{K_0} \to \mathbb{C}[\log x_1, \dots, \log x_\N]\{x_1, \dots, x_\N \}.
\end{align*}
Using this map on \eqref{eq:L(-1)-eq} followed by the $ L(-1) $-derivative property, we obtain
\begin{equation}  
\begin{aligned}
& \left(\sum_{j=1}^\N b_{\el j} \frac{\partial}{\partial x_j} \right)^{m_\el} \langle w_0, \mathcal{Y}_1(w_1,x_1) \cdots \mathcal{Y}_\N(w_\N,x_\N) w_{\N+1} \rangle \\
&\quad + \sum_{k = 1}^{m_\el} \iota_{|x_1| > \dots > |x_\N|} a_{k,\el}(x_1, \dots, x_\N) \left(\sum_{j=1}^\N b_{\el j} \frac{\partial}{\partial x_j} \right)^{m_\el-k} \langle w_0, \mathcal{Y}_1(w_1,x_1) \cdots \mathcal{Y}_\N(w_\N,x_\N) w_{\N+1} \rangle  = 0.
\end{aligned} 
\end{equation}
As these arguments are formal and we have used $ \iota_{|x_1| > \dots > |x_\N|} $ to expand the coefficients, when we evaluate $ x_1 = z_1 , \dots, x_\N = z_\N $ we obtain the differential equations \eqref{eq:DE} formally satisfied by the multi-series $ \langle w_0, \mathcal{Y}_1(w_1,z_1) \cdots \mathcal{Y}_\N(w_\N,z_\N) w_{\N+1} \rangle $ in the region $ |z_1| > \cdots > |z_\N | > 0 $.
\end{proof}

\begin{remark}
Recall that $ R $ is Noetherian when both $ V $ and $ G $ are finitely generated, for which there are many natural examples. Even if $ V $ is not finitely generated, $ R $ will still be Noetherian when $ G $ is finite, which is the important case studied in Section~\ref{sec:convergence}.

If $ R $ is not Noetherian, we can extend $ R $ to the field of fractions $ \Frac(R) $, since $ R $ is an integral domain\footnote{We credit Jishen Du for the idea of extending $ R $ to the field of fractions $ \Frac(R) $.}. Since $ T_{K_0} $ was defined as a freely generated $ R $-module, we can extend it to a $ \Frac(R) $-vector space. The preceding arguments still hold with $ R $ replaced with $ \Frac(R) $ (and some minor changes such as a change in codomain of $ \phi $ that we will not detail here), namely $ \Frac(R) \otimes_R (T_{K_0}/J_{K_0}) $ is a finite-dimensional vector space and we obtain a differential equation similar to that of \eqref{eq:DE} but with the coefficients now appropriate expansions of elements in $ \Frac(R) $. However, we emphasize that these coefficients have unknown domains as they may contain denominators that are arbitrary elements in $ R $. Since the genus-zero correlation functions in an orbifold conformal field theory should be defined on $ M^{\N} $, we would like the coefficients of the differential equation to also be defined on $ M^{\N} $, which naturally happens when $ R $ is Noetherian. For this reason, we still highlight the case where $ R $ is Noetherian.
\end{remark}

\section{Improving coefficients} \label{sec:improve coefficients}
In this section, we improve upon the form of the coefficients of the differential equations derived in the last section. To derive such differential equations with improved coefficients, we will need further Noetherian assumptions on certain subrings $ R^{(c=0)} $ of $ R $. 

We define a $ \mathbb{C} $-gradation on $ R $. Define the degree of $ (x_i - x_j)^\alpha $ to be $ \alpha $. Define the degree of $ \log (x_i - x_j) $ to be zero. The $ \mathbb{C} $-gradation on the ring $ R $ by \emph{degree} is defined by extending the notion of degree from these generators to all monomials.

We now introduce some filtrations to study the types of singularities at certain points $ c = 0 $. Let $ c = \ x_i - x_j $ for $ i,j \in \{ 1, \dots, \N + 1\} $ with $ i < j $. We define a filtration of $ R $, dependent on a choice of $ c = x_i - x_j $, by defining $ F_r^{(c = 0)}(R) $ to be the $ \mathbb{C} $-span of 
\begin{equation}
 \prod_{1 \leq s<t \leq \N+1} (x_s - x_t)^{\alpha_{st}} \prod_{1 \leq s<t \leq \N+1} (\log( x_s - x_t))^{k_{st}} \quad \text{with } \alpha_{st} \in A, \ \Re(\alpha_{ij}) \geq - r, \ k_{st} \in \mathbb{Z}_{\geq 0} ,
\end{equation}
for each $ r \in \mathbb{R}_{\geq 0} $. We define the subring $ R^{(c = 0)} $ to be the subring $ F_0^{(c = 0)}(R) $ of $ R $. Define $ R^{(c = 0)}_{(\deg 0)} $ to be the subspace of degree zero elements in $ R^{(c = 0)} $.

Define $ F_r^{(c = 0)} (T_{K_0}) $, $ r \in \mathbb{R} $, to be the vector subspace of $ T_{K_0} $ spanned by $ f(x_1,\dots, x_{\N}) w_0 \otimes \cdots \otimes w_{\N+1} $ for all $ f(x_1,\dots, x_{\N}) \in F_m^{(c=0)}(R) $, $ m \in \mathbb{R}_{\geq 0} $ and weight-homogeneous $ w_i \in W_i $ with $ \Re(\wt w_0) \leq K_0 $ satisfying $ m + \sigma \leq r $. Here, we continue to use $ \sigma = \Re( \wt w_0 + \cdots + \wt w_{\N+1}) $. Define $ F_r^{(c = 0)} (J_{K_0}) = F_r^{(c = 0)} (T_{K_0}) \cap J_{K_0} $. These are compatible filtrations on the $ R $-modules $ T_{K_0} $ and $ J_{K_0} $ in the sense that $ F_m^{(c=0)}(R) F_r^{(c=0)}(T_{K_0}) \subseteq F_{m+r}^{(c=0)}(T_{K_0}) $ and $ F_m^{(c=0)}(R) F_r^{(c=0)}(J_{K_0}) \subseteq F_{m+r}^{(c=0)}(J_{K_0}) $. Define
\begin{equation}
T_{K_0}^{(c = 0)} = R^{(c = 0)} \otimes \left(  \coprod_{\substack{ m \in \mathbb{C} \\ \Re(m) \leq K_0}} (W_0)_{[m]} \right) \otimes W_1 \otimes \cdots \otimes W_{\N+1},
\end{equation}
with the grading $ T_{K_0}^{(c=0)} = \coprod_{m \in \mathbb{R}} (T_{K_0}^{(c=0)})_{(m)} $ given by the sum of the real parts of the weights for $ W_i $, with trivial grading for $ R^{(c = 0)} $. Then $ F_m^{(c=0)}(R) (T_{K_0}^{(c=0)})_{(r)} \subseteq F_{m+r}^{(c=0)}(T_{K_0}) $. 

\begin{lemma}\label{lem:(c=0) grade shift c5}
Let $ K_0\in \mathbb{R} $. Let $ u \in V_+ $ be a weight-homogeneous simultaneous generalized eigenvector for $ G $ and let $ w_i \in W_i $, $ i = 0, \dots, \N+1 $, be weight-homogeneous vectors with $ \Re(\wt w_0) \leq K_0 $. Let 
\begin{align*}
w^{(\el)} &=  w_0 \otimes w_1 \otimes \cdots \otimes u_{\alpha_\el - n_\el,0}  w_\el \otimes \cdots \otimes w_{\N+1}, \qquad \text{ for } \el = 1, \dots, \N + 1
\end{align*} 
where $ n_\el = n $ if $ g_\el $ has an eigenvalue not in $ \mathbb{R}_{>0} $, and $ n_\el = 1 $ if $ g_\el $ only has  eigenvalues in $ \mathbb{R}_{>0} $. Let $ \Re (\wt w^{(\el)}) = s $. If $ n \geq  2 $, then $ w^{(\el)} \in F_s^{(c = 0)}(J_{ K_0}) + \coprod_{m > 0} F_m^{(c = 0)}(R) (T_{ K_0}^{(c = 0)})_{(s-m)} $.
\end{lemma}

\begin{proof}
Let $ \el = 1, \dots, \N+1 $. By choosing $ (n_1, \dots, n_{\N+1}) \in \mathbb{Z}^{\N+1} $ with $ n_\el = n $ or $ n_\el = 1 $ (depending on $ g_\el $) and satisfying the conditions 
\begin{gather}
\Re(\alpha_1 + \cdots + \alpha_{\N+1} +  \alpha_\el)  < 2 \wt u + n_1 + \cdots + n_{\N+1} + n_\el - 2,  \\
\Re(\alpha_\el - \alpha_i) < n_\el - n_i , \quad \text{ for } i \neq \el, \\
\Re(\alpha_1 + \cdots + \alpha_{\N+1}) \leq n_1 + \cdots + n_{\N+1}, 
\end{gather}
from Lemma~\ref{lem: generator filter c5}, we know $ \mathcal{A}_\el (u,w_0, \dots ,w_{\N+1};n_1, \dots, n_{\N+1}) - w^{(\el)} \in F_{s-m}(T_{K_0}) $ for some $ m > 0 $. We are left to perform a case-by-case check of the requirements that show that there exist generators $ \mathcal{A}_\el(u,w_0, \dots ,w_{\N+1}; n_1, \dots, n_{\N+1}) $ of $ J_{K_0} $ that are in $ F_s^{(c = 0)}(T_{K_0}) $. 

\textbf{Case 1.} Assume that $ c = x_p - x_q $ for some distinct $ p, q \in \{ 1, \dots, \N+1 \} $ with $ p,q \neq \el $. The only terms of $ \mathcal{A}_\el(u,w_0, \dots ,w_{\N+1}; n_1, \dots, n_{\N+1}) $ that we need to analyze are those containing a factor of $ c^{m} $ for some $ m \in \mathbb{C} $ with $ \Re(m) < 0 $. These terms are proportional to 
\begin{gather*}
(x_p - x_q)^{\alpha_q-n_q-k_q}    w_0 \otimes w_1 \otimes \cdots \otimes v_{\alpha_p - n_p + k_1 + \cdots + \widehat{k_p} + \cdots + k_{\N+1},0 } w_p \otimes  \cdots \otimes w_{\N+1} \text{ or } \\
(x_q - x_p)^{\alpha_p-n_p-k_p}    w_0 \otimes w_1 \otimes \cdots \otimes v_{\alpha_q - n_q + k_1 + \cdots + \widehat{k_q} + \cdots + k_{\N+1},0 } w_q \otimes  \cdots \otimes w_{\N+1}
\end{gather*}
for all $ (k_1, \dots, k_{\N+1} ) \in \mathbb{Z}_{\geq 0}^{\N+1} $, up to some element in $ R^{(c=0)} $, for some $ v \in V_{(\wt u)} $ with the same eigenvalues as $ u $.
When $ \Re(\alpha_q)-n_q-k_q < 0 $, we need to satisfy the condition
\begin{align*}
&\sigma + \wt u - \Re(\alpha_p) + n_p - k_1 - \cdots - \widehat{k_p} - \cdots - k_{\N+1} - 1 - \Re(\alpha_q) + n_q + k_q \\
&\quad \leq \sigma + \wt u - \Re(\alpha_\el) + n_\el - 1,
\end{align*}
for all $ (k_1, \dots, k_{\N+1} ) \in \mathbb{Z}_{\geq 0}^{\N+1} $, which is equivalent to
\begin{equation} \label{eq:ineq lpq c5}
\Re(\alpha_\el - \alpha_p - \alpha_q) \leq n_\el - n_p - n_q .
\end{equation}
Note that this is invariant under the exchange of $ p $ and $ q $, so this is the same condition to be satisfied when $ \Re(\alpha_p)-n_p-k_p < 0 $. 

In the case that $ g_\el $ has an eigenvalue not in $ \mathbb{R}_{>0} $, we can let $ \Re(\alpha_\el),\Re(\alpha_p),\Re(\alpha_q) \in [0,1) $ be arbitrary so that a necessary and sufficient condition for \eqref{eq:ineq lpq c5} is
\begin{align*}
n_p + n_q \leq n_\el - 1 .
\end{align*}
Hence, if $ n_\el = n \geq 2 $, we recall by Remark~\ref{rem:choices c5} that we can ensure that all conditions are met by picking either $ n_p = 0 $ or $ n_q = 0 $, and $ n_i = 1 $ for the remaining $ i $.  

In the case that $ g_\el $ only has eigenvalues in $ \mathbb{R}_{>0} $, we can let $ \Re(\alpha_p),\Re(\alpha_q) \in [0,1) $ be arbitrary so that a necessary and sufficient condition for \eqref{eq:ineq lpq c5} is
\begin{align*}
n_p + n_q \leq n_\el .
\end{align*}
By Remark~\ref{rem:choices c5}, we always pick either $ n_p = 0 $ or $ n_q = 0 $, with the other being $ 0 $ or $ 1 $. Either way, the condition is satisfied for $ n_\el = 1 $.

\textbf{Case 2.} Assume that $ c = x_\el - x_p $ for some $ p \in \{ 1, \dots, \N+1 \}$ with $ p \neq \el $. The only terms we need to look at are proportional to 
\begin{gather*}
(x_\el - x_p)^{-\alpha_p+n_p + k_\el}  v_{\alpha_1 + \cdots + \alpha_{\N+1} - n_1 - \cdots - n_{\N+1} - k_1 - \cdots - k_\N,0}^* w_0 \otimes \cdots  \otimes w_{\N+1}, \text{ if } p = \N+1, \\
 (x_\el - x_p)^{-\alpha_p+n_p}  v_{\alpha_1 + \cdots + \alpha_{\N+1} - n_1 - \cdots - n_{\N+1} - k_1 - \cdots - k_\N,0}^* w_0 \otimes \cdots  \otimes w_{\N+1} , \text{ if } p \neq \N+1,\\
 (x_\el - x_p)^{-\alpha_p+n_p}  w_0 \otimes w_1 \otimes \cdots \otimes v_{\alpha_i - n_i + k_1 + \cdots + \widehat{k_i} + \cdots + k_{\N+1},0} w_i \otimes  \cdots \otimes w_{\N+1}, \text{ for } i \neq \el, p,  \\
 (x_\el - x_p)^{-k_p}    w_0 \otimes w_1 \otimes \cdots \otimes v_{\alpha_\el - n_\el + k_1 + \cdots + \widehat{k_\el} + \cdots + k_{\N+1},0} w_\el \otimes  \cdots \otimes w_{\N+1} , \text{ or } \\
 (x_\el - x_p)^{-\alpha_p+n_p + \alpha_\el-n_\el-k_\el}    w_0 \otimes w_1 \otimes \cdots \otimes v_{\alpha_p - n_p + k_1 + \cdots + \widehat{k_p} + \cdots + k_{\N+1},0 } w_p \otimes  \cdots \otimes w_{\N+1} ,
\end{gather*}
for all $ (k_1, \dots, k_{\N+1} ) \in \mathbb{Z}_{\geq 0}^{\N+1} $, up to some element in $ R^{(c=0)} $, for some $ v \in V_{(\wt u)} $ with the same eigenvalues as $ u $. 
If the powers of $ c = x_\el - x_p $ have negative real part, then the inequalities we need to check are
\begin{align*}
&\sigma - \wt u + \Re(\alpha_1 + \cdots + \alpha_{\N+1}) - n_1 - \cdots - n_{\N+1} - k_1 - \cdots - k_\N + 1 + \Re(\alpha_p) - n_p - k_\el \\
	&\quad \leq  \sigma + \wt u - \Re(\alpha_\el) + n_\el - 1, \\
&\sigma - \wt u + \Re(\alpha_1 + \cdots + \alpha_{\N+1}) - n_1 - \cdots - n_{\N+1} - k_1 - \cdots - k_\N + 1 + \Re(\alpha_p) - n_p \\
	&\quad \leq  \sigma + \wt u - \Re(\alpha_\el) + n_\el - 1, \\
&\sigma + \wt u - \Re(\alpha_i) + n_i - k_1 - \cdots - \widehat{k_i} - \cdots - k_{\N+1} - 1 + \Re(\alpha_p) - n_p \\
	&\quad \leq  \sigma + \wt u - \Re(\alpha_\el) + n_\el - 1, 
\end{align*}
and
\begin{align*}
&\sigma + \wt u - \Re(\alpha_\el) + n_\el - k_1 - \cdots - \widehat{k_\el} - \cdots - k_{\N+1} - 1 + k_p \\
	&\quad \leq \sigma + \wt u - \Re(\alpha_\el) + n_\el - 1 , \\
&\sigma + \wt u - \Re(\alpha_p) + n_p - k_1 - \cdots - \widehat{k_p} - \cdots - k_{\N+1} - 1 + \Re(\alpha_p) - n_p -  \Re(\alpha_\el) + n_\el + k_\el  \\
	&\quad \leq \sigma + \wt u - \Re(\alpha_\el) + n_\el - 1,
\end{align*}
for all $ (k_1, \dots, k_{\N+1} ) \in \mathbb{Z}_{\geq 0}^{\N+1} $. The last two inequalities are immediately satisfied for any choice of $ ( n_1, \dots, n_{\N+1}) $. 

Consider the case that $ g_\el $ has an eigenvalue not in $ \mathbb{R}_{>0} $. If we pick $ n_p = 1 $ so that $  -\Re(\alpha_p)+n_p > 0 $, we do not need to check the first three inequalities. Hence, when $  n_\el = n \geq 2 $, we can still pick $ n_i = 1 $ for all $ i \neq \el $.  

Consider the case that $ g_\el $ only has eigenvalues in $ \mathbb{R}_{>0} $. Then we choose $ n_i $ as in Lemma~\ref{lem: generator filter c5} so that $ \Re(\alpha_i) \leq n_i $. Then $ - \Re(\alpha_i) + n_i \geq 0 $ and the first three kinds of terms do not have powers of $ c $ with negative real parts. Hence, there are no inequalities to check. 

Thus, for all $ c = x_i - x_j $, we have $ w^{(\el)} \in F_s^{(c = 0)}(J_{ K_0}) + \coprod_{m > 0} F_m^{(c = 0)}(R) (T_{ K_0}^{(c = 0)})_{(s-m)} $ when $ n \geq  2 $.
\end{proof}

\begin{lemma}\label{lem:(c=0) grade shift C1 case c5}
Let $ K_0 \in \mathbb{R} $. Let $ u \in V_+ $ be a weight-homogeneous simultaneous generalized eigenvector for $ G $ and let $ w_i \in W_i $, $ i = 0, \dots, \N+1 $, be weight-homogeneous with $ \Re(\wt w_0) \leq K_0 $. Let 
\begin{align*}
w^{(\el)} &=  w_0 \otimes w_1 \otimes \cdots \otimes u_{\alpha_\el - 1,0}  w_\el \otimes \cdots \otimes w_{\N+1}, \qquad \text{ for } \el = 1, \dots, \N + 1
\end{align*} 
and denote $ \Re (\wt w^{(\el)}) = s $. Assume $ n = 1 $ and at least $ \N $ of the automorphisms $ g_0, \dots, g_{\N+1} $ only have eigenvalues in $ \mathbb{R}_{>0} $. Then $ w^{(\el)} \in F_s^{(c = 0)}(J_{ K_0}) + \coprod_{m > 0} F_m^{(c = 0)}(R) (T_{ K_0}^{(c = 0)})_{(s-m)} $. 
\end{lemma}

\begin{proof}
Using the beginning of the proof for Lemma~\ref{lem:(c=0) grade shift c5}, we check the case that $ n_\el = 1 $ and $ n_i = 0 $ for all $ i \neq \el $. In Case 1, inequality \eqref{eq:ineq lpq c5} still holds. In Case 2, we cannot pick $ n_p = 1 $, so we must check 
\begin{align*}
&\sigma - \wt u + \Re(\alpha_1 + \cdots + \alpha_{\N+1} )- n_1 - \cdots - n_{\N+1} - k_1 - \cdots - k_\N + 1 + \Re(\alpha_p) - n_p - k_\el \\
	&\quad \leq  \sigma + \wt u - \Re(\alpha_\el) + n_\el - 1, \text{ if } p = \N+1\\
&\sigma - \wt u + \Re(\alpha_1 + \cdots + \alpha_{\N+1}) - n_1 - \cdots - n_{\N+1} - k_1 - \cdots - k_\N + 1 + \Re(\alpha_p) - n_p \\
	&\quad \leq  \sigma + \wt u - \Re(\alpha_\el) + n_\el - 1, \text{ if } p \neq \N+1 \\
&\sigma + \wt u - \Re(\alpha_i) + n_i - k_1 - \cdots - \widehat{k_i} - \cdots - k_{\N+1} - 1 + \Re(\alpha_p) - n_p \\
	&\quad \leq  \sigma + \wt u - \Re(\alpha_\el) + n_\el - 1, \text{ for } i \neq \el,p ,
\end{align*}
for all $ (k_1, \dots, k_{\N+1} ) \in \mathbb{Z}_{\geq 0}^{\N+1} $. These are satisfied if
\begin{gather*}
 \Re(\alpha_1 + \cdots + \alpha_{\N+1} + \alpha_p + \alpha_\el) \leq  2 , \\
  \Re(-\alpha_i + \alpha_p + \alpha_\el) \leq  1.
\end{gather*}
Say there are $ g_j $ and $ g_k $ possibly with eigenvalues not in $ \mathbb{R}_{>0} $. Then $ \Re(\alpha_i) = 0 $ for all $ i \neq j,k $, and $ \Re(\alpha_j + \alpha_k) = 0 $ or $ 1 $. From this, we see that these inequalities are met. 
\end{proof}

\begin{lemma}\label{lem: F_r(T) 2}
Let $ K_0 \in \mathbb{R} $ and assume that either 
\begin{enumerate}
\item[(i)] $ n \geq 2 $ and $ g_0, \dots, g_{\N+1} \in G $ are general; or  
\item[(ii)] $ n = 1 $ and at least $ \N $ of $ g_0, \dots, g_{\N+1} $ only have eigenvalues in $ \mathbb{R}_{>0} $.
\end{enumerate} 
Then there exists $ M \in \mathbb{Z} $ such that $ F_r^{(c = 0)}(T_{K_0}) \subseteq F_r^{(c = 0)}(J_{K_0}) + F_M(T_{K_0}) $, for all $ r \in \mathbb{R} $. 
\end{lemma} 

\begin{proof}
Let $ n_i = n $ if $ g_i $ has an eigenvalue not in $ \mathbb{R}_{>0} $, and $ n_i = 1 $ if $ g_i $ only has eigenvalues in $ \mathbb{R}_{>0} $. Since $ \dim W_i / C_{n_i}(W_i) < \infty $, for $ i = 1,\dots , \N + 1 $, there is a finite-dimensional subspace $ U_i \subseteq W_i $ such that $ W_i = U_i + C_{n_i}(W_i) $. Hence, we can find some $ M_i \in \mathbb{Z} $ sufficiently large so that any element in $ (W_i)_{[m]} $ is necessarily in $  C_{n_i}(W_i) $ whenever $ \Re(m) > M_i $. Taking $ M = \lceil K_0 \rceil + \sum_{i = 1}^{\N+1} M_i $, we see that any weight-homogeneous element in $ T_{K_0} $ that has real part of its weight larger than $ M $ must be in
\begin{align*}
\sum_{i = 1}^{\N + 1} R \otimes \left(  \coprod_{\substack{ m \in \mathbb{C} \\ \Re(m) \leq K_0}} (W_0)_{[m]} \right) \otimes W_1 \otimes \cdots \otimes C_{n_i}(W_i) \otimes \cdots \otimes W_{\N+1} .
\end{align*}
This space is $ R $-spanned by elements in Lemma~\ref{lem:(c=0) grade shift c5} or Lemma~\ref{lem:(c=0) grade shift C1 case c5} since $ V $ is assumed to be of CFT type.

For all $ r \leq M $, we immediately have 
\begin{align*}
F_r^{(c = 0)}(T_{K_0}) \subseteq F_r(T_{K_0}) \subseteq F_M(T_{K_0}) \subseteq F_r^{(c = 0)}(J_{K_0}) + F_M(T_{K_0}) .
\end{align*} 
Since $ T_{K_0} $ is $ R $-spanned by a discretely graded space, we can use induction on the real parts $ s \geq M $ of the weights of $ T_{K_0} $. Assume that we have some real part $ s > M $ of a weight of $ T_{K_0} $ such that $ F_r^{(c = 0)}(T_{K_0}) \subseteq F_r^{(c = 0)}(J_{K_0}) + F_M(T_{K_0}) $ for all $ r \leq t $, where $ t $ is the highest real part of a weight of $ T_{K_0} $ that is strictly below $ s $. 

Consider an element $ w \in (T_{K_0}^{(c = 0)})_{(s)} $, which we can assume is an $ R^{(c = 0)} $-linear combination of the elements from Lemma~\ref{lem:(c=0) grade shift c5} or Lemma~\ref{lem:(c=0) grade shift C1 case c5}. Hence, we have $ w \in F_{s}^{(c = 0)}(J_{K_0}) + \coprod_{m > 0} F^{(c=0)}_{m}(R) (T_{K_0}^{(c = 0)})_{(s - m)} $. Since $ s - m < s $ when $ m >0 $, we have $ (T_{K_0}^{(c=0)})_{(s-m)} \subseteq F^{(c=0)}_{s-m}(J_{K_0}) + F_M(T_{K_0}) $ by the induction hypothesis. So, $ F^{(c=0)}_{m}(R) (T_{K_0}^{(c=0)})_{(s-m)} \subseteq F^{(c=0)}_{s}(J_{K_0}) + F_M(T_{K_0}) $. Thus, we have shown that $ (T_{K_0}^{(c = 0)})_{(s)} \subseteq F^{(c=0)}_{s}(J_{K_0}) + F_M(T_{K_0}) $. Let $ r \leq s $. Since we have $ (T_{K_0}^{(c = 0)})_{(\sigma)} \subseteq F^{(c=0)}_{\sigma}(J_{K_0}) + F_M(T_{K_0}) $ for all $ \sigma \leq r $, we see that $ F^{(c=0)}_{m}(R)(T_{K_0}^{(c = 0)})_{(\sigma)} \subseteq F^{(c=0)}_{m + \sigma}(J_{K_0}) + F_M(T_{K_0}) \subseteq F^{(c=0)}_{r}(J_{K_0}) + F_M(T_{K_0}) $ for all $ m \geq 0 $ with $ m + \sigma \leq r $. Thus, $ F_r^{(c=0)}(T_{K_0}) \subseteq  F^{(c=0)}_{r}(J_{K_0}) + F_M(T_{K_0}) $ for all $ r \leq s $.
\end{proof}

\begin{lemma}\label{lem:sigma+S c5}
Let $ K_0 \in \mathbb{R} $ and assume that either 
\begin{enumerate}
\item[(i)] $ n \geq 2 $ and $ g_0, \dots, g_{\N+1} \in G $ are general; or  
\item[(ii)] $ n = 1 $ and at least $ \N $ of $ g_0, \dots, g_{\N+1} $ only have eigenvalues in $ \mathbb{R}_{>0} $.
\end{enumerate} 
For any $ s \in [0,1) $, there exists $ S \in \mathbb{R} $ such that $ s + S \in \mathbb{Z}_{>0} $, and for any weight-homogeneous vectors $ w_i \in W_i $ with $ \Re(\wt w_0) \leq K_0 $, and any $ \mathcal{W}_1 \in F_\sigma^{(c = 0)}(J_{K_0}) $, $ \mathcal{W}_2 \in  F_M(T_{K_0}) $ satisfying $ \sigma \in s + \mathbb{Z} $ and
\begin{align*}
w_0 \otimes \cdots \otimes w_{\N+1} = \mathcal{W}_1 + \mathcal{W}_2,
\end{align*} we have
\begin{align*}
c^{\sigma + S} \mathcal{W}_2 \in T_{K_0}^{(c = 0)}.
\end{align*}
\end{lemma}

\begin{proof}
Pick $ S \in \mathbb{R} $ with $ s + S \in \mathbb{Z}_{> 0} $ to be sufficiently large so that $ (T_{K_0})_{(r)} = 0 $ when $ r < - S $. Let $ w_i \in W_i $ be weight-homogeneous vectors with $ \sigma \in s + \mathbb{Z} $ and let $ \mathcal{W}_1 \in F_\sigma^{(c = 0)}(J_{K_0}), \mathcal{W}_2 \in  F_M(T_{K_0}) $ satisfying
\begin{align*}
w_0 \otimes \cdots \otimes w_{\N+1} = \mathcal{W}_1 + \mathcal{W}_2.
\end{align*}
Then $ \sigma + S \in \mathbb{Z} $ and $ \sigma \geq - S $, so $ \sigma + S \in \mathbb{Z}_{\geq 0} $. Then
\begin{align*}
c^{\sigma + S} \mathcal{W}_2 = c^{\sigma + S} w_0 \otimes \cdots \otimes w_{\N+1} - c^{\sigma + S} \mathcal{W}_1 \in T_{K_0}^{(c=0)} + c^{\sigma + S} F_\sigma^{(c = 0)}(J_{K_0}) \subseteq T_{K_0}^{(c=0)}.
\end{align*} 
To show the last inclusion, consider $ f(x_1,\dots,x_{\N}) \widetilde{w}_0 \otimes \cdots \otimes \widetilde{w}_{\N+1} \in F_\sigma^{(c = 0)}(T_{K_0}) $ with $ f(x_1, \dots, x_\N) \in F^{(c=0)}_m(R) $ and $ \Re(\wt \widetilde{w}_0 \otimes \cdots \otimes \widetilde{w}_{\N+1} ) = \widetilde{\sigma} $. Then
\begin{align*}
m - S - \sigma \leq m + \widetilde{\sigma} - \sigma  \leq \sigma - \sigma = 0 
\end{align*}
implies that $ f(x_1,\dots,x_\N) c^{\sigma + S} \in R^{(c=0)} $.
\end{proof}

\begin{theorem}  \label{thm: DEs c c5}
Assume that either 
\begin{enumerate}
\item[(i)] $ n \geq 2 $ and $ g_0, \dots, g_{\N+1} \in G $ are general; or  
\item[(ii)] $ n = 1 $ and at least $ \N $ of $ g_0, \dots, g_{\N+1} $ only have eigenvalues in $ \mathbb{R}_{>0} $.
\end{enumerate} 
Consider a change of variables from $ (z_1, \dots, z_\N) $ to $ (\zeta_1, \dots, \zeta_\N) $ given by
\begin{equation}
z_j = \sum_{i=1}^\N b_{ij} \zeta_i + \gamma_j,
\end{equation}
for $ B = (b_{ij})_{i,j = 1}^\N \in \GL(\N, \mathbb{C}) $ and $ (\gamma_1, \dots, \gamma_\N) \in \mathbb{C}^\N $. Let $ c = \ x_i - x_j $ for $ i < j \in \{ 1, \dots, \N+1 \} $. Assume that $ R^{(c=0)} $ is Noetherian. Then, for each $ (N+2)$-tuple of weight-homogeneous $ w_i \in W_i $, there exists a system of differential equations
\begin{equation}\label{eq:DE with coefficients c5}
\begin{aligned}
& (c(z_1, \dots, z_\N) )^{m_\el} \left(\frac{\partial}{\partial \zeta_\el} \right)^{m_\el} \psi + \sum_{k = 1}^{m_\el}  d^{(c)}_{k,\el}(z_1, \dots, z_\N) (c(z_1, \dots, z_\N) )^{m_\el-k} \left(\frac{\partial}{\partial \zeta_\el} \right)^{m_\el-k} \psi  = 0,
\end{aligned} 
\end{equation}
for $ \el = 1, \dots, \N $, for some $ m_{\el} \in \mathbb{Z}_{>0} $ and $ d^{(c)}_{k,\el}(x_1, \dots, x_\N) \in R^{(c=0)}_{(\deg 0)} $, satisfied formally by 
\begin{align*}
\langle w_0, \mathcal{Y}_1(w_1,z_1) \cdots \mathcal{Y}_\N(w_\N,z_\N) w_{\N+1} \rangle
\end{align*}
in the region $ |z_1| > \dots > | z_\N | > 0 $.
\end{theorem}

\begin{proof}
Let $ w_i \in W_i $ be of homogeneous weight. Pick $ K_0 = \Re( \wt w_0) $. Based on the change of variables $ (\zeta_1, \dots, \zeta_\N) $, recall the operators $ \widehat{L}_{\zeta_\el} $, for each $ \el = 1, \dots, \N $, defined by
\begin{align*}
\widehat{L}_{\zeta_\el} (\overline{w}_0 \otimes \cdots \otimes \overline{w}_{\N+1}) = \sum_{j = 1}^\N b_{\el j} \overline{w}_0 \otimes \cdots \otimes L(-1) \overline{w}_j \otimes \cdots \otimes \overline{w}_{\N+1},
\end{align*}
for all $ \overline{w}_i \in W_i $, extended by linearity to all $ W_0 \otimes \cdots \otimes W_{\N+1} $. 

Observe that $ \widehat{L}_{\zeta_\el} $ increases the weight by one, and $ T_{K_0} $ is closed under $ \widehat{L}_{\zeta_\el} $. Hence, by Lemma~\ref{lem: F_r(T) 2}, there exist $ \mathcal{W}_1^{(k)} \in F_{\sigma + k}^{(c=0)}(J_{K_0}) $ and $ \mathcal{W}_2^{(k)} \in F_M(T_{K_0}) $ such that
\begin{align*}
\widehat{L}_{\zeta_\el}^k (w_0 \otimes \dots \otimes w_{\N+1}) = \mathcal{W}_1^{(k)}  + \mathcal{W}_2^{(k)},
\end{align*} 
for each $ k \in \mathbb{Z}_{\geq 0} $.

By Lemma~\ref{lem:sigma+S c5}, there exists $ S \in \mathbb{R} $ such that $ \sigma + S \in \mathbb{Z}_{>0} $ and $
c^{\sigma + k + S} \mathcal{W}_2^{(k)} \in T_{K_0}^{(c=0)} $. So,
\begin{align*}
c^{\sigma + k + S} \mathcal{W}_2^{(k)} \in \coprod_{r \leq M} (T_{K_0}^{(c=0)})_{(r)}.
\end{align*}
Now we consider the finitely generated $ R^{(c=0)} $-module $ \coprod_{r \leq M} (T_{K_0}^{(c=0)})_{(r)} $. Consider the ascending chain $ (M_j)_{j=0}^\infty $ of submodules $ M_j $ generated by
\begin{align*}
\{c^{\sigma + k + S} \mathcal{W}_2^{(k)} : k \in \mathbb{Z}, 0 \leq k \leq j \}.
\end{align*} 
Since $ R^{(c=0)} $ is Noetherian, there exists $ m_\el \in \mathbb{Z}_{> 0} $ and there exist $ d^{(c)}_{1,\el}(x_1, \dots, x_\N) $, $\dots$, $d^{(c)}_{m_\el,\el}(x_1,\dots,x_\N) \in R^{(c=0)} $ such that
\begin{equation} \label{eq: W_2 relation c5}
c^{\sigma + m_\el + S} \mathcal{W}_2^{(m_\el)} + \sum_{k=0}^{m_\el-1} d^{(c)}_{m_\el-k,\el}(x_1, \dots, x_\N) c^{\sigma + k + S} \mathcal{W}_2^{(k)} = 0.
\end{equation}
Since the elements $ \mathcal{W}_1^{(j)} $ are in the kernel of $ \phi $ (recall \eqref{eq:correlation map} and its restriction to $ T_{K_0} $), we can divide \eqref{eq: W_2 relation c5} by $ c^{\sigma + S} $ and apply $ \phi $ to obtain
\begin{align*}
&\phi \left( c^{m_\el} \widehat{L}_{\zeta_\el}^{m_\el} (w_0 \otimes \dots \otimes w_{\N+1}) \right) \\
&\quad + \sum_{k =1}^{m_\el} \phi \left( d^{(c)}_{k,\el}(x_1, \dots, x_\N) c^{m_\el -k}  \widehat{L}_{\zeta_\el}^{m_\el - k} (w_0 \otimes \dots \otimes w_{\N+1}) \right) = 0 ,
\end{align*}
which is the formal version of the desired differential equation \eqref{eq:DE with coefficients c5} in the desired region. When $ T_{K_0} $ is given the grading
\begin{align*}
&\deg \left( f(x_1, \dots, x_\N) \overline{w}_0 \otimes \cdots \otimes \overline{w}_{\N+1} \right) \\
&\quad := - \deg f(x_1, \dots, x_\N) - \wt \overline{w}_0 + \wt \overline{w}_1 + \cdots + \wt \overline{w}_{\N+1},
\end{align*}
and the generators for $ J_{K_0} $ are checked to be homogeneous with respect to this new grading, we can see that $ d^{(c)}_{k,\el}(x_1, \dots, x_\N) $ can be chosen to be elements of zero degree.
\end{proof}

We can immediately  extend the previous theorem to the case where $ c $ is any non-zero element in $ \Span_{\mathbb{C}}\{x_1, \dots, x_\N\} $. If $ c $ is a scalar multiple of $ x_i $ or $ x_i - x_j $, the result follows from Theorem~\ref{thm: DEs c c5} by rescaling $ c $. Otherwise, no negative powers of $ c $ occur in $ R $, so we can define $ R^{(c =0)} = R $.  The grading
\begin{align*}
&\deg \left( f(x_1, \dots, x_\N) \overline{w}_0 \otimes \cdots \otimes \overline{w}_{\N+1} \right) \\
&\quad := - \deg f(x_1, \dots, x_\N) - \wt \overline{w}_0 + \wt \overline{w}_1 + \cdots + \wt \overline{w}_{\N+1}
\end{align*}
on $ T_{K_0} $ can again be used to show that the coefficients $ a_{k,\el}(x_1, \dots, x_\N) $ in the differential equation \eqref{eq:DE} have degree $ -k $. We then multiply \eqref{eq:DE} by $ c^{m_\el} $ to obtain our desired differential equation.

The following special case of the previous theorem will be used to prove the convergence property in the next section.
\begin{corollary}
Assume that either 
\begin{enumerate}
\item[(i)] $ n \geq 2 $ and $ g_0, \dots, g_{\N+1} \in G $ are general; or  
\item[(ii)] $ n = 1 $ and at least $ \N $ of $ g_0, \dots, g_{\N+1} $ only have eigenvalues in $ \mathbb{R}_{>0} $.
\end{enumerate}
Assume that $ R^{(x_{\N}=0)} $ is Noetherian. Then, for each $ (N+2)$-tuple of weight-homogeneous $ w_i \in W_i $, there exists a differential equation
\begin{equation}\label{eq:DE for convergence c5}
\begin{aligned}
\sum_{k = 0}^{m_\N} a_{k,\N}(z_1, \dots, z_\N) \left(z_\N \frac{\partial}{\partial z_\N} \right)^{k} \psi = 0,
\end{aligned} 
\end{equation}
with $ a_{k,\N}(x_1, \dots, x_\N) \in R^{(x_\N = 0)} $ and $ a_{m_\N, \N}(x_1, \dots, x_\N) = 1 $, formally satisfied by the series 
\begin{align*}
\langle w_0, \mathcal{Y}_1(w_1,z_1) \cdots \mathcal{Y}_\N(w_\N,z_\N) w_{\N+1} \rangle
\end{align*}
in the region $ |z_1| > \cdots > |z_\N| > 0 $.
\end{corollary}

\section{Convergence and regular singularities for finite abelian groups} \label{sec:convergence}
From here on, we will \textit{assume that $ G $ is a finite abelian group}. The additive group $ A $ of eigenvalues (recall  Section~\ref{sec:setting c5}) then becomes $ \frac{1}{t_0} \mathbb{Z} $, for some $ t_0 \in \mathbb{Z}_{>0} $ that divides $ |G| $. From this, we see that $ R $ and $ R^{(c=0)} $ are automatically Noetherian for all non-zero $ c \in \Span_{\mathbb{C}}\{x_1,\dots,x_\N\} $, and we can omit the logarithmic generators from $ R $ and $ R^{(c=0)} $. Also observe that an automorphism with only positive real eigenvalues must be the identity of $ V $. As a result, the condition ``$ \N $ of $ g_0, \dots, g_{\N+1} $ only have eigenvalues in $ \mathbb{R}_{>0} $" becomes ``two of $ g_0, \dots, g_{\N+1} $ are $ g $ and $ g^{-1} $ for some $ g \in \Aut V $ of finite order, and the remaining automorphisms are $ \id_V $".

\begin{remark}
We have not yet used the cofiniteness assumptions~\eqref{eq:cofinite conditions 2} nor discrete grading assumptions for the intermediate modules $ \W_1, \dots ,\W_{N-1} $; only the cofiniteness assumptions \eqref{eq:cofinite conditions 1} for the in-state modules have been used to derive the differential equations above. The induction step in the proof of the following theorem will require the cofiniteness assumptions \eqref{eq:cofinite conditions 2} when proving convergence. Also note that at no point do we need to assume any cofiniteness conditions for the module $ W_0 $. 
\end{remark}

In Theorem~\ref{thm: DEs c c5}, we have shown the existence of the differential equations for products of twisted intertwining operators among discretely graded twisted $ V $-modules satisfying the cofiniteness conditions \eqref{eq:cofinite conditions 1} and \eqref{eq:cofinite conditions 2} and twisted by elements in $ G $. We will use these differential equations to show that such twisted intertwining operators satisfy the convergence property for products of $ \N $ twisted intertwining operators, and the convergence and extension property for products of two twisted intertwining operators. 

\begin{lemma} \label{lem:form of series}
Let $ W_0, \dots, W_{\N+1} $, $ \widetilde{W}_0, \dots, \widetilde{W}_\N $, $ \mathcal{Y}_1, \dots, \mathcal{Y}_\N $ be as in Section~\ref{sec:setting c5}. Assume that $ G $ is a finite abelian group and let $ t = |G| $. Then the multi-series 
\begin{equation*}
\langle w_0, \mathcal{Y}_1(w_1,x_1) \cdots \mathcal{Y}_\N(w_\N,x_\N) w_{\N+1} \rangle
\end{equation*}
is of the form
\begin{equation}\label{eq:form of series}
\sum_{s_1 \in S_1} \sum_{k_1 = 0}^{K_1} x_1^{s_1} (\log x_1)^{k_1} \cdots \sum_{s_\N \in S_\N} \sum_{k_\N = 0}^{K_\N} x_\N^{s_\N} (\log x_\N)^{k_\N}  \sum_{n_1, \dots, n_\N \in \frac{1}{t}\mathbb{Z} } a_{n_1, \dots, n_\N}^{(s_1,\dots, s_\N; k_1, \dots, k_\N)} x_1^{n_1} \cdots x_\N^{n_\N}
\end{equation}
for some $ a_{n_1, \dots, n_\N}^{(s_1,\dots, s_\N; k_1, \dots, k_\N)} \in \mathbb{C} $, finite sets $ S_1, \dots, S_\N \subseteq \mathbb{C} $ with distinct $ s,s' \in S_i $ satisfying $ s-s' \notin \frac{1}{t}\mathbb{Z} $, and $ K_1, \dots, K_\N \in \mathbb{Z}_{\geq 0} $.
\end{lemma}

\begin{proof}
Let $ W $ be a $ C_m $-cofinite discretely graded $ g $-twisted $ V $-module, for an automorphism $ g \in G $ and $ m \in \mathbb{Z}_{>0} $. Then $ W = U + C_m(W) $ for some finite-dimensional space $ U $. We can enlarge $ U $ so that it has a basis of weight-homogeneous vectors, say with finitely many weights $ r_1, \dots, r_p \in \mathbb{C}$. Let $ w \in W $ be a homogeneous vector. Then $ w $ is a linear combination of vectors in $ U $ and vectors of the form $ u_{\alpha - m,0} \overline{w} $ for some weight homogeneous $ u \in V_+^{[\alpha]} $ with $ \alpha = 0, \frac{1}{t}, \dots, \frac{t-1}{t} $, and some weight homogeneous $ \overline{w} \in W $. We have
\begin{equation*}
\Re (\wt \overline{w}) = \Re (\wt w )- \wt u + \alpha - m + 1 \leq \Re(\wt w) - \frac{1}{t}.
\end{equation*} 
Since $ W $ is discretely graded and hence lower bounded, this repeated argument will eventually terminate. This shows that the weights of $ W $ are contained in $ \{ r_1, \dots, r_p \} + \frac{1}{t} \mathbb{Z}_{\geq 0} $. 

Furthermore, since $ g $ has finite order, it is semisimple, so the twisted vertex operator $ Y_W(u,x) $ has no logarithms. 
Hence, the usual identity $ [L(0), Y_W(u,x)] = Y_W(L(0) u ,x) + x \frac{d}{dx} Y_W(u,x) $ implies that $ [L(0), u_{m,0}] = (\wt u - m - 1) u_{m,0} $. So, $ [L(0) - L(0)_S, u_{m,0}] = 0 $, where $ L(0)_S $ is the semisimple part of $ L(0) $. Since $ U $ is a  finite-dimensional space, there exists some $ J \in \mathbb{Z}_{\geq 0} $ such that $ (L(0) - L(0)_S)^J U = 0 $. Since $ U $ generates all of $ W $ and $ [L(0) - L(0)_S, u_{m,0}] = 0 $ for each $ u \in V $, we have $ (L(0) - L(0)_S)^J W = 0 $ for the same $ J $.

Let $ i = 2, \dots, \N $. Let $ w_i \in W_i $, $ \widetilde{w}_i \in \widetilde{W}_i $ be weight-homogeneous. Then the summands of
\begin{align*}
\mathcal{Y}_i(w_i, x) \widetilde{w}_i = \sum_{k = 0}^K \sum_{n \in \mathbb{C}} (\mathcal{Y}_i)_{n,k}(w_i) \widetilde{w}_i x^{-n-1} (\log x)^k,
\end{align*}
must have $ \wt ( (\mathcal{Y}_i)_{n,k}(w_i) \widetilde{w}_i ) = \wt w_i - n - 1 + \wt \widetilde{w}_i $ being a weight of $ \widetilde{W}_{i-1} $ if $ (\mathcal{Y}_i)_{n,k}(w_i) \widetilde{w}_i $ is to be non-zero. Hence, $ n $ can only take values in $ \{ h_1, \dots, h_l \} + \frac{1}{t}\mathbb{Z} $, for some finitely many numbers $ h_1, \dots, h_l \in \mathbb{C}$. We can further assume that no two elements of $ h_1, \dots, h_l $ are congruent modulo $ \frac{1}{t} \mathbb{Z} $. This set of numbers is dependent on $ W_i , \widetilde{W}_i , \widetilde{W}_{i-1} $, and in fact independent of $ w_i $ and $  \widetilde{w}_i $. Following Remark~3.24 of \cite{HLZ2}, the upper bound $ K $ is independent of $ w_i $ and $  \widetilde{w}_i $ if we can show that the operators $ L(0) - L(0)_S $ on $ W_i $, $ \W_i $ and $ \W_{i-1} $ are each nilpotent. Indeed, we have shown this sufficient condition above.

Since we do not assume cofiniteness conditions on $ W_0 $, we consider the remaining case $ i = 1 $ separately. Let $ w_0 \in W_0 $, $ w_1 \in W_1 $, $ \widetilde{w}_1 \in \widetilde{W}_1 $ be weight-homogeneous. We see that the summands of 
\begin{align*} 
\langle w_0, \mathcal{Y}_1(w_1, x) \widetilde{w}_1 \rangle = \sum_{k = 0}^K \sum_{n \in \mathbb{C}} \langle w_0, (\mathcal{Y}_1)_{n,k}(w_1) \widetilde{w}_1 \rangle x^{-n-1} (\log x)^k
\end{align*}
must have $ \wt ( (\mathcal{Y}_1)_{n,k}(w_1) \widetilde{w}_1 ) = \wt w_1 - n - 1 + \wt \widetilde{w}_1 = \wt w_0 $ if $ \langle w_0, (\mathcal{Y}_1)_{n,k}(w_1) \widetilde{w}_1 \rangle $ is to be non-zero. Hence, we still find that $ n $ can only take values in $ \{ h_1, \dots, h_l \} + \frac{1}{t}\mathbb{Z} $, for some finitely many numbers $ h_1, \dots, h_l \in \mathbb{C}$. This set of numbers is now dependent on $ w_0 \in W_0 $, but still independent of $ w_1 $ and $ \widetilde{w}_1 $. 

We can again show that the upper bound $ K $ is independent of $ w_1 $ and $  \widetilde{w}_1 $.  Since we do not assume cofiniteness conditions on $ W_0 $, we cannot assume there is some $ J \in \mathbb{Z}_{\geq 0} $ such that $ (L(0) - L(0)_S)^J W_0 = 0 $. However, Proposition~3.20(a) of \cite{HLZ2} still provides the upper bound $ K $ given some $ J \in \mathbb{Z}_{\geq 0} $ such that $ (L(0) - L(0)_S)^J w_0 = 0 $. This upper bound is now dependent on $ w_0 \in W_0 $, but still independent of $ w_1 $ and $ \widetilde{w}_1 $.

Since the modules $ W_1, \dots, W_{\N+1}, \W_1, \dots, \W_{\N-1} $ are assumed to be discretely graded and $ C_1 $- or $ C_n $-cofinite, and there are only finitely many $ \mathcal{Y}_1, \dots, \mathcal{Y}_\N $, we obtain the form \eqref{eq:form of series}. Observe that cofiniteness conditions on $ W_0 $ are not needed to obtain the desired finite set $ S_1 $ and upper bound $ K_1 $ since $ w_0 $ is fixed. Finally, by linearity, this form remains true when $ w_0, \dots, w_{\N+1} $ are not assumed to be weight-homogeneous.
\end{proof}

The following theorem is the main result of this paper, and is proved with a method originally from \cite{HuangDEs}, and then adapted in \cite{TanMostlyUntwisted} to deal with monodromy around the origin. The reader should be familiar with Section~6 of \cite{TanMostlyUntwisted}. 

The main difficulty in applying Huang's method here is to find the correct generators $ \mathcal{A}_i $ that provide the ordinary differential equations \eqref{eq:DE for convergence c5} with a regular singularity at $ z_\N = 0 $ (possibly after a change of variables). After this problem is solved, the absolute convergence and analytic continuation of $ \langle w_0, \mathcal{Y}_1(w_1,z_1) \cdots \mathcal{Y}_\N(w_\N,z_\N) w_{\N+1} \rangle $ follow by induction on $ \N $ using the theory of regular singular differential equations (see Appendix B in \cite{Knapp+1986}, and Sections~5 and 6 of \cite{TanMostlyUntwisted}). The coefficients of \eqref{eq:DE for convergence c5} are, however, not necessarily single-valued in any open disk centered at $ z_\N = 0 $. The same trick of changing variables to $ \zeta_\N = z_\N^{1/t} $, for some appropriate $ t \in \mathbb{Z}_{>0} $, used in \cite{TanMostlyUntwisted} resolves this issue. In the present case, $ t = |G| $ will make all coefficients $ a_{k,N} (z_1, \dots, \zeta_\N^t) $ in \eqref{eq:DE for convergence c5} single-valued in $ \zeta_\N $, for all $ |z_1|, \dots, |z_{\N-1}| > |z_\N| $, since $ a_{k,N} (x_1, \dots, x_\N) \in R^{(x_\N = 0)} $.

\begin{theorem}[The convergence property for products of $ \N $ twisted intertwining operators among modules twisted by elements in a finite abelian group]  \label{thm: convergence c5} 
Let $ V $ be a CFT-type vertex operator algebra with a finite abelian subgroup $ G $ of $ \Aut V $. Let $ W_0, \dots, W_{\N+1} $, $ \widetilde{W}_0, \dots, \widetilde{W}_\N $, $ \mathcal{Y}_1, \dots, \mathcal{Y}_\N $, $ w_0, \dots, w_{\N+1} $ be as in Section~\ref{sec:setting c5}. Assume that either 
\begin{enumerate}
\item[(i)] $ n \geq 2 $ and $ g_0, \dots, g_{\N+1} \in G $ are general; or  
\item[(ii)] $ n = 1 $ and two of $ g_0, \dots, g_{\N+1} $ are $ g $ and $ g^{-1} $ for some $ g \in \Aut V $ of finite order, and the remaining automorphisms are $ \id_V $.
\end{enumerate} 
Then we have the following:
\begin{itemize}
\item[(1)] The multi-series 
\begin{equation*}
\langle w_0, \mathcal{Y}_1(w_1,z_1) \cdots \mathcal{Y}_\N(w_\N,z_\N) w_{\N+1} \rangle
\end{equation*}
in $ z_1, \dots, z_\N $ is absolutely convergent in the region $ |z_1| > \cdots > |z_\N| > 0 $.
\item[(2)] The sum of
\begin{equation*}
\langle w_0, \mathcal{Y}_1(w_1,z_1) \cdots \mathcal{Y}_\N(w_\N,z_\N) w_{\N+1} \rangle
\end{equation*}
in $ |z_1|> \dots >|z_\N|>0 $ can be analytically continued to a multivalued analytic function 
\begin{equation}
F(\langle w_0, \mathcal{Y}_1(w_1,z_1) \cdots \mathcal{Y}_\N(w_\N,z_\N) w_{\N+1} \rangle)
\end{equation}
defined on all of 
\begin{equation*}
 M^\N := \{ (z_1, \dots, z_\N) \in \widehat{\mathbb{C}}^\N : z_i \neq  0, \infty \text{ and } z_i \neq z_j \text{ when } i \neq j  \}.  
\end{equation*}
\end{itemize}
\end{theorem}

\begin{proof}
The proof is the same as that for Theorem~6.10 in \cite{TanMostlyUntwisted} now using the differential equation \eqref{eq:DE for convergence c5}, and replacing $ t = o(g) $ with $ t = |G| $. Lemma~\ref{lem:form of series} here replaces Lemma~6.7 in \cite{TanMostlyUntwisted}.
\end{proof}

Recall the notions of a component-isolated singularity and regular singularity from Definition~5.1 and Definition~5.2 of \cite{TanMostlyUntwisted}, respectively. (See also \cite{DuHuang} where the notion of a component-isolated singularity was introduced.)

\begin{proposition}\label{prop:sing special c5}
Let $ V $ be a CFT-type vertex operator algebra with a finite abelian subgroup $ G $ of $ \Aut V $. Let $ W_0, W_1, W_2 , W_{3} $, $ \widetilde{W}_0, \widetilde{W}_1, \widetilde{W}_2 $, $ \mathcal{Y}_1, \mathcal{Y}_2 $, $ w_0, w_1, w_2, w_{3} $ be as in Section~\ref{sec:setting c5}. Assume that either 
\begin{enumerate}
\item[(i)] $ n \geq 2 $ and $ g_0, \dots, g_{3} \in G $ are general; or  
\item[(ii)] $ n = 1 $ and two of $ g_0, \dots, g_{3} $ are $ g $ and $ g^{-1} $ for some $ g \in \Aut V $ of finite order, and the remaining automorphisms are $ \id_V $.
\end{enumerate}
Then every component-isolated singularity of $ F( \langle w_0, \mathcal{Y}_1(w_1, z_1) \mathcal{Y}_2(w_2, z_2) w_{3} \rangle) $ with a change of coordinates of the form $ (\zeta_1, \zeta_2) = (z_1, z_2)A $, for $ A \in \GL(2, \mathbb{C}) $, is a regular singularity. 
\end{proposition}

\begin{proof}
The proof is the same as that for Proposition~6.12 of \cite{TanMostlyUntwisted} after replacing $ t = o(g) $ with $ t = |G| $, using Theorem~\ref{thm: DEs c c5}, and checking that the differential equations have single-valued coefficients in $ (\eta_1, \eta_2) \in \D_{r_1}(0) \times \D_{r_2}(0) $ as explained in what follows.

Assume we have a change of variables $ (\zeta_1, \zeta_2) = (z_1,z_2) A $. Consider a component-isolated singularity at $ (\zeta_1, \zeta_2) = (\delta_1, \delta_2) $ for a function defined on $ M^2 $, for some $ \delta_1, \delta_2 \in \{0,\infty\} $. Recall the change of variables $ \eta_j^{s(\delta_j) t} = \zeta_j $ where $ s(0)= 1 $ and $ s(\infty) = -1 $.

The only part of the proof for Proposition~6.12 of \cite{TanMostlyUntwisted} that could be different is Case 2. So, we have to check that elements in $ R^{(\zeta_i = 0)}_{(\deg 0)} $ are still single-valued and defined for all $ (\eta_1, \eta_2) \in \D_{r_1}(0) \times \D_{r_2}(0) $ when $ \delta_i = 0 $ and $ \delta_j = \infty $.

Elements in $ R^{(\zeta_i = 0)}_{(\deg 0)} $ are linear combinations of monomials of the form $ z_1^{p/t} z_2^{q/t} (z_1 - z_2)^{r/t} $, where $ p, q,r \in \mathbb{Z} $ and $ p + q + r = 0 $. We rewrite these monomials as
\begin{equation*}
(\alpha_1 \eta_i^t + \beta_1 \eta_j^{-t})^{p/t} (\alpha_2 \eta_i^t + \beta_2 \eta_j^{-t})^{q/t} (\alpha_3 \eta_i^t + \beta_3 \eta_j^{-t})^{r/t} ,
\end{equation*}
for some $ \alpha_1,\alpha_2,\alpha_3, \beta_1, \beta_2, \beta_3 \in \mathbb{C} $. Since $ p + q + r = 0 $, we can rewrite this as
\begin{equation*}
(\alpha_1 \eta_i^t \eta_j^t + \beta_1)^{p/t} (\alpha_2 \eta_i^t \eta_j^t + \beta_2)^{q/t} (\alpha_3 \eta_i^t \eta_j^t + \beta_3)^{r/t}.
\end{equation*}
Observe that this is still single-valued in $ (\eta_1, \eta_2) $. If none of $ \beta_1, \beta_2, \beta_3 $ are equal to zero, then this expression is defined for $ \eta_i = 0 $ or $ \eta_j = 0 $. Otherwise, if some of $ \beta_1, \beta_2, \beta_3 $ are equal to zero, then the respective variables $ z_1, z_2 ,z_1 - z_2 $ are proportional to $ \zeta_i $. The sum of the powers of the variables $ z_1, z_2 ,z_1 - z_2 $ that are proportional to $ \zeta_i $ must be non-negative to ensure the monomial lies in $ R^{(\zeta_i = 0)} $. So, this expression is still defined for $ \eta_i = 0 $ or $ \eta_j = 0 $. 

Hence, we have checked that all elements in $ R^{(\zeta_i = 0)}_{(\deg 0)} $ are still single-valued and defined for all $ (\eta_1, \eta_2) \in \D_{r_1}(0) \times \D_{r_2}(0) $. Thus, we have shown the proposition.
\end{proof}

The previous proposition includes the singularities made from the combinations of $ z_1, z_2 = 0, \infty $ and $ z_1 - z_2 = 0 $ that provide component-isolated singularities. Note that this includes the important singularity $ (z_1 - z_2, z_2) = ( 0 , \infty) $ used for proving the associativity and operator product expansions of twisted intertwining operators. That is, when $ V $ is of CFT type, we have the \emph{convergence and extension property for products of two twisted intertwining operators} among $ C_1 $-cofinite discretely graded untwisted modules and $ C_2 $-cofinite discretely graded twisted modules twisted by elements in a finite abelian group. 

We now consider the special case where $ V $ is $ C_2 $-cofinite and of CFT type. 

\begin{lemma}\label{lem:c2 submodules c5}
Let $ V $ be a $ C_2 $-cofinite CFT-type vertex operator algebra with an automorphism $ g $, and let $ W $ be a $ g $-twisted $ V $-module. Let $ w \in W $ have homogeneous weight. Then the submodule $ \langle w \rangle $ generated by $ \{w\} $ is discretely graded and $ C_n $-cofinite for all $ n > 0 $.
\end{lemma}

\begin{proof}
Using Proposition~2.12 of \cite{TanCofinite}, we know that $ \langle w \rangle $ is spanned by elements of the form
\begin{equation} \label{eq:spanning element form}
b^{(i_1)}_{\alpha_{i_1} - n_1,0} \cdots b^{(i_k)}_{\alpha_{i_k} - n_k,0} w
\end{equation}
with $ k \in \mathbb{Z}_{\geq 0 }$, $ b^{(i_1)} , \dots, b^{(i_k)} \in B $ and $ n_1 , \dots,  n_k \in \mathbb{Z}_{> L} $ satisfying
\begin{enumerate}
\item[(i)] (ordering condition) \quad $ n_1 \geq \cdots \geq n_k $; and
\item[(ii)](repeat condition) \qquad $ n_j = n_{j+1} = \cdots = n_{j + M - 2} \leq K $ implies $ n_{j + M - 2} \neq n_{j + M-1} $ (i.e. $ n_j $ cannot repeat $ M $ times when $ n_j \leq K $).
\end{enumerate}

Recall that $ B = \{b^{(i)}\} $ is a finite set of homogeneous elements in $ V_+ $ that are $ \mathcal{S}_g $-eigenvectors with eigenvalue $ \alpha_i $ corresponding to $ b^{(i)} $, $ M \in \mathbb{Z}_{>0} $ is fixed and depends on $ V $, $ L \in \mathbb{Z} $ is fixed and depends on $ B $ and $ \langle w \rangle $, and $ K \in \mathbb{Z}_{> L} $ can be chosen to be sufficiently high.

We choose $ K $ to be a positive integer higher than $ L $. By the form of the spanning set, we see that there are only finitely many spanning elements \eqref{eq:spanning element form} with all $ n_j \leq K $. Then every other spanning element has factors $ b^{(i_j)}_{\alpha_{i_j} - n_j,0} $, for some $ n_j > K $, whose weight has real part $ \wt b^{(i_j)} - \Re{\alpha_{i_j}} + n_j - 1 \geq 1 $ which increase the real part of the weight of the spanning element. Hence, $ \langle w \rangle $ is discretely graded. From Theorem~2.14 of \cite{TanCofinite}, we know that $ \langle w \rangle $ is $ C_n $-cofinite for all $ n > 0 $.
\end{proof}
 
\begin{corollary} \label{cor:C2 convergence}
Let $ V $ be a $ C_2 $-cofinite CFT-type vertex operator algebra. Let $ G $ be a finite abelian subgroup of $ \Aut V $. Then twisted intertwining operators among lower bounded twisted $ V $-modules twisted by elements in $ G $ satisfy the convergence property for products of $ \N $ twisted intertwining operators, and the convergence and extension property for products of two twisted intertwining operators.
\end{corollary}

\begin{proof}
For each ordered $ (\N+2) $-tuple of weight-homogeneous vectors $ w_i \in W_i $, we consider the discretely graded and $ C_2 $-cofinite submodules $ \langle w_i \rangle $ by Lemma~\ref{lem:c2 submodules c5}. Let $ \widehat{W}_{\N} = \langle w_{\N+1} \rangle $. Let $ \widehat{W}_{\N-1} $ be the submodule of $ \W_{\N -1} $ spanned by all the coefficients of $ \mathcal{Y}_\N(\overline{w}_{\N},x) \overline{w}_{\N+1} $ for all $ \overline{w}_{\N} \in \langle w_\N \rangle $ and all $ \overline{w}_{\N+1} \in \widehat{W}_{\N} $. The twisted Nahm inequality of \cite{HuangNahmInequality} shows that $ \widehat{W}_{\N-1} $ is $ C_2 $-cofinite. We repeat this process recursively to obtain the remaining $ C_2 $-cofinite modules $ \widehat{W}_{i} $ for $ i = 0, \dots, \N - 2 $. Since $ \widehat{W}_{i} $ are submodules of lower bounded modules, they are also lower bounded. Hence, by Proposition~3.3 of \cite{HuangNahmInequality}, these $ \widehat{W}_{i} $ are also discretely graded (i.e., quasi-finite dimensional). 

We can restrict each twisted intertwining operator $ \mathcal{Y}_i $ to the discretely graded and $ C_2 $-cofinite modules $ \langle w_i \rangle $ and $ \widehat{W}_{i} $. Then Theorem~\ref{thm: convergence c5} and Proposition~\ref{prop:sing special c5} hold.
\end{proof}

\begin{remark}
Our results---Theorem~\ref{thm: convergence c5}, Proposition~\ref{prop:sing special c5}, and Corollary~\ref{cor:C2 convergence}---establish sufficient conditions for the convergence (and extension) property of certain chiral genus-zero correlation functions defined by products of twisted intertwining operators. These results contribute to the rigorous mathematical foundations of orbifold conformal field theory and provide certain conditions assumed in \cite{DuHuang} for the construction of tensor categorical structures arising from the twisted representation theory of vertex operator algebras.

We briefly compare our results with existing convergence results for products of twisted intertwining operators, and highlight their different limitations. For each $ g \in G $, every $ g $-twisted $ V $-module restricts to a module for the fixed-point subalgebra $ V^G $. In this context, the twisted intertwining operators become untwisted intertwining operators among $ V^G $-modules. If $ V $ is assumed to be simple, $ C_2 $-cofinite and of CFT type, and $ G $ is assumed to be finite abelian (or more generally solvable), then $ V^G $ is $ C_2 $-cofinite by results of Miyamoto in \cite{Miyamoto} (see also \cite{CM}). In this case, the results of \cite{HuangDEs} (and its logarithmic extensions \cite{HLZ7, HuangC1}) provide the conclusion of Corollary~\ref{cor:C2 convergence}. The corresponding categorical approach through $ V^G $, including the construction of braided $ G $-crossed tensor categories, is developed in \cite{McRaeTwisted} for a general finite group $ G $.

Observe that Corollary~\ref{cor:C2 convergence} does not require $ V $ to be simple. Since non-simple $ C_2 $-cofinite vertex operator algebras exist---for example, the $ \mathcal{W}_{2,p} $ algebras studied in \cite{AdamovicMilas}---our result is an improvement of the existing convergence result when $ V $ is assumed to be $ C_2 $-cofinite and $ G $ is assumed to be finite abelian, by extending to the case that $ V $ is non-simple.

Our approach has been to lift Huang's method \cite{HuangDEs} of deriving differential equations to a twisted setting. This approach requires a Jacobi identity for twisted intertwining operators, which is relatively straightforward to find in the possibly infinite abelian case, but presently unknown in the general case (not even the case of general finite groups). With the Jacobi identity \eqref{eq:twisted jacobi c5} that we have shown to be equivalent to the duality axiom for twisted intertwining operators among modules twisted by commuting elements, we have derived differential equations with controlled coefficients under certain Noetherian assumptions. Furthermore, our differential equations have regular singularities when $ G $ is finite. Consequently, limited to the theory of regular singular differential equations, our convergence result is restricted to the case that $ G $ is a finite abelian group. Our results hold for certain mixtures of $ C_1 $-cofinite and $ C_2 $-cofinite discretely graded modules.

On the other hand, when $ V $ is assumed to be simple, $ C_2 $-cofinite and of CFT type, Miyamoto's result further covers the case that $ G $ is solvable, which our result does not. At present, there is no result that guarantees that $ V^G $ is $ C_2 $-cofinite, when $ V $ is simple, $ C_2 $-cofinite and of CFT type, and $ G $ is a general finite group---such a result is likely difficult to prove, which presents a substantial obstacle to obtaining convergence results via fixed-point subalgebras.

Our work suggests a possible path to overcome this issue provided there is a successful development of an appropriate Jacobi-identity-type relation for intertwining operators among modules twisted by elements in a general finite group. We are interested in such a development. With the results of \cite{TanCofinite}, we know that $ C_2 $-cofinite CFT-type vertex operator algebras have $ C_2 $-cofinite twisted modules, without requiring any cofiniteness assumption on fixed point subalgebras. This suggests a promising approach to proving convergence when $ G $ is a general finite group.

Another key feature is that our method imposes cofiniteness conditions on the modules themselves, rather than on the vertex operator algebra. This perspective clarifies conceptual connections between cofiniteness conditions for twisted modules and the convergence of correlation functions in orbifold conformal field theory.
\end{remark}

\begin{remark}
An important step when using the theory of regular singular differential equations to prove the convergence of multi-series defined by matrix coefficients $ \langle w_0, \mathcal{Y}_1(w_1,z_1) \cdots \mathcal{Y}_\N(w_\N,z_\N) w_{\N+1} \rangle $ is to check that the multi-series already has form \eqref{eq:form of series}. Observe that the proof for Lemma~\ref{lem:form of series} will not supply the index set $ \frac{1}{t} \mathbb{Z} $ if the set of combined eigenvalues for $ \{ \mathcal{S}_g \}_{g \in G} $ is not a finite set of rational numbers. The proof for Lemma~\ref{lem:form of series} will supply upper bounds for the powers of the logarithms if the elements of $ G $ are semisimple, but will not extend beyond this case. (When $ g $ is not semisimple, $ Y_W(u,x) $ can contain logarithms, which means $ [L(0), Y_W(u,x)] = Y_W(L(0) u ,x) + x \frac{d}{dx} Y_W(u,x) $ will not imply $ [L(0), u_{m,k}] = (\wt u - m - 1) u_{m,k} $ in general. This may prevent the action of $ L(0) - L(0)_S $ on $ W $ from being nilpotent.) 

When proving convergence beyond the finite abelian case, it is interesting to observe that these two obstructions are consistent with the behavior of the coefficients of the differential equations \eqref{eq:DE with coefficients c5} that we have derived. In the first obstruction, a change of variables $ \zeta_\N = z_\N^{1/t} $ may not be available to turn the coefficients into functions of $ \zeta_\N $ that are holomorphic in sufficiently small open disks centered at $ 0 $. In the second obstruction, logarithms of $ z_\N $ in the coefficients may produce solutions without an upper bound for the powers of logarithms. For example, $ x\frac{dy}{dx} + (\log x) y = 0 $ has the solution $ y = a e^{\frac{-1}{2} (\log x)^2} $, which is a power series in $ \log x $. If we can constrain the coefficients even further, we can still obtain solutions in $ \mathbb{C}[\log x]\{x\} $. For example, $ x\frac{dy}{dx} + x(\log x) y = 0 $ has the solution $ y = a e^{x(1 - \log x)} \in  \mathbb{C}[\log x]\{x\} $. However, observe that this solution is still not in $ \mathbb{C}\{x\}[\log x] $.

An investigation into the behavior of these ``non-regular-singular" differential equations may give insight into a method that can be used to provide convergence, at least in certain cases, when $ G $ is not a finite group.
\end{remark}

\begin{remark}
If an analogue of Theorem~\ref{thm: convergence c5} can be proven when $ G $ is possibly infinite, there is an interesting subcase of Assumption (i), namely the case when all elements of $ G $ are locally unipotent. When this occurs, all elements of $ G $ will only have the eigenvalue $ 1 \in \mathbb{R}_{>0} $, so the modules considered are only required to be $ C_1 $-cofinite. Development of such a theorem would advance the understanding of $ G $-fixed-point subalgebras, their representation theory, and their extensions.

If $ \{Q_i\}_i $ is a collection of commuting locally nilpotent derivations of $ V $ that annihilate the conformal vector, then their exponentials $ e^{2 \pi \iu Q_i} $ would generate such a group. The fixed point algebra would then be equal to the intersection of the kernels of $ \{Q_i\}_i $. Such data can occur when constructing certain $ \mathcal{W} $-algebras, such as in \cite{AMlatticeconstruction} and \cite{ALMtripletvertexalgebra}. The connection between $ \mathcal{W} $-algebras, logarithmic conformal field theory, and orbifold conformal field theory for infinite groups is discussed in   \cite{Gentwistmod}.
\end{remark}

Observe that $ V $ is in fact not required to be of CFT type in Assumption (ii) for Theorem~\ref{thm: convergence c5} and Proposition~\ref{prop:sing special c5}. This is because the definition of $ C_1 $-cofiniteness already assumes $ u \in V_+ $ (see Definition~\ref{def:Cn-cofinite twisted modules}). Assumption (ii) is actually an extension of the special case considered in \cite{TanMostlyUntwisted} where $ W_{\N+1} $ was assumed to be $ g $-twisted and $ W_0 $ was assumed to be $ g^{-1} $-twisted, with the remaining in-state modules assumed to be untwisted. We refer to the case of Assumption (ii) as the \textit{mostly untwisted} case, since all but two of $ W_0, \dots, W_{\N+1} $ are untwisted.  

The convergence results of \cite{TanMostlyUntwisted} for $ C_1 $-cofinite modules in fact hint at similar convergence results in the mostly untwisted case. We can develop this intuition by referring back to the notion of a chiral genus-zero correlation function in (orbifold) conformal field theory. The matrix coefficient $ \langle w_0, \mathcal{Y}_1(w_1,z_1) \cdots \mathcal{Y}_\N(w_\N,z_\N) w_{\N+1} \rangle $ of a product of twisted intertwining operators should indeed converge and define a chiral genus-zero correlation function in the region $ |z_1| > \cdots > |z_{\N}| > 0 $ (recall the notion of radial ordering in conformal field theory). Since we are permitted to permute the field insertions on the Riemann sphere, we also expect for all permutations $ \sigma $ on $ \{0, \dots, \N+1\} $ that the multi-series
\begin{equation*}
 \langle w_{\sigma(0)}, \mathcal{Y}_{\sigma(1)}'(w_{\sigma(1)},z'_{\sigma(1)}) \cdots \mathcal{Y}_{\sigma(\N)}'(w_{\sigma(\N)},z'_{\sigma(\N)}) w_{{\sigma(\N+1)}} \rangle 
\end{equation*}
should converge in the region $ |z'_{\sigma(1)}| > \cdots > |z'_{\sigma(\N)}| > 0 $, for some possibly different twisted intertwining operators $ \mathcal{Y}_i' $ and new coordinates $ z'_i $. Furthermore, these two sums should be related by analytic continuation. For this reason, we should expect convergence in the mostly untwisted case for $ C_1 $-cofinite modules from the results of \cite{TanMostlyUntwisted}. 

We stress that we cannot use commutativity and associativity of twisted intertwining operators to prove convergence in all permutations of the mostly untwisted case; it is in fact the other way around. As explained in \cite{Reptheoryandorbifoldcft}, we must first prove convergence for products of $ \N $ intertwining operators, and the convergence and extension property for products of two intertwining operators. The next step in the mostly untwisted case is to apply the theory of \cite{DuHuang}. 

It is interesting to note that convergence in the mostly untwisted case does not require $ C_2 $-cofiniteness of $ V $; the only cofiniteness assumption required is $ C_1 $-cofiniteness of the in-state and intermediate modules. As far as we know, this is a new insight as all existing convergence results in orbifold conformal field theory that involve intertwining operators among twisted modules (other than \cite{TanMostlyUntwisted}) require $ V $ to be $ C_2 $-cofinite (or at least assume strong conditions on $ V^G $ that have not yet been verified when $ V $ is not $ C_2 $-cofinite). At the moment we do not have evidence of convergence in the case that all modules are $ C_1 $-cofinite and allowed to be twisted by arbitrary elements in $ G $---the best is $ C_2 $-cofiniteness. So we may not expect the category of all $ C_1 $-cofinite discretely graded $ g $-twisted $ V $-modules, for all $ g \in G $, to naturally admit the structure of a $ G $-crossed tensor category. However, since convergence and analytic extension in the mostly untwisted case among $ C_1 $-cofinite discretely graded modules have now been proven, these correlation functions provide evidence for the following conjecture about bimodule categories (see Chapter 7 of \cite{EGNO} for a definition). 

\begin{conjecture}\label{conj:bimodule categories}
Let $ V $ be a vertex operator algebra (not necessarily $ C_2 $-cofinite). Let $ g $ be an automorphism of $ V $ of finite order. Denote by $ V\textsf{-Mod}^{C_1} $ the tensor category of $ C_1 $-cofinite discretely graded (or lower bounded) $ V $-modules. Denote by $ V\textsf{-Mod}_g^{C_1} $ the category of $ C_1 $-cofinite discretely graded (or lower bounded) $ g $-twisted $ V $-modules. Then $ V\textsf{-Mod}_g^{C_1} $ is naturally a $( V\textsf{-Mod}^{C_1}, V\textsf{-Mod}^{C_1}) $-bimodule category, with fusion structure induced from the twisted intertwining operators in the mostly untwisted case. 
\end{conjecture}

If this conjecture is true, the twisted representation theory of vertex operator algebras can be applied to construct families of examples of module categories for tensor categories.

\section{Examples of $ C_n $-cofinite twisted modules} \label{sec:examples}

In this section, we provide examples of $ C_n $-cofinite twisted modules for CFT-type vertex operator algebras that can be used in the applications of Theorem~\ref{thm: DEs c5}, Theorem~\ref{thm: DEs c c5}, Theorem~\ref{thm: convergence c5}, Proposition~\ref{prop:sing special c5}, and Corollary~\ref{cor:C2 convergence}. We should remark that we do not need to know explicit examples of twisted intertwining operators to find benefit in using these results. These results verify hypotheses in the tensor-product constructions of [DH26], which do not require explicit formulas for the twisted intertwining operators. So, we are mainly interested in providing examples of CFT-type vertex operator algebras and their modules that satisfy conditions \ref{eq:cofinite conditions 1} and \ref{eq:cofinite conditions 2} as outlined in Section~\ref{sec:setting c5}.  
 
First we recall the following theorem from \cite{TanCofinite}, which is an extension of Buhl's classical result \cite{Buhl} lifted to the twisted setting. Recall that all $ g $-twisted $ V $-modules are weak $ g $-twisted $ V $-modules.

\begin{theorem}\label{thm:Cn cofiniteness}
Let $ V $ be a $ C_2 $-cofinite CFT-type vertex operator algebra with an automorphism $ g $. Let $ W $ be a finitely generated weak $ g $-twisted $ V $-module. Then $ W $ is $ C_n $-cofinite for all $ n \in \mathbb{Z}_{>0} $.
\end{theorem}

\begin{example}
Using Theorem~\ref{thm:Cn cofiniteness}, we have many examples of $ C_n $-cofinite twisted modules since there are many well-known examples of $ C_2 $-cofinite CFT-type vertex operator algebras. Examples include the simple affine vertex operator algebras at positive integral levels associated to WZW models \cite{FrenkelZhuAffine}, the simple Virasoro vertex operator algebras associated to minimal models \cite{FrenkelZhuAffine, Wang}, lattice vertex operator algebras associated to positive-definite even lattices \cite{FLM, DongLattice}, and the Moonshine Module $ V^\natural $ \cite{FLM}. These vertex operator algebras are explicitly seen to be of CFT type, and the proofs of their $ C_2 $-cofiniteness can be found in \cite{DLM_Modular}. Thus, by taking any of these $ C_2 $-cofinite CFT-type vertex operator algebras and any finite abelian group of its automorphisms, we can apply Corollary~\ref{cor:C2 convergence}.
\end{example}

\begin{example}
Corollary~\ref{cor:C2 convergence} does not require $ V $ to be simple, so it applies to non-simple $ C_2 $-cofinite CFT-type vertex operator algebras. For example, the vertex operator algebra $ \mathcal{W}_{2,3} $ was shown to be $ C_2 $-cofinite in \cite{AdamovicMilas} by Adamovi\'c and Milas. Furthermore, in Proposition~8.2, they showed that the maximal proper ideal of $ \mathcal{W}_{2,3} $ is generated by the conformal vector $ \omega $, and $ \mathcal{W}_{2,3}/\langle \omega\rangle \cong \mathbb{C}\one $. Hence, $ \mathcal{W}_{2,3} $ is an example of a non-simple $ C_2 $-cofinite vertex operator algebra, furthermore, with a one-dimensional simple quotient. The vertex operator algebra $ \mathcal{W}_{2,3} $ is defined as a subalgebra of the lattice vertex operator algebra $ V_D $, with $ D = \mathbb{Z}(2 \alpha) $ and $ (\alpha, \alpha) = 3 $, and a modified conformal vector $ \omega = \frac{1}{6}(\alpha(-1)^2 + \alpha(-2))\one $. So, $ L(0) $ acts as $ \frac{1}{6}(\alpha(0)^2 -\alpha(0)) $ on the lowest-weight vectors $ e^{2 m \alpha} $ giving them modified weight $ 6m^2 - m $, and $ L(0) $ assigns weight one to $ \alpha(-1) \one $ as usual. From this, we see that $ V_D $ has only non-negative weights, with weight-zero space equal to $ \mathbb{C}\one $. Since $ \mathcal{W}_{2,3} $ is a subalgebra of $ V_D $, it is of CFT type.

Given a lattice vertex operator algebra, recall that the zero modes of the weight-one elements in the Heisenberg subalgebra exponentiate to automorphisms. In this case, there is an involution $ \sigma = \exp\left( \frac{\pi \iu}{6} \alpha(0) \right) $ defined on $ V_D $ that restricts to
\begin{equation*}
\sigma(\one) = \one, \quad \sigma(\omega) = \omega, \quad \sigma(E) = -E, \quad \sigma(F) = -F, \quad \sigma(H) = H
\end{equation*}
on the generators $ \one, \omega , E, F, H $ of $ \mathcal{W}_{2,3} $ (here, we use the same notation as in Proposition~5.4 of \cite{AdamovicMilas}). Hence, $ \sigma $ restricts to an involution of $ \mathcal{W}_{2,3} $. We thus have an example of an abelian group $ \langle \sigma \rangle $ for our non-simple $ C_2 $-cofinite CFT-type vertex operator algebra to apply Corollary~\ref{cor:C2 convergence}. By Theorem~4.1 of \cite{HuangExistTwist}, there exists a non-zero lower-bounded $ \sigma $-twisted $ \mathcal{W}_{2,3} $-module. (We cannot use the existence theorem in \cite{DLM_Modular} since $ \mathcal{W}_{2,3} $ is not simple.)  So, convergence is not immediately obvious without Corollary~\ref{cor:C2 convergence}. Also note that  Corollary~\ref{cor:C2 convergence} proves convergence in cases that \cite{TanMostlyUntwisted} and the mostly untwisted case cannot, for example a correlation function among four $ \sigma $-twisted $ \mathcal{W}_{2,3} $-modules. 
\end{example}

One may ask whether there are examples of $ C_1 $-cofinite twisted modules for non-$ C_2 $-cofinite vertex operator algebras. We provide a family of examples now.

Consider a finite-dimensional simple or abelian Lie algebra $ \mathfrak{g} $ with a (possibly normalized) invariant non-degenerate symmetric bilinear form $ ( \cdot, \cdot ) $. If $ \mathfrak{g} $ is simple, let $ h^\vee $ be the dual Coxeter number of $ \mathfrak{g} $, and let $ ( \cdot, \cdot ) $ be the normalized form with square length $ 2 $ long roots.  If $ \mathfrak{g} $ is abelian, let $ h^\vee = 0 $. Every automorphism $ g $ of $ \mathfrak{g} $ that leaves $ (\cdot, \cdot) $ invariant  induces an automorphism of the universal affine vertex operator algebra $ V_{\hat{\mathfrak{g}}} (\lvl,0) $ at level $ \lvl \in \mathbb{C} \backslash \{ - h^\vee \} $ (or Heisenberg vertex operator algebra when $ \mathfrak{g} $ is abelian). Note that $ V_{\hat{\mathfrak{h}}} (\lvl,0) $ for $ \mathfrak{h} $ an abelian Lie algebra, i.e., a Heisenberg vertex operator algebra, is not $ C_2 $-cofinite. 

Let $ g $ be an automorphism of $ (\mathfrak{g},  ( \cdot, \cdot ) ) $. Denote by $ \mathcal{L}_{g} $, $ \mathcal{S}_{g} $ and $ \mathcal{N}_{g} $ the operators on $ \mathfrak{g} $ such that $ g = e^{2 \pi \iu \mathcal{L}_{g} }$, $ \mathcal{S}_{g} $ and $ \mathcal{N}_{g} $ are the commuting semisimple and nilpotent parts of $ \mathcal{L}_{g} $, and the eigenvalues of $ \mathcal{S}_{g} $ have real part in $ [0,1) $. The automorphism $ g $ naturally induces an automorphism of $ V_{\hat{\mathfrak{g}}} (\lvl,0) $, which we also denote by $ g $.

Recall from \cite{Huang_twisted_af_modules} (see also \cite{bakalovtwist}) that the \emph{twisted affine Lie algebra} $ \hat{\mathfrak{g}}^{[g]} $ is spanned by the elements 
\begin{equation}
\mathbf{k}, \ a \otimes t^{n} = a_n, \quad \text{for all } a \in \mathfrak{g} \text{ and } n \in \mathbb{C} \text{ such that } e^{2 \pi \iu \mathcal{S}_{g}} a = e^{2 \pi \iu n} a,
\end{equation} 
with Lie bracket given by 
\begin{align}
[a_m, b_n] &= [a,b]_{m + n} + ((m+\mathcal{N}_{g})a,b)\delta_{m+n,0} \mathbf{k} , \\
[\mathbf{k},a_m] &= 0.
\end{align}
We have the subspaces
\begin{gather*}
\hat{\mathfrak{g}}^{[g]}_+ := \Span\{a_n \in \hat{\mathfrak{g}}^{[g]} : \Re(n) > 0 \}, \quad \hat{\mathfrak{g}}^{[g]}_- := \Span\{a_n \in \hat{\mathfrak{g}}^{[g]} : \Re(n) < 0 \}, \\
\hat{\mathfrak{g}}^{[g]}_\mathbb{I} := \Span\{a_n \in \hat{\mathfrak{g}}^{[g]} : \Re(n) = 0 \}, \quad \hat{\mathfrak{g}}^{[g]}_0 := \hat{\mathfrak{g}}^{[g]}_\mathbb{I} \oplus \mathbb{C} \mathbf{k}.
\end{gather*}
Furthermore, $ \hat{\mathfrak{g}}^{[g]}_+ $, $ \hat{\mathfrak{g}}^{[g]}_- $ and $ \hat{\mathfrak{g}}^{[g]}_0 $ are subalgebras of $ \hat{\mathfrak{g}}^{[g]} $, giving a triangular decomposition
\begin{equation}
\hat{\mathfrak{g}}^{[g]} = \hat{\mathfrak{g}}^{[g]}_+ \oplus \hat{\mathfrak{g}}^{[g]}_0 \oplus \hat{\mathfrak{g}}^{[g]}_- .
\end{equation}

In \cite{Huang_twisted_af_modules}, certain $ g $-twisted $ V_{\hat{\mathfrak{g}}} (\lvl,0) $-modules called $ \widetilde{M}_\lvl^{[g]} $ were constructed and were shown to be equivalent to the induced $ \hat{\mathfrak{g}}^{[g]} $-modules 
\begin{equation}
U(\hat{\mathfrak{g}}^{[g]}) \otimes_{U(\hat{\mathfrak{g}}^{[g]}_+ \oplus \hat{\mathfrak{g}}^{[g]}_0)} M,
\end{equation}
where $ M $ is a finite-dimensional $ \hat{\mathfrak{g}}^{[g]}_0 $-module on which $ \mathbf{k} $ acts as $ \lvl $, and $ \hat{\mathfrak{g}}^{[g]}_+ $ acts trivially. Let $ a \in (V_{\hat{\mathfrak{g}}} (\lvl,0) )_{(1)} \cong \mathfrak{g} $ be such that $ \mathcal{S}_{g} a = \alpha a $ for some $ \alpha \in \mathbb{C} $ with $ \Re(\alpha) \in [0,1) $. To avoid confusion with subscripts, we denote the modes $ a_{n,k} $ by $ a(n,k) $ and also introduce $ a(n) = a(n,0) $ for brevity. Then the actions of the modes $ a(n) $ of the twisted vertex operator
\begin{equation}
(Y_{\widetilde{M}_\lvl^{[g]}})_0( a, x) = \sum_{n \in \alpha + \mathbb{Z}} a(n) x^{-n-1} , 
\end{equation} 
are given by the actions of $ a_n \in \hat{\mathfrak{g}}^{[g]} $ on $ \widetilde{M}_\lvl^{[g]} $. Using $ Y_{\widetilde{M}_\lvl^{[g]}}(a, x) = (Y_{\widetilde{M}_\lvl^{[g]}})_0( x^{-\mathcal{N}_{g}} a, x) $, we can determine the action of $ a(n,k) $ for $ k > 0 $. Since $ (V_{\hat{\mathfrak{g}}} (\lvl,0))_{(1)} \cong \mathfrak{g} $ generates $ V_{\hat{\mathfrak{g}}} (\lvl,0) $, this is sufficient to determine $ Y_{\widetilde{M}_\lvl^{[g]}}( \cdot , x) \cdot : V_{\hat{\mathfrak{g}}} (\lvl,0) \otimes \widetilde{M}_\lvl^{[g]} \to \widetilde{M}_\lvl^{[g]}\{x\} [\log x] $.

By the PBW-theorem, we know that $ \widetilde{M}_\lvl^{[g]} $ is spanned by elements of the form
\begin{equation} \label{eq:spanning set element 2}
a_1(\alpha_1 - n_1) \cdots a_k(\alpha_k - n_k) m , 
\end{equation}
for $ a_i \in (V_{\hat{\mathfrak{g}}} (\lvl,0))_{(1)}^{[\alpha_i]}$, $ n_1 \geq \cdots \geq n_k \in \mathbb{Z}_{> 0} $, $ k \in \mathbb{Z}_{\geq 0} $, and $ m \in M $.

Consider such an element for $ k > 0 $. If $ n_1 = 1 $, the element \eqref{eq:spanning set element 2} is immediately in $ C_1(\widetilde{M}_\lvl^{[g]}) $. Otherwise, $ n_1 > 1 $ and we can use the $ L(-1) $-derivative property (see Lemma~1.6 of \cite{TanCofinite} for a precise statement) to show that the element \eqref{eq:spanning set element 2} is in $ C_1(\widetilde{M}_\lvl^{[g]}) $. Hence, $ \widetilde{M}_\lvl^{[g]} = C_1(\widetilde{M}_\lvl^{[g]}) + M $ and $ \dim \widetilde{M}_\lvl^{[g]}/C_1(\widetilde{M}_\lvl^{[g]}) \leq \dim M < \infty $. Thus, $ \widetilde{M}_\lvl^{[g]} $ is a $ C_1 $-cofinite (weak) twisted module. 

Furthermore, $ \widetilde{M}_\lvl^{[g]} $ has the structure of a twisted module for $ V_{\hat{\mathfrak{g}}} (\lvl,0) $ as a vertex \emph{operator} algebra. Namely, it has an $ L(0) $-action given by 
\begin{equation}
\begin{aligned}
L(0) &= \frac{1}{2(\lvl + h^\vee)} \sum_{i \in I} \Big( \sum_{p \in \alpha^i + \mathbb{Z}_{\geq 0}} a^{i'}(-p) a^i(p)   + \sum_{p \in \alpha^i + \mathbb{Z}_{< 0}}  a^{i}(p) a^{i'}(-p) \\
&\quad -  [(\mathcal{N}_{g} - \alpha^i + 1) a^{i'}, a^i ] (0) - \frac{\lvl}{2} (( \mathcal{N}_{g}-\alpha^i) (\mathcal{N}_{g} - \alpha^i + 1) a^{i'} , a^i ) \Big),
\end{aligned}
\end{equation}
for some basis $ \{a^i \}_{i \in I} $ of $ \mathfrak{g} $ consisting of generalized eigenvectors for $ g $ with $ \mathcal{S}_g a^i = \alpha^i a^i $, $ \Re(\alpha^i) \in [0,1) $, and a dual basis $ \{ a^{i'} \}_{i \in I} $ (see Proposition~4.8 of \cite{Huang_twisted_af_modules} for details; note that we have rewritten the right-hand side to make the modes normally ordered). From the ordering of the modes, we can immediately see that the $ L(0) $-action preserves $ M $ as a subspace of $ \widetilde{M}_\lvl^{[g]} $. Since $ M $ is finite-dimensional, it decomposes into generalized eigenspaces for $ L(0) $. We also know that the action of $ a(n) $ shifts the generalized $ L(0) $-eigenvalue by $ - n $. So, from \eqref{eq:spanning set element 2}, we can see that $ \widetilde{M}_\lvl^{[g]} $ decomposes into generalized  $ L(0) $-eigenspaces. Since $ M $ is finite-dimensional,  $ \mathfrak{g} $ is finite-dimensional, and $ \Re(\wt a_i(\alpha_i - n_i)) = \wt a_i - \Re(\alpha_i) + n_i - 1 > 0 $ for all $ n_i \in \mathbb{Z}_{>0} $, we see that $ \widetilde{M}_\lvl^{[g]} $ is furthermore discretely graded. Thus, we have the following result.

\begin{proposition} \label{prop:C1 cofinite affine modules}
Let $ \mathfrak{g} $ be a finite-dimensional simple or abelian Lie algebra with an invariant non-degenerate symmetric bilinear form $ (\cdot, \cdot) $ (normalized with square length $ 2 $ long roots if $ \mathfrak{g} $ is simple). Let $ g $ be an automorphism of $ ( \mathfrak{g},(\cdot, \cdot)) $. Let $ M $ be a finite-dimensional $ \hat{\mathfrak{g}}^{[g]}_0 $-module on which $ \mathbf{k} $ acts as $ \lvl $. Let $ \lvl \in \mathbb{C}\backslash \{-h^\vee\} $. Then the $ g $-twisted $ V_{\hat{\mathfrak{g}}} (\lvl,0) $-module $ \widetilde{M}_\lvl^{[g]} $ is $ C_1 $-cofinite and discretely graded. \qed
\end{proposition}

\begin{example}
Take any finite-dimensional simple or abelian Lie algebra $ \mathfrak{g} $ with an invariant non-degenerate symmetric bilinear form $ (\cdot, \cdot) $ (normalized with square length $ 2 $ long roots if $ \mathfrak{g} $ is simple), and an  automorphism $ g $ of $(\mathfrak{g},(\cdot, \cdot)) $ of finite order. Choose modules $ W_0, \dots, W_{\N+1} $ from Proposition~\ref{prop:C1 cofinite affine modules} such that one is $ g $-twisted, one is $ g^{-1} $-twisted and the remaining are untwisted. Choose compatible intermediate modules that are also $ C_1 $-cofinite and discretely graded. Then the mostly untwisted case (i.e., Assumption (ii)) of Theorem~\ref{thm: convergence c5} applies, as does Proposition~\ref{prop:sing special c5} when $ \N = 2 $. We suspect that this application can be used to construct families of bimodule categories as in Conjecture~\ref{conj:bimodule categories} above.

Even when the automorphism $ g $ is not of finite order, we can still apply Theorems~\ref{thm: DEs c5} and~\ref{thm: DEs c c5} in certain cases. For example, consider $ \mathfrak{g} = \mathfrak{sl}(2,\mathbb{C}) $ with the standard generators $ E $, $ F $ and $ H $, and let $ (\cdot, \cdot) $ be the standard normalized Killing form. The derivation $ \ad(E) $ is nilpotent, so the automorphism $ e^{2 \pi \iu \ad(E)} $ is unipotent. Now consider the abelian group $ G $ generated by $ e^{2 \pi \iu \ad(E)} $ lifted to an automorphism of $ V_{\widehat{\mathfrak{sl}(2,\mathbb{C})}} (\lvl,0) $. Since every element $ g \in G $ has $ 1 $ as its only eigenvalue, it suffices to assume that the $ g $-twisted modules are $ C_1 $-cofinite and discretely graded to satisfy conditions \ref{eq:cofinite conditions 1} and \ref{eq:cofinite conditions 2} in the setting of this paper. Furthermore, the ring $ R $ is equal to $ \mathbb{C}[(x_i - x_j)^{\pm 1} , \log (x_i - x_j) : i < j \in  \{1, \dots, \N +1\} ] $. So, $ R $ and $ R^{(c=0)} $ are explicitly Noetherian. Thus, we can apply Theorems~\ref{thm: DEs c5} and \ref{thm: DEs c c5} to all $ g $-twisted modules from Proposition~\ref{prop:C1 cofinite affine modules}, for all $ g \in G = \langle e^{2 \pi \iu \ad(E)} \rangle $.

This result is easily generalizable to other finite-dimensional simple algebras, and their nilpotent derivations that respect the bilinear form. 
\end{example} 

\begin{example}
For the sake of concreteness, we will provide an explicit well-known example of a twisted module (see \cite{LepowskyWilson, FLM} and \cite{DFMS}) and show that it is $ C_1 $-cofinite. Consider the Heisenberg vertex operator algebra (of level $ \lvl = 1 $) corresponding to the one-dimensional abelian Lie algebra $ \mathfrak{h} = \mathbb{C} a $ with bilinear form given by $ (a,a) = 1 $. This is the vertex operator algebra corresponding to the theory of a bosonic string propagating on the one-dimensional non-compact space-time $ \mathcal{M} = \mathbb{R} $. The symmetry $ r \mapsto -r $ of $ \mathcal{M} $ induces a symmetry $ g $ of the Heisenberg vertex operator algebra characterized by $ a \mapsto - a $. Then the weight-one subspace $ (V_{\hat{\mathfrak{h}}} (\lvl,0))_{(1)} $ decomposes into eigenspaces for $ \mathcal{S}_{g} $; in this case, the entire subspace corresponds to the eigenvalue $ \alpha = \frac{1}{2} $ since $ g a = - a = e^{2 \pi \iu \frac{1}{2}} a $. Then $ \hat{\mathfrak{h}}^{[g]} = \Span \{ a_{n} , \mathbf{k} : n \in \frac{1}{2} + \mathbb{Z} \} $  and the subalgebra $ \hat{\mathfrak{h}}_{\mathbb{I}}^{[g]} $ is trivial. So, we can consider the one-dimensional $ \hat{\mathfrak{h}}_{0}^{[g]} $-module $ M = \mathbb{C}m $ spanned by $ m $, with $ \mathbf{k} $ acting as the identity. The corresponding induced module
\begin{equation*}
 \widetilde{M}_\lvl^{[g]} = U(\hat{\mathfrak{h}}^{[g]}) \otimes_{U(\hat{\mathfrak{h}}^{[g]}_+ \oplus \hat{\mathfrak{h}}^{[g]}_0)} M
\end{equation*}
has a basis consisting of the elements 
\begin{equation*} 
a(\tfrac{1}{2} - n_1) \cdots a(\tfrac{1}{2}  - n_k) m , \qquad n_1 \geq \cdots \geq n_k \in \mathbb{Z}_{> 0}, \ k \in \mathbb{Z}_{\geq 0} .
\end{equation*}
From this basis, we can see that $\widetilde{M}_\lvl^{[g]} $ is $ C_1 $-cofinite. Furthermore, 
\begin{equation*}
L(0) m = \frac{1}{2}\left(0 + 0 + 0 - \frac{1}{2}\left( -\frac{1}{2} \right)\left(\frac{1}{2} \right) \right)m = \frac{1}{16} m.
\end{equation*}
So, $ \widetilde{M}_\lvl^{[g]} $ is discretely graded into $ L(0) $-eigenspaces with eigenvalues in $ \frac{1}{16} + \frac{1}{2} \mathbb{Z}_{\geq 0} $. The dimension of $ (\widetilde{M}_\lvl^{[g]})_{(\frac{1}{16} + n)} $ is equal to the number of possible ways to partition $ n $ with numbers in $ \frac{1}{2} + \mathbb{Z}_{\geq 0} $.

Now define two of $ W_0, \dots, W_{\N+1} $ to be the $ g $-twisted module $ \widetilde{M}_\lvl^{[g]} $, and define the remaining modules to be the usual irreducible untwisted modules generated by a lowest-weight vector with $ a(0) $-eigenvalue $ \lambda \in \mathbb{C} $. For any twisted intertwining operators $ \mathcal{Y}_i(\cdot, x) \cdot: W_i \otimes \W_i \to \W_{i-1}\{x\}[\log x] $ with compatible $ C_1 $-cofinite and discretely graded modules $ \W_i $, and vectors $ w_i \in W_i $, the multi-series
\begin{equation*}
\langle w_0, \mathcal{Y}_1(w_1,z_1) \cdots \mathcal{Y}_\N(w_\N,z_\N) w_{\N+1} \rangle
\end{equation*}
converges absolutely in the region $ |z_1| > \cdots > | z_\N| > 0 $ and analytically extends to all of $ M^\N $. Furthermore, when $ \N = 2 $, the component-isolated singularities of $ F( \langle w_0, \mathcal{Y}_1(w_1, z_1) \mathcal{Y}_2(w_2, z_2) w_{3} \rangle) $, corresponding to a change of coordinates of the form $ (\zeta_1, \zeta_2) = (z_1, z_2)A $ for $ A \in \GL(2, \mathbb{C}) $, are regular singularities. 
\end{example}

The arguments used to show Proposition~\ref{prop:C1 cofinite affine modules} do not work to show that the twisted modules are $ C_n $-cofinite for $ n \geq 2 $. So, one may ask whether there are $ C_n $-cofinite twisted modules for non-$ C_2 $-cofinite vertex operator algebras for some $ n \geq 2 $. We show that this is also possible.

We now restrict to the case that $ \mathfrak{g} $ is simple and the level $ \lvl $ is a positive integer. In this case, the simple affine vertex operator algebras $ L_{\hat{\mathfrak{g}}} (\lvl,0) $ are $ C_2 $-cofinite. There is a surjective vertex operator algebra homomorphism $ \pi: V_{\hat{\mathfrak{g}}} (\lvl,0)  \to L_{\hat{\mathfrak{g}}} (\lvl,0) $ with kernel equal to the unique maximal proper ideal of $ V_{\hat{\mathfrak{g}}} (\lvl,0) $. 

Let $ g $ be an automorphism of $ (\mathfrak{g}, (\cdot,\cdot)) $ and consider the induced automorphism $ g $ on $ V_{\hat{\mathfrak{g}}} (\lvl,0) $. We claim that this induces an automorphism of $ L_{\hat{\mathfrak{g}}} (\lvl,0) $, also denoted by $ g $, according to 
\begin{equation} \label{eq:pi g = g pi}
\begin{tikzcd}
V_{\hat{\mathfrak{g}}} (\lvl,0) \ar[r,"\pi"] \ar[d,"g",swap] &L_{\hat{\mathfrak{g}}} (\lvl,0) \ar[d,"g"]\\
V_{\hat{\mathfrak{g}}} (\lvl,0) \ar[r,"\pi"] &L_{\hat{\mathfrak{g}}} (\lvl,0) .
\end{tikzcd}
\end{equation}

The subspace $ g \ker(\pi) $ is a proper ideal of $ V_{\hat{\mathfrak{g}}} (\lvl,0) $, so is contained in $ \ker (\pi ) $. So, $ \ker( \pi ) \subseteq \ker(\pi \circ g) $. By the universal property of the quotient, we obtain a vertex operator algebra homomorphism $ g : L_{\hat{\mathfrak{g}}} (\lvl,0) \to L_{\hat{\mathfrak{g}}} (\lvl,0) $. Since $ g $ is invertible, we can repeat the same process for $ g^{-1} $ and see that 
\begin{align*}
\id_{L_{\hat{\mathfrak{g}}(\lvl,0)}} \circ \pi &= \pi \circ \id_{V_{\hat{\mathfrak{g}}(\lvl,0)}} = \pi \circ g \circ g^{-1} = g \circ g^{-1} \circ \pi, \\
\id_{L_{\hat{\mathfrak{g}}(\lvl,0)}} \circ \pi &= \pi \circ \id_{V_{\hat{\mathfrak{g}}(\lvl,0)}} = \pi \circ g^{-1} \circ g = g^{-1} \circ g \circ \pi.
\end{align*}  
So, $ g \circ g^{-1} = \id_{L_{\hat{\mathfrak{g}}(\lvl,0) }} = g^{-1} \circ g $ by the universal property of the quotient. Hence, $ g : L_{\hat{\mathfrak{g}}} (\lvl,0) \to L_{\hat{\mathfrak{g}}} (\lvl,0) $ is indeed an automorphism.

We claim that every $ g $-twisted $ L_{\hat{\mathfrak{g}}} (\lvl,0) $-module is naturally a $ g $-twisted $ V_{\hat{\mathfrak{g}}} (\lvl,0) $-module. Let $ (W, Y_W) $ be a $ g $-twisted $ L_{\hat{\mathfrak{g}}} (\lvl,0) $-module. Define $ \widetilde{W} = W $ and
\begin{equation}
Y_{\widetilde{W}} = Y_W \circ \pi : V_{\hat{\mathfrak{g}}} (\lvl,0) \to (\End W) \{x\} [\log x] .
\end{equation}
Each condition in the definition of a $ g $-twisted module (see Definition~3.1 of \cite{HUANG2018346}, but omit any need for $ g $-actions) can be checked for $ (\widetilde{W}, Y_{\widetilde{W}}) $. For example, we use the fact that $ \pi $ is a homomorphism to obtain
\begin{equation*}
Y_{\widetilde{W}}(Y(u,x_0)v,x_2) = Y_{W}(\pi Y(u,x_0)v,x_2) = Y_{W}( Y( \pi u,x_0) \pi v,x_2)
\end{equation*}
while showing the duality property, and we use \eqref{eq:pi g = g pi} to show the equivariance property
\begin{equation*}
e^{2 \pi \iu x\frac{d}{dx}} Y_{\widetilde{W}}(gu,x) = e^{2 \pi \iu x \frac{d}{dx}} Y_{W}(\pi g u,x) = e^{2 \pi \iu x\frac{d}{dx}} Y_{W}(g \pi u,x) = Y_{W}(\pi u,x) = Y_{\widetilde{W}}(u,x).
\end{equation*}
The $ L(0) $-grading, and any grading restrictions such as discrete grading, follow from the fact that $ \pi $ preserves conformal vectors.

By Theorem~\ref{thm:Cn cofiniteness}, we know that finitely generated $ g $-twisted $ L_{\hat{\mathfrak{g}}} (\lvl,0)  $-modules are $ C_n $-cofinite for all $ n > 0 $. Moreover, when these modules have an $ L(0) $-grading, they are finitely generated by generalized $ L(0) $-eigenvectors and hence are discretely graded by arguments similar to those used in the proof of Lemma~\ref{lem:c2 submodules c5}.

We claim that a $ C_n $-cofinite $ g $-twisted $ L_{\hat{\mathfrak{g}}} (\lvl,0)  $-module $ W $ is also $ C_n $-cofinite when viewed as a $ g $-twisted $ V_{\hat{\mathfrak{g}}} (\lvl,0)  $-module. Let $ u \in L_{\hat{\mathfrak{g}}} (\lvl,0) $. Then there exists $ v \in V_{\hat{\mathfrak{g}}} (\lvl,0) $ such that $ \pi v = u $ and so
\begin{align*}
\res_x x^{-n} Y_W(x^{\mathcal{L}_{g}} u, x) &= \res_x x^{-n} Y_W(x^{\mathcal{L}_{g}} \pi v, x) \\
&= \res_x x^{-n} Y_W( \pi x^{\mathcal{L}_{g}} v, x) = \res_x x^{-n} Y_{\widetilde{W}}( x^{\mathcal{L}_{g}} v, x). 
\end{align*}
Hence, $ C_n(W) \subseteq C_n(\widetilde{W}) $. Therefore, 
\begin{align*}
\dim \widetilde{W}/C_n(\widetilde{W}) \leq \dim W/ C_n(W) < \infty.
\end{align*}
Thus, we have the following examples of $ C_n $-cofinite twisted modules for non-$ C_2 $-cofinite vertex operator algebras.

\begin{proposition}\label{prop:C_2-cofinite for non-C_2}
Let $ \mathfrak{g} $ be a finite-dimensional simple Lie algebra with a normalized invariant non-degenerate symmetric bilinear form $ (\cdot, \cdot) $. Let $ g $ be an automorphism of $ ( \mathfrak{g},(\cdot, \cdot)) $. Let $ \lvl \in \mathbb{Z}_{>0} $. Then the finitely generated $ g $-twisted $ L_{\hat{\mathfrak{g}}} (\lvl,0) $-modules are $ C_n $-cofinite and discretely graded $ g $-twisted $ V_{\hat{\mathfrak{g}}} (\lvl,0) $-modules for all $ n \in \mathbb{Z}_{>0} $. \qed
\end{proposition}

\begin{example}
Consider a finite-dimensional simple Lie algebra and a normalized invariant non-degenerate symmetric bilinear form pair $ ( \mathfrak{g},(\cdot, \cdot)) $. Let $ G $ be a finitely generated abelian group of automorphisms of $ ( \mathfrak{g},(\cdot, \cdot)) $, and lift $ G $ to a group of automorphisms of $ V_{\hat{\mathfrak{g}}} (\lvl,0) $, for some level $ \lvl \in \mathbb{Z}_{>0} $. With Propositions~\ref{prop:C1 cofinite affine modules} and \ref{prop:C_2-cofinite for non-C_2}, we have a collection of $ C_1 $-cofinite or $ C_2 $-cofinite $ g $-twisted $ V_{\hat{\mathfrak{g}}} (\lvl,0) $-modules, for each $ g \in G $ depending on whether all eigenvalues of $ g $ are in $ \mathbb{R}_{>0} $. Since $ V_{\hat{\mathfrak{g}}} (\lvl,0) $ and $ G $ are finitely generated, the ring $ R $ is Noetherian and we can apply Theorem~\ref{thm: DEs c5}. In specific cases, we can also apply Theorem~\ref{thm: DEs c c5}. 

When $ G $ is finite, we now have a collection of $ C_1 $-cofinite untwisted $ V_{\hat{\mathfrak{g}}} (\lvl,0) $-modules  and $ C_2 $-cofinite $ g $-twisted $ V_{\hat{\mathfrak{g}}} (\lvl,0) $-modules for each non-identity $ g \in G $. With this data, we can apply Theorem~\ref{thm: convergence c5} and, when $ \N = 2 $, Proposition~\ref{prop:sing special c5}.
\end{example}

\printbibliography

\noindent {\small \sc Yau Mathematical Sciences Center, Tsinghua University, Beijing 100084, China}

\noindent {\em E-mail address}: \url{danieltan@mail.tsinghua.edu.cn}

\end{document}